\documentclass[11pt, a4paper]{article}
\usepackage{amsmath}
\usepackage{amsfonts}
\usepackage{amssymb}
\usepackage{ulem}
\usepackage{amsthm}
\usepackage{url}
\usepackage{dcolumn}
\usepackage{cite}
\usepackage{fancyhdr}

\usepackage{enumitem}

\usepackage{mathrsfs}

\usepackage{graphicx, color}
\usepackage{hyperref}
\usepackage{xcolor}

\usepackage{secdot}
\sectiondot{subsection}

\usepackage{multirow}

\newtheorem*{rem*}{Remark}

\newtheorem*{ques}{Question}

\newtheorem{nameth}{Sugimoto's Perspective}
\newtheorem*{nameth2}{Amfirifma's Perspective}

\renewcommand{\thefigure}{\arabic{section}{.}\arabic{figure}}
\renewcommand{\thetable}{\arabic{section}{.}\arabic{table}}
\makeatletter
\@addtoreset{figure}{section}
\@addtoreset{table}{section}
\makeatother

\usepackage{here}

\begin{document}

\title{Tile sets consisting of two types of concave polygons derived from periodic tilings corresponding to non-periodic tilings with hat and turtle tiles}
\author{ Teruhisa SUGIMOTO$^{ 1), 2)}$ }
\date{}
\maketitle

{\footnotesize
\begin{center}
$^{1)}$ The Interdisciplinary Institute of Science, Technology and Art

$^{2)}$ Japan Tessellation Design Association

E-mail: ismsugi@gmail.com
\end{center}
}

{\small
\begin{abstract}
\noindent
Using a convex pentagonal monotile belonging to the Type 5 family, we investigate the 
relationships among the hat tile, turtle tile, and \mbox{Tile$(1, 1)$}. By applying 
Sugimoto's Perspective and Amfirifma's Perspective, we obtain four types of concave 
polygons, AH-tile, BH-tile, AT-tile, and BT-tile, each having Heesch number 1 under the 
conditions considered. We show that these polygons correspond to clusters used to 
generate the non-periodic tilings $\mathscr{T}_h$ and $\mathscr{T}_s$. We further 
discuss the possibility that tile sets consisting of pairs selected from these polygons may 
correspond to \mbox{$\textit{ASPmr}\{\text{A-tile}, \text{B-tile}\}$}.
\end{abstract}
}


\section{Introduction}
\label{section1}

A \textit{tiling} (or \textit{tessellation}) of the plane is a collection of regions (diagrams), called 
tiles that covers the plane without gaps or overlaps, except for the boundaries of the tiles. 
The term ``tile'' refers to a bounded topological disk, whose boundary is a simple closed 
curve \cite{ref1, ref3}. 

In general, a set containing a finite number of tile types (i.e., a set consisting of finitely many tile 
types) is considered when studying tilings. Each tile in such a set is called a \textit{prototile}. 
If tilings of the plane can be formed using these prototiles, then the \textit{set of prototiles} 
(simply referred to as a ``set of tiles'' or ``tile set'' in this manuscript) is said to admit tilings 
of the plane \cite{ref1, ref2, ref3}. Note that, when a tile set capable of generating tilings is given, 
the tile set is not required to admit a tiling that uses all of the tiles in the set (that is, there 
is no requirement that a tiling using all of the tiles exist), nor is each tile required to admit 
tilings on its own. A tile set is said to admit tilings if a tiling can be generated using some of 
the tiles in the set. In this study, however, we impose the following additional condition on a 
tile set: every tile in the set must appear in at least one tiling generated by the tile set.

A tiling is said to be monohedral if all tiles have the same size and shape (i.e., are congruent). 
Reflected tiles (i.e., posterior-side tiles) are allowed in monohedral tilings. In other words, 
in monohedral tilings, anterior-side and posterior-side tiles are treated as the same 
tile. \cite{ref1, ref3, ref5, ref10}. Therefore, in this case, the tiling is regarded as having only 
one tile type, and the corresponding set of prototiles consists of a single prototile. A tile 
that admits a monohedral tiling is called a monotile.

\begin{figure}[t]
 \centering\includegraphics[width=14cm,clip]{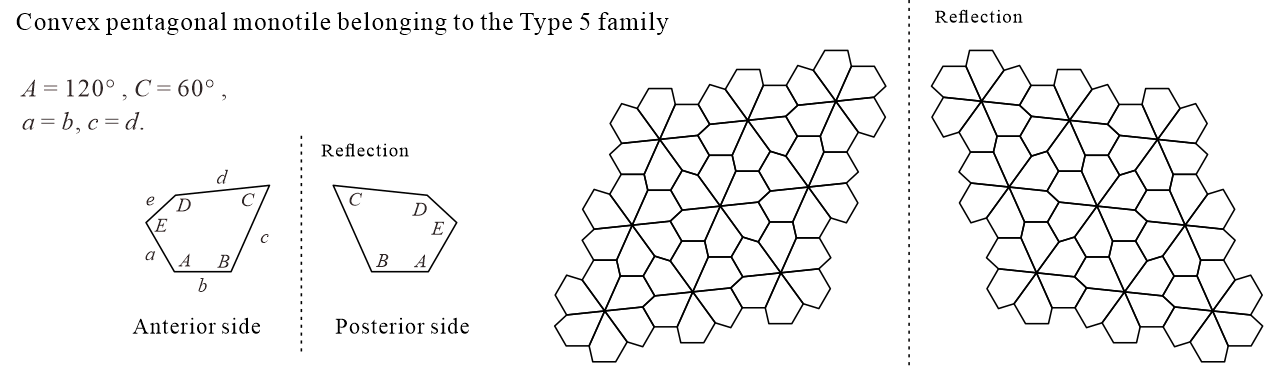} 
  \caption{{\small Convex pentagonal monotiles belonging to the Type 5 family and representative tilings.} }
\label{Fig.1-1}
\end{figure}

Convex pentagonal monotiles belonging to the Type 5 family (also referred to as convex 
pentagonal monotiles belonging to Type 5, or simply as convex pentagonal monotiles of Type 5) 
are convex pentagons that satisfy the tile conditions shown in Figure~\ref{Fig.1-1} and admit tilings. 
A convex pentagonal monotile belonging to the Type 5 family can generate a periodic tiling\footnote{ 
A tiling is said to have translational symmetry if its translation by a non-zero vector coincides 
with itself. In general, a tiling of the plane is said to be \textit{periodic} (i.e., a periodic tiling) if 
it coincides with itself under translations by two linearly independent vectors. However, 
throughout this study, a tiling with translational symmetry is referred to as periodic, and a 
tiling without translational symmetry is referred to as \textit{non-periodic}.
} such as that shown in Figure~\ref{Fig.1-1}, and this tiling is referred to as a representative 
tiling of a convex pentagonal monotile belonging to the Type 5 family. A representative tiling is 
a tiling formed solely from geometric relationships derived from the tile conditions for each 
Type \cite{ref3, ref10}. Hereafter, the representative tiling of a convex pentagonal monotile 
belonging to the Type 5 family is simply referred to as the Type 5 tiling.  
Note that the Type 5 tiling shown in Figure~\ref{Fig.1-1} is formed using only one side of the 
convex pentagonal monotile.

\renewcommand{\figurename}{{\small Figure}}

\begin{figure}[H]
 \centering\includegraphics[width=12cm,clip]{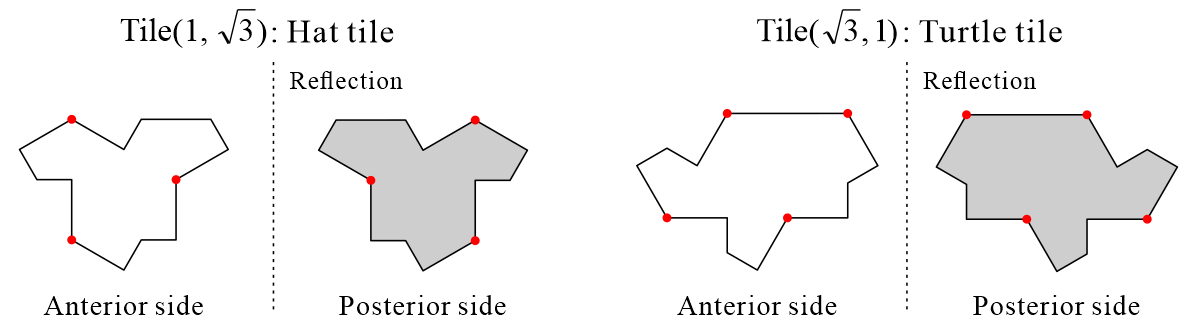} 
  \caption{{\small Hat tile and turtle tile.} }
\label{Fig.1-2}
\end{figure}

\begin{figure}[H]
 \centering\includegraphics[width=12.5cm,clip]{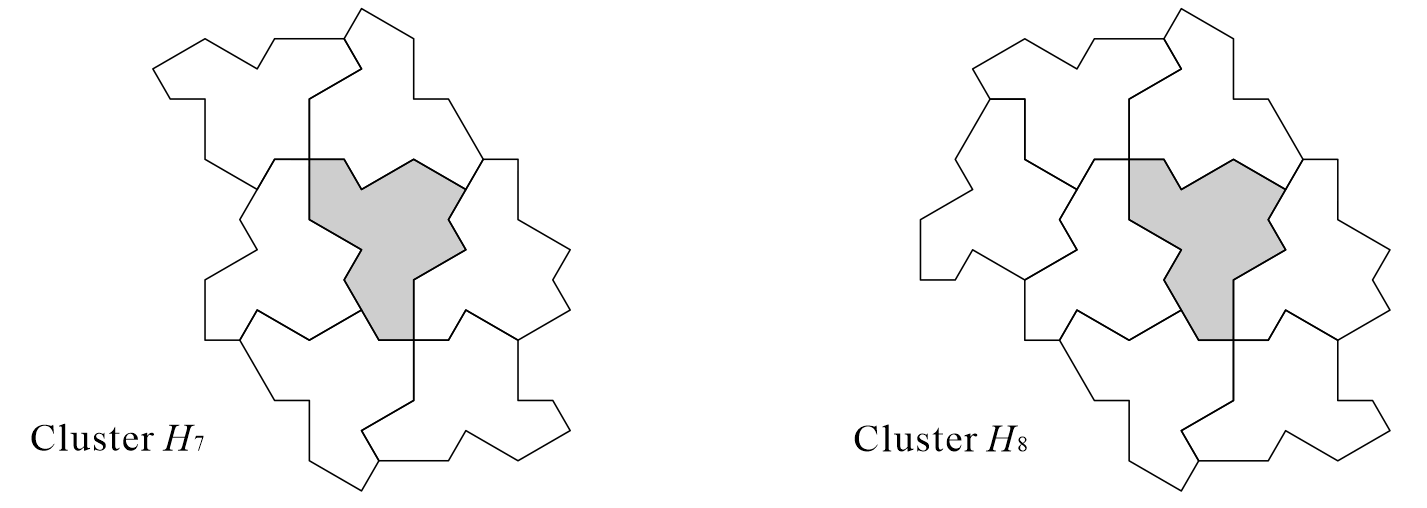} 
  \caption{{\small Clusters $H_7$ and $H_8$ formed by hat tiles.} }
\label{Fig.1-3}
\end{figure}

As shown in Figure~\ref{Fig.1-2}, the hat tile and the turtle tile are aperiodic monotiles\footnote{ 
A set of prototiles is said to be \textit{aperiodic} (i.e., an ``aperiodic set of prototiles,'' 
``aperiodic set of tiles,'' or ``aperiodic tile set'') if it admits tilings of the plane and all tilings 
generated from those prototiles are non-periodic. If there exists an aperiodic tile set consisting 
of a single prototile that admits only non-periodic tilings, then the prototile is called an 
aperiodic monotile \cite{ref1, ref2, ref3, ref5, ref6, ref13, ref14}.  \label{footnote2}
} corresponding to \mbox{Tile$(1, \sqrt{3})$} and \mbox{Tile$(\sqrt{3}, 1)$}, respectively, 
which admit only non-periodic tilings as shown by Smith et al. \cite{ref5}. Note that 
these tiles have no matching rules\footnote{ 
Matching rules prescribe how tiles must be joined to form a valid tiling and are often 
represented by assigning colors or orientations of arrows to designated edges of the 
tiles. For example, Penrose tiles can be regarded as an aperiodic tile set consisting of 
two types of polygons with matching rules \cite{ref1, ref2, ref3}.
}. In this manuscript, the substitution method using the clusters $H_7$ and $H_8$ 
(Figure~\ref{Fig.1-3}), shown in Figure 2.11 of \cite{ref5}, is referred to as the $H_7/H_8$ 
substitution. The non-periodic tiling generated using the corresponding clusters and the 
$H_7/H_8$ substitution is denoted by $\mathscr{T}_h$.

Around September 2023, I introduced the relationship between hat and turtle tiles and 
the Type 5 tiling on a mailing list for tiling enthusiasts and on my webpage \cite{ref8}.

The main idea presented in \cite{ref8} can be summarized as follows and is referred 
to here as Sugimoto's Perspective.

\begin{nameth}
A pattern of periodic tilings, such as the Type 5 tiling, can be made to correspond to 
$\mathscr{T}_h$ using hat or turtle tiles. As a result, one obtains hat or turtle tiles whose 
anterior and posterior sides bear different surface patterns. By combining these patterned 
hat or turtle tiles so that their surface patterns form periodic tilings, the non-periodic 
tiling $\mathscr{T}_h$ with hat or turtle tiles can be generated.
\end{nameth}

\noindent\textbf{Note:} As mentioned in \cite{ref8}, the periodic tiling pattern used 
in this correspondence is not restricted to the Type 5 tiling. We expect the same 
idea to apply to tilings belonging to the wallpaper group such as p6m and p6 
(see \cite{ref1}, \cite{ref3}, etc. for details on these wallpaper groups).

\vspace{\baselineskip}

On the other hand, from an email sent by Smith to a mailing list of tiling enthusiasts 
in January 2026, we learned that Amfirifma had noticed the same property and 
obtained a result concerning two types of concave polygons that can be formed 
from such a pattern. We interpreted Amfirifma's result and inferred the following perspective.

\begin{nameth2}
A pattern of periodic tilings, such as the Type 5 tiling, can be made to correspond 
to $\mathscr{T}_h$ using hat or turtle tiles. By combining the polygons that constitute 
the tiling pattern, one obtains a tile set consisting of two types of tiles (concave polygons) 
 corresponding to the anterior and posterior sides of a hat or turtle tile. 
Using the $H_7/H_8$ substitution, this tile set can generate $\mathscr{T}_h$ in which 
the densities of the two types of tiles are mutually reversed. Moreover, 
in some cases, each concave polygon in the tile set does not admit monohedral tilings 
(cannot generate a tiling on its own)\footnote{ 
A ``pattern'' is a design formed by geometric diagrams. For example, a pattern of a Type 5 
tiling (or a Type 5 tiling pattern) is a design on the plane corresponding to a Type 5 tiling 
formed by convex pentagonal monotiles belonging to the Type 5 family. A patterned tile is a 
tile whose surface (i.e., the faces of the tile, which are referred to as the anterior side and 
posterior side) bears a design formed by geometric diagrams. In this manuscript, a hat tile 
whose surface bears a design consisting of convex pentagonal monotiles belonging to the 
Type 5 family or line segments forming parts of such monotiles is referred to as a patterned 
hat tile. See Section~\ref{section2} and subsequent sections for details.
}.
\end{nameth2}

The purpose of this manuscript is to explain the relationship between Sugimoto's Perspective 
and Amfirifma's Perspective and to present related observations and investigations.

In this study, when a tile set contains $n$ types of tiles, we use the expression ``tile set 
consisting of $n$ types of tiles.'' This is because we take into consideration the fact that 
this study focuses on the property described in Amfirifma's Perspective that ``the densities 
of the two types of tiles are mutually reversed.'' In general, when a tile set contains $n$ types 
of tiles, it may also be described as a ``tile set consisting of $n$ tiles.'' However, using this 
expression can result in a phrase such as ``tile sets consisting of two tiles that have the 
property that the densities of the two types of tiles are mutually reversed.'' In this expression, 
although the same two types of tiles are being referred to, the expressions ``two types of tiles'' 
and ``two tiles'' are mixed, which may make the meaning difficult to understand. On the other 
hand, if ``types'' is replaced by ``number,'' the expression becomes ``the property that the 
densities of the two tiles are mutually reversed,'' which may not accurately convey the 
property described in Amfirifma's Perspective. Therefore, in this study, we do not use the 
expression ``tile set consisting of $n$ tiles'' and consistently use ``tile set consisting of 
$n$ types of tiles.''

Furthermore, it should be noted that, as is customary in many studies of tilings, we consider 
the anterior-side and posterior-side of a tile as the same type of tile in this study. 
For example, an anterior-side hat tile and a posterior-side hat tile are not regarded as 
different types of tiles. Consequently, a pair consisting of an anterior-side hat tile and a 
posterior-side hat tile is not regarded as a tile set consisting of two types of tiles.

\begin{figure}[H]
 \centering\includegraphics[width=12.5cm,clip]{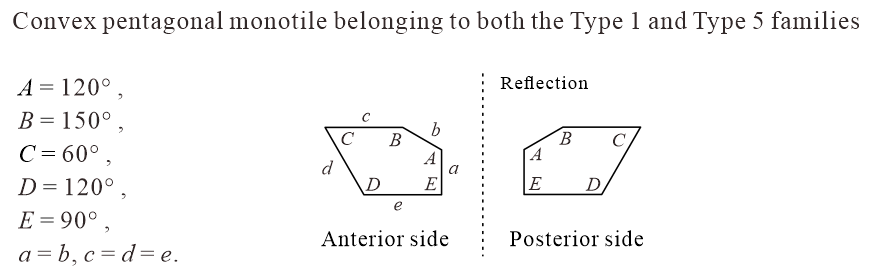} 
  \caption{{\small Convex pentagonal monotile belonging to both the Type 1 and Type 5 families.} }
\label{Fig.1-4}
\end{figure}

In this manuscript, the convex pentagonal monotile shown in Figure~\ref{Fig.1-4}, which belongs 
to both the Type 1 and Type 5 families, was used for the investigation, analysis, and drawing 
of diagrams and patterns related to convex pentagonal monotiles of Type 5. We chose this convex 
pentagon because we are familiar with it and because its vertex angles make it convenient 
for investigation, analysis, and geometric construction \cite{ref10, ref11, ref12}.

\vspace{2\baselineskip}
\section{Relationship between $\mathscr{T}_h$ with hat or turtle tiles and Type 5 tiling}
\label{section2}

\subsection{Explanation for the case of the hat tile} 
\label{subsection2.1}

\begin{figure}[t]
 \centering\includegraphics[width=15cm,clip]{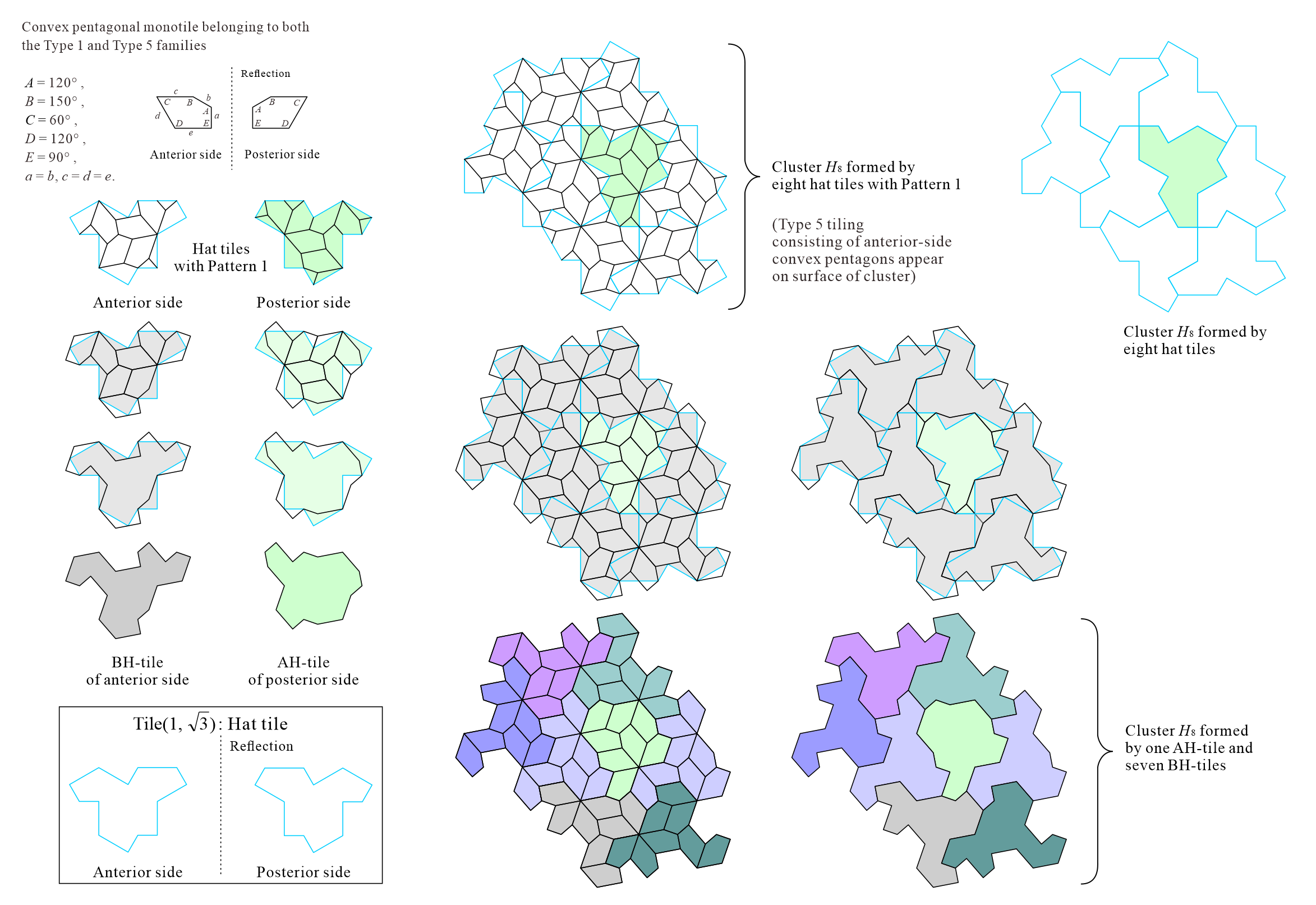} 
  \caption{{\small Relationship between hat tiles with Pattern 1 and the corresponding AH-tiles and BH-tiles.} }
\label{Fig.2-1}
\end{figure}

\begin{figure}[t]
 \centering\includegraphics[width=15cm,clip]{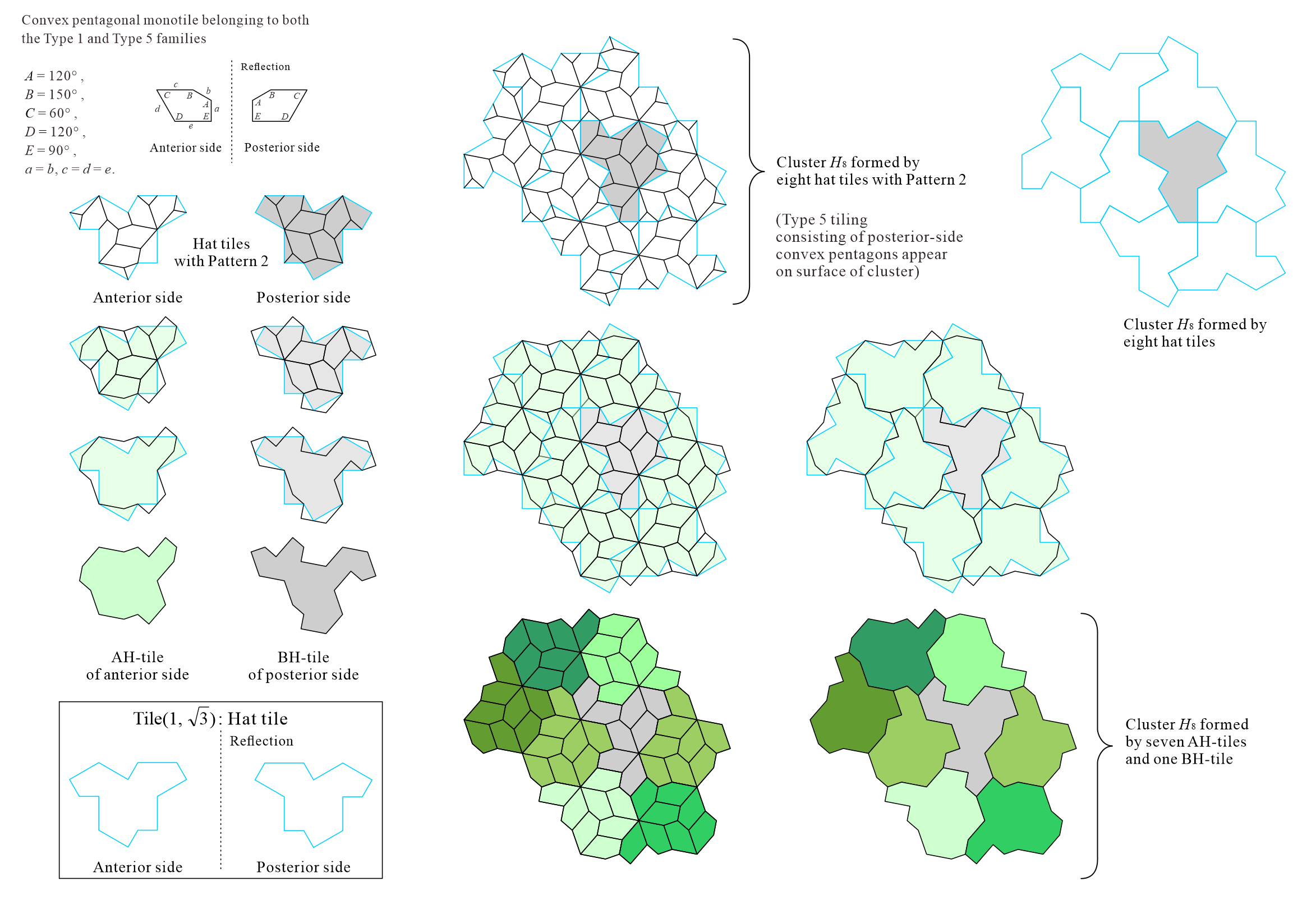} 
  \caption{{\small Relationship between hat tiles with Pattern 2 and the corresponding AH-tiles and BH-tiles.} }
\label{Fig.2-2}
\end{figure}

Figures~\ref{Fig.2-1} and \ref{Fig.2-2} show the correspondence between the Type 5 tiling 
(formed by the convex pentagonal monotiles belonging to both the Type 1 and Type 5 families 
shown in the upper-left corner of each figure) and the hat tile and cluster $H_8$. These figures 
also illustrate the relationship between Sugimoto's Perspective and Amfirifma's Perspective. 
By placing the red-circled vertices of the hat tiles shown in Figure~\ref{Fig.1-2} at the 
valence-six vertices of the Type 5 tiling (i.e., vertices at which six convex pentagons meet) 
and superimposing the resulting hat tilings on the Type 5 tiling, one obtains the correspondence 
shown in Figures~\ref{Fig.2-1} and \ref{Fig.2-2}.

In Figure~\ref{Fig.2-1}, the patterns applied to the hat tiles correspond to the Type 5 tiling 
formed by the convex pentagons designated as the anterior side in the figure. In Sugimoto's 
Perspective, the case shown in Figure~\ref{Fig.2-1} is called Pattern 1, and a hat tile bearing 
the corresponding pattern is referred to as a hat tile with Pattern 1. In other words, 
if cluster $H_8$ is formed by hat tiles with Pattern 1 as shown in Figure~\ref{Fig.2-1}, 
the Type 5 tiling consisting of anterior-side convex pentagons appears on its surface.

In Figure~\ref{Fig.2-2}, the patterns applied to the hat tiles correspond to the Type 5 tiling 
formed by the convex pentagons designated as the posterior side in the figure. In Sugimoto's 
Perspective, the case shown in Figure~\ref{Fig.2-2} is called Pattern 2, and a hat tile bearing 
the corresponding pattern is referred to as a hat tile with Pattern 2. In other words, 
if cluster $H_8$ is formed by hat tiles with Pattern 2 as shown in Figure~\ref{Fig.2-2}, 
the Type 5 tiling consisting of posterior-side convex pentagons appears on its surface.

When forming cluster $H_8$ with hat tiles, seven anterior-side hat tiles and one posterior-side 
hat tile are used, as shown in the upper-right portions of Figures~\ref{Fig.2-1} and \ref{Fig.2-2}.

\begin{figure}[H]
 \centering\includegraphics[width=15cm,clip]{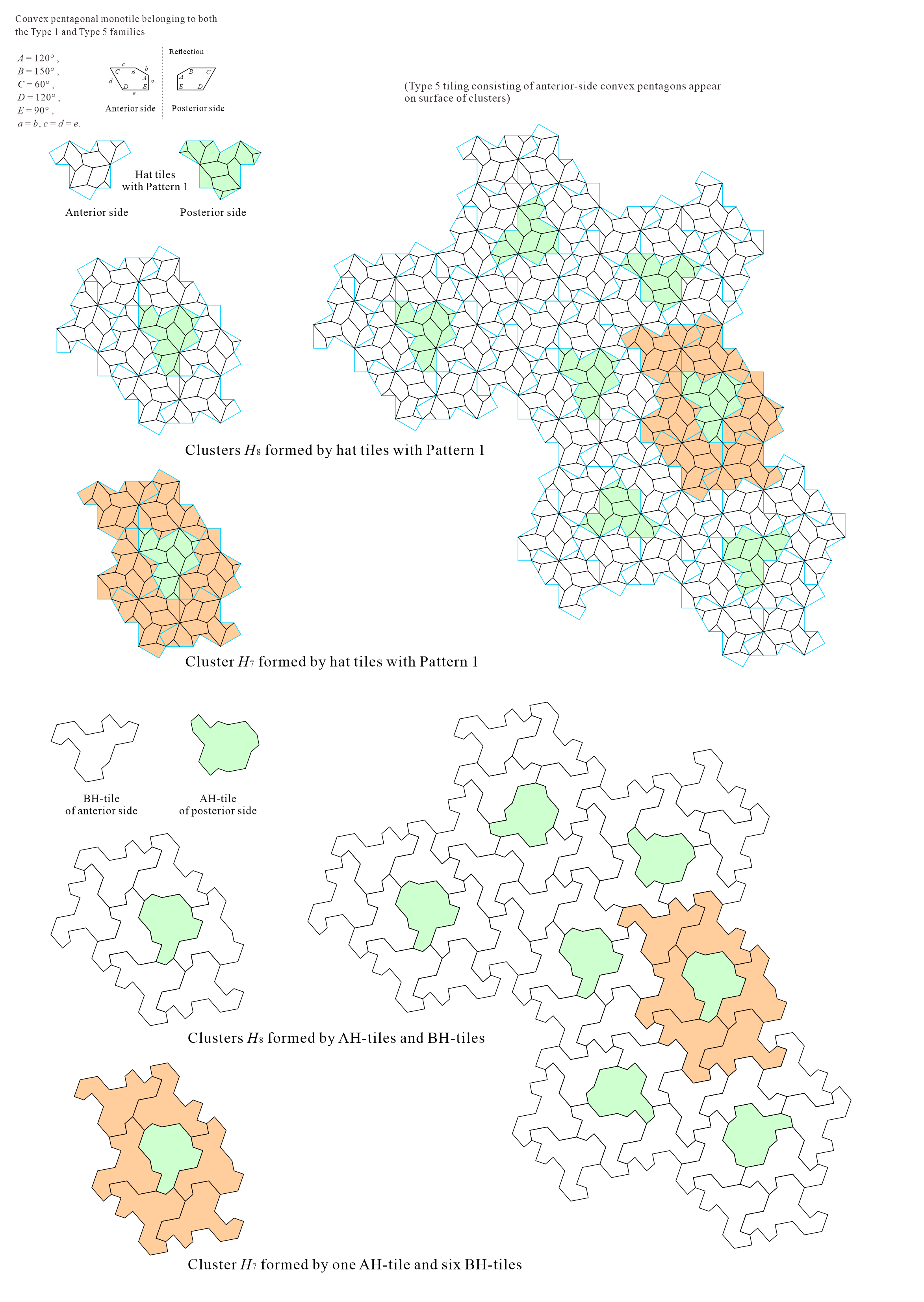} 
  \caption{{\small Cluster $H_8$ formed by hat tiles with Pattern 1, the corresponding 
  cluster $H_8$ formed by AH-tiles and BH-tiles and the next-step cluster $H_8$.} }
\label{Fig.2-3}
\end{figure}

When Amfirifma's Perspective is applied to the hat tile with Pattern 1 as shown in 
Figure~\ref{Fig.2-1}, two types of concave polygons, each consisting of eight convex pentagons, 
are obtained. The concave polygon corresponding to the anterior-side hat tile with Pattern 1 
is called the BH-tile, whereas the concave polygon corresponding to the posterior-side hat 
tile with Pattern 1 is called the AH-tile. The distinction between the anterior and posterior 
sides of the AH-tile and BH-tile corresponds to that of the original hat tile. As shown in 
Figure~\ref{Fig.2-1}, the cluster $H_8$ formed by hat tiles with Pattern 1 corresponds exactly 
to the cluster $H_8$ formed by one AH-tile and seven BH-tiles. Therefore, the $H_7/H_8$ 
substitution can be applied to the cluster $H_8$ formed by one AH-tile and seven BH-tiles, 
thereby producing the next-step cluster $H_8$ as shown in Figure~\ref{Fig.2-3}. Repeated 
application of the substitution generates $\mathscr{T}_h$.

\begin{figure}[H]
 \centering\includegraphics[width=15cm,clip]{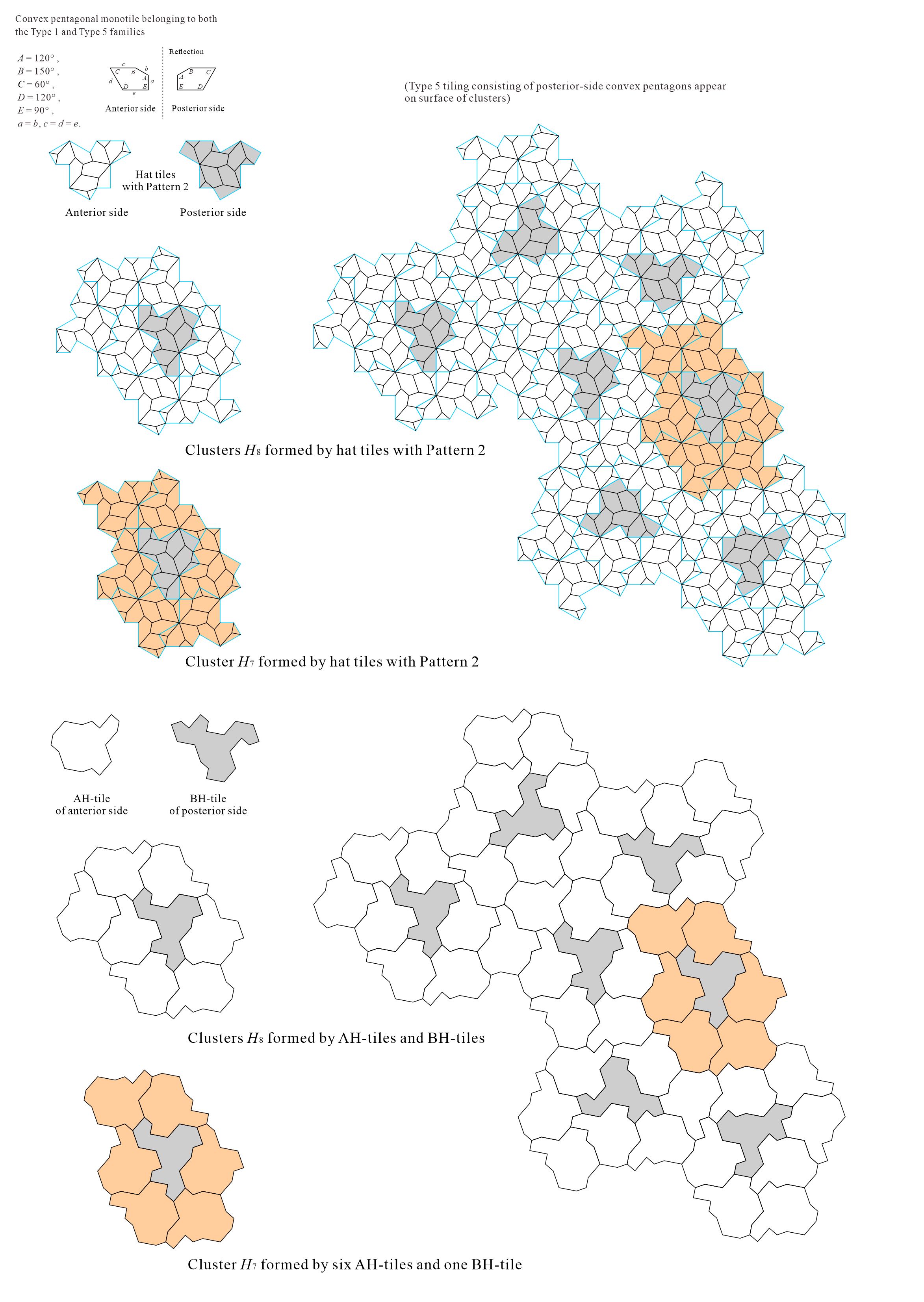} 
  \caption{{\small Cluster $H_8$ formed by hat tiles with Pattern 2, the corresponding 
  cluster $H_8$ formed by AH-tiles and BH-tiles and the next-step cluster $H_8$.} }
\label{Fig.2-4}
\end{figure}

Similarly, when Amfirifma's Perspective is applied to the hat tile with Pattern 2 shown in 
Figure~\ref{Fig.2-2}, the concave polygon corresponding to the anterior-side hat tile with 
Pattern 2 becomes the AH-tile, whereas the concave polygon corresponding to the 
posterior-side hat tile with Pattern 2 becomes the BH-tile. Note that the AH-tile and 
BH-tile in Figure~\ref{Fig.2-2} are mirror images of the AH-tile and BH-tile in Figure~\ref{Fig.2-1}, 
respectively. As a result, as shown in Figure~\ref{Fig.2-2}, the cluster $H_8$ formed by 
hat tiles with Pattern 2 corresponds exactly to the cluster $H_8$ formed by seven AH-tiles 
and one BH-tile. Therefore, the $H_7/H_8$ substitution can be applied to the cluster $H_8$ 
formed by seven AH-tiles and one BH-tile, thereby producing the next-step cluster $H_8$ as 
shown in Figure~\ref{Fig.2-4}. Repeated application of the substitution generates $\mathscr{T}_h$.

\subsection{Explanation for the case of the turtle tile} 
\label{subsection2.2}

\begin{figure}[H]
 \centering\includegraphics[width=14.75cm,clip]{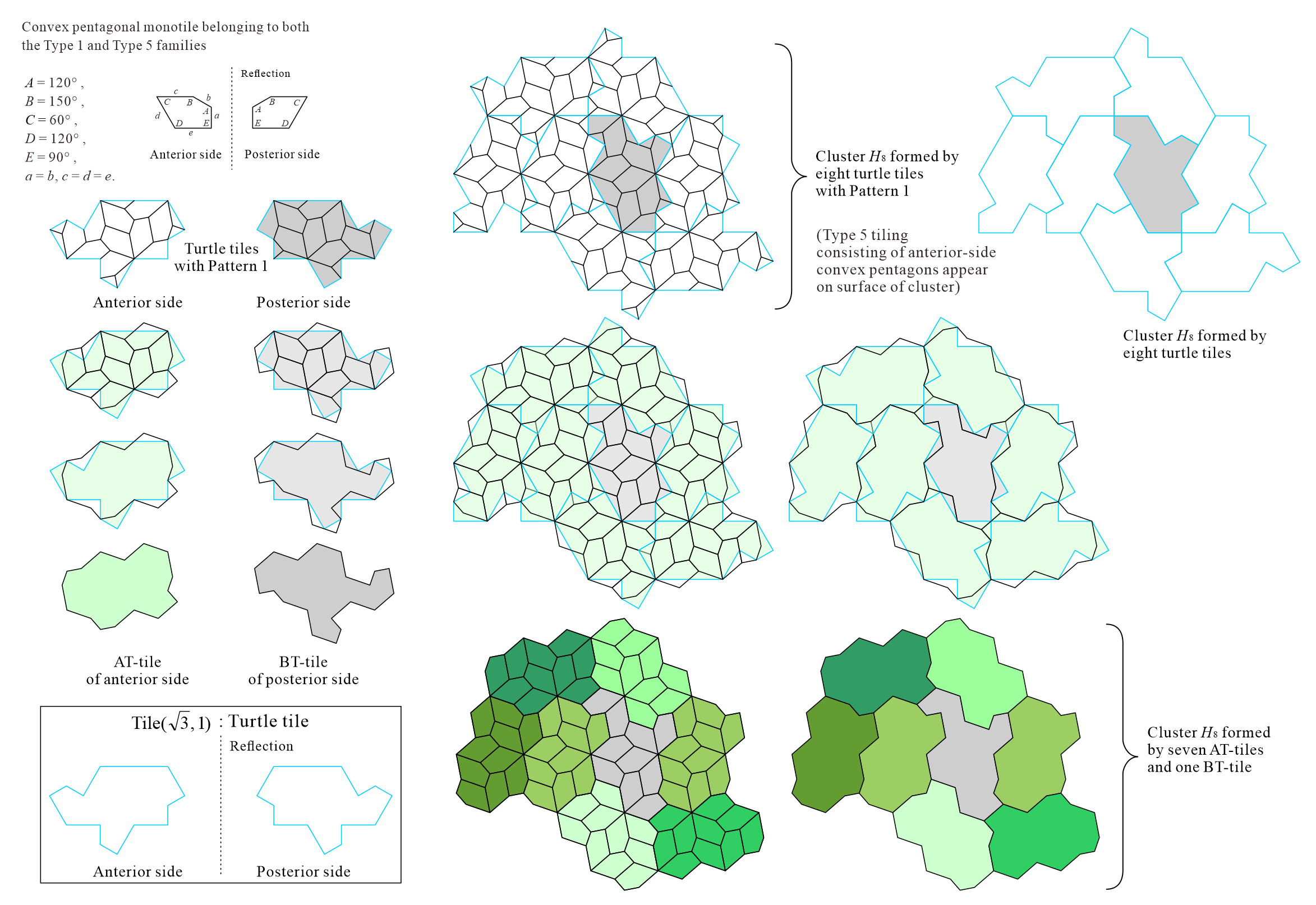} 
  \caption{{\small Relationship between turtle tiles with Pattern 1 and the corresponding AT-tiles and BT-tiles.} }
\label{Fig.2-5}
\end{figure}

\begin{figure}[H]
 \centering\includegraphics[width=14.75cm,clip]{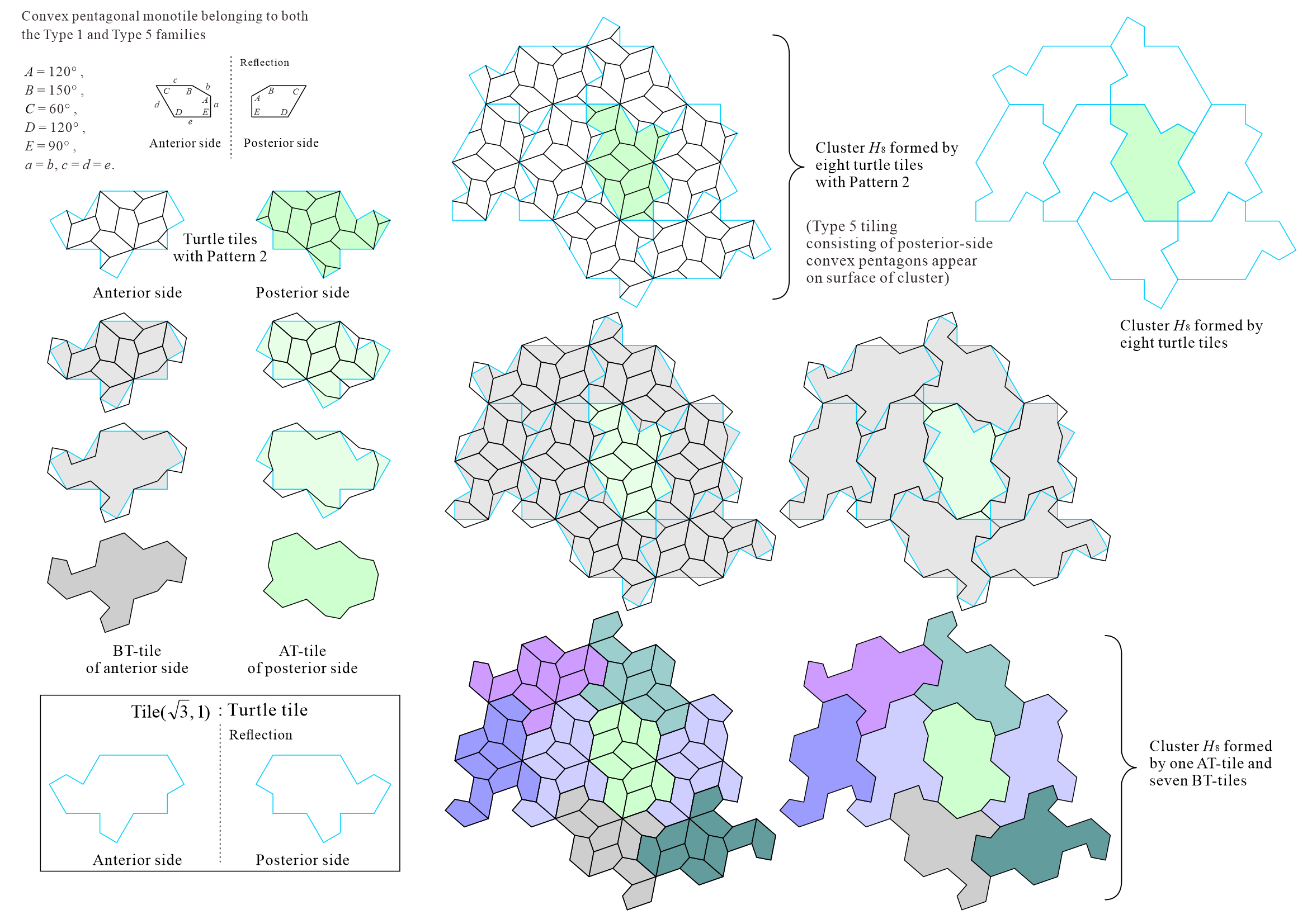} 
  \caption{{\small Relationship between turtle tiles with Pattern 2 and the corresponding AT-tiles and BT-tiles.} }
\label{Fig.2-6}
\end{figure}

Figures~\ref{Fig.2-5} and \ref{Fig.2-6} show the correspondence between the Type 5 tiling 
(formed by convex pentagonal monotiles belonging to both the Type 1 and Type 5 families 
shown in the upper-left corner of each figure) and the turtle tile and cluster $H_8$. 
These figures also illustrate the relationship between Sugimoto's Perspective and 
Amfirifma's Perspective. By placing the red-circled vertices of the turtle tiles shown in 
Figure~\ref{Fig.1-2} at the valence-six vertices of the Type 5 tiling (i.e., vertices at 
which six convex pentagons meet) and superimposing the resulting turtle tilings on the 
Type 5 tiling, one obtains the correspondence shown in Figures~\ref{Fig.2-5} and \ref{Fig.2-6}.

In Figure~\ref{Fig.2-5}, the patterns applied to the turtle tiles correspond to the Type 5 
tiling formed by convex pentagons designated as the anterior side in the figure. In Sugimoto's 
Perspective, the case shown in Figure~\ref{Fig.2-5} is called Pattern 1, and a turtle tile 
bearing the corresponding pattern is referred to as a turtle tile with Pattern 1. 
Owing to the properties of Pattern 1, a Type 5 tiling consisting of anterior-side convex 
pentagons appears on the surface of $H_8$.

In Figure~\ref{Fig.2-6}, the patterns applied to the turtle tiles correspond to the Type 5 
tiling formed by convex pentagons designated as the posterior side in the figure. In Sugimoto's 
Perspective, the case shown in Figure~\ref{Fig.2-6} is called Pattern 2, and a turtle tile 
bearing the corresponding pattern is referred to as a turtle tile with Pattern 2. 
Owing to the properties of Pattern 2, a Type 5 tiling consisting of posterior-side convex 
pentagons appears on the surface of $H_8$.

When forming cluster $H_8$ with turtle tiles, seven anterior-side turtle tiles and one 
posterior-side turtle tile are used, as shown in the upper-right portions of 
Figures~\ref{Fig.2-5} and \ref{Fig.2-6}.

\begin{figure}[H]
 \centering\includegraphics[width=15cm,clip]{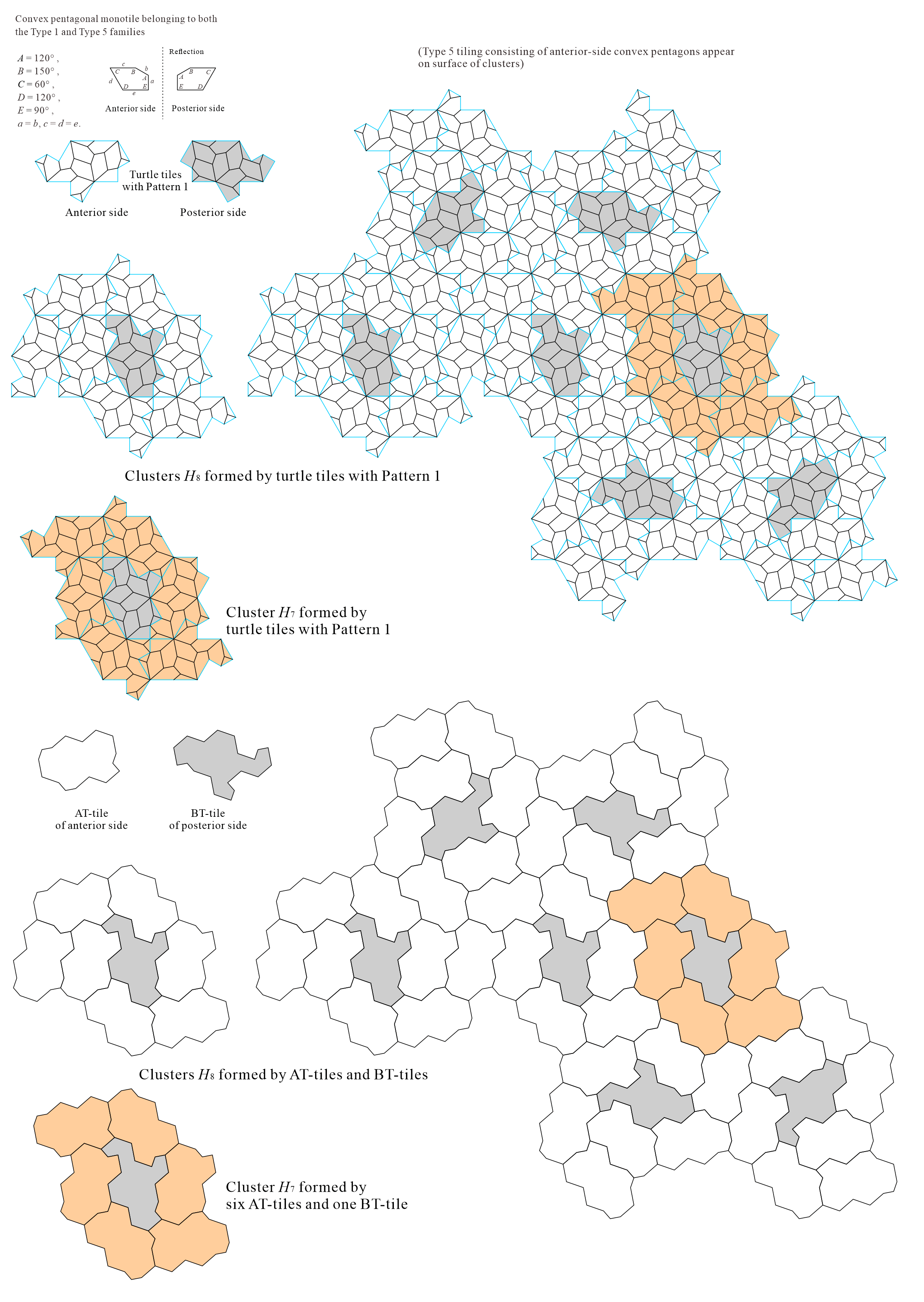} 
  \caption{{\small Cluster $H_8$ formed by turtle tiles with Pattern 1 and the corresponding 
  cluster $H_8$ formed by AT-tiles and BT-tiles.} }
\label{Fig.2-7}
\end{figure}

When Amfirifma's Perspective is applied to the turtle tile with Pattern 1 as shown in 
Figure~\ref{Fig.2-5}, two types of concave polygons, each consisting of 10 convex pentagons, 
are obtained. The concave polygon corresponding to the anterior-side turtle tile with Pattern 1 
is called the AT-tile, whereas the concave polygon corresponding to the posterior-side turtle 
tile with Pattern 1 is called the BT-tile. The distinction between the anterior and posterior 
sides of the AT-tile and BT-tile corresponds to that of the original turtle tile. As shown in 
Figure~\ref{Fig.2-5}, the cluster $H_8$ formed by turtle tiles with Pattern 1 corresponds exactly 
to the cluster $H_8$ formed by seven AT-tiles and one BT-tile. Therefore, the $H_7/H_8$ 
substitution can be applied to the cluster $H_8$ formed by seven AT-tiles and one BT-tile, 
thereby producing the next-step cluster $H_8$ as shown in Figure~\ref{Fig.2-7}. Repeated 
application of the substitution generates $\mathscr{T}_h$.

\begin{figure}[H]
 \centering\includegraphics[width=15cm,clip]{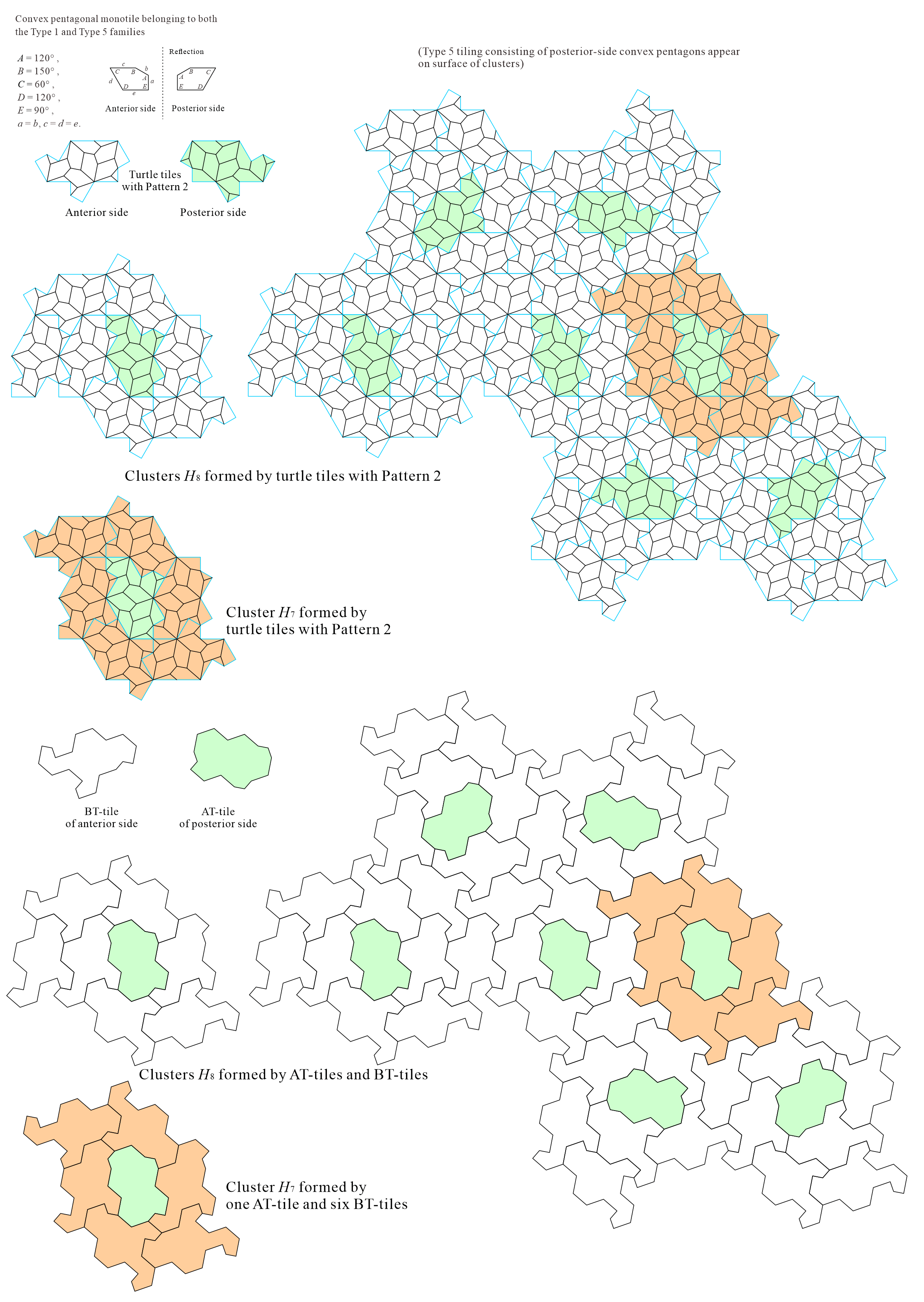} 
  \caption{{\small Cluster $H_8$ formed by turtle tiles with Pattern 2 and the corresponding 
  cluster $H_8$ formed by AT-tiles and BT-tiles.} }
\label{Fig.2-8}
\end{figure}

Similarly, when Amfirifma's Perspective is applied to the turtle tile with Pattern 2 shown in 
Figure~\ref{Fig.2-6}, the concave polygon corresponding to the anterior-side turtle tile with 
Pattern 2 becomes the BT-tile, whereas the concave polygon corresponding to the posterior-side 
turtle tile with Pattern 2 becomes the AT-tile. Note that the AT-tile and BT-tile in 
Figure~\ref{Fig.2-6} are mirror images of the AT-tile and BT-tile in Figure~\ref{Fig.2-5}, respectively. 
As a result, as shown in Figure~\ref{Fig.2-6}, the cluster $H_8$ formed by the turtle tiles with 
Pattern 2 corresponds exactly to the cluster $H_8$ formed by one AT-tile and seven BT-tiles. 
Therefore, the $H_7/H_8$ substitution can be applied to the cluster $H_8$ formed by one 
AT-tile and seven BT-tiles, thereby producing the next-step cluster $H_8$ as shown in 
Figure~\ref{Fig.2-8}. Repeated application of the substitution generates $\mathscr{T}_h$.

\subsection{Reversal of the densities of two types of tiles in $\mathscr{T}_h$}
\label{subsection2.3}

Let $H_8(n)$ denote the cluster corresponding to $H_8$ obtained after $n$ steps of $H_7/H_8$ 
substitution. However, the cluster $H_8$ in Figure~\ref{Fig.1-3} is regarded as Step 2 and is 
designated as $H_8(2)$.

In Figure~\ref{Fig.2-1}, the Step-2 cluster $H_8(2)$ which is formed by one AH-tile and seven 
BH-tiles is denoted by \mbox{$H_8(\text{AH}+7\text{BH}, 0)$}, and in Figure~\ref{Fig.2-2}, 
the cluster $H_8(2)$ formed by seven AH-tiles and one BH-tile is denoted by 
\mbox{$H_8(7\text{AH}+\text{BH}, 0)$}. \mbox{$H_8(\text{AH}+7\text{BH}, m)$} denotes the cluster 
corresponding to $H_8$ obtained after $m$ steps of the $H_7/H_8$ substitution starting from 
\mbox{$H_8(\text{AH}+7\text{BH}, 0)$}. Thus, \mbox{$H_8(\text{AH}+7\text{BH}, 1)$} corresponds 
to the Step-3 cluster $H_8(3)$. Similarly, \mbox{$H_8(7\text{AH}+\text{BH}, m)$} denotes the cluster 
corresponding to $H_8$ obtained after $m$ steps of the $H_7/H_8$ substitution starting from 
\mbox{$H_8(7\text{AH}+\text{BH}, 0)$}. Accordingly, \mbox{$H_8(\text{AH}+7\text{BH}, \infty)$} 
and \mbox{$H_8(7\text{AH}+\text{BH}, \infty)$} denote $\mathscr{T}_h$ generated by AH-tiles and 
BH-tiles, respectively.

In Figure~\ref{Fig.2-5}, the cluster $H_8(2)$ formed by seven AT-tiles and one BT-tile is denoted 
by \mbox{$H_8(7\text{AT}+\text{BT}, 0)$}, and in Figure~\ref{Fig.2-6},  $H_8(2)$ formed by 
one AT-tile and seven BT-tiles is denoted by \mbox{$H_8(\text{AT}+7\text{BT}, 0)$}. 
Similar to the case of AH-tile and BH-tile, \mbox{$H_8(7\text{AT}+\text{BT}, m)$} and 
\mbox{$H_8(\text{AT}+7\text{BT}, m)$} denote the cluster corresponding to $H_8$ obtained after 
performing $H_7/H_8$ substitution $m$ steps. Therefore, \mbox{$H_8(7\text{AT}+\text{BT}, \infty)$} 
and \mbox{$H_8(\text{AT}+7\text{BT}, \infty)$} denote $\mathscr{T}_h$ generated by 
AT-tiles and BT-tiles, respectively.

The $H_7/H_8$ substitution can be represented by the matrix 
\[
\begin{pmatrix}
1 & 1 \\
5 & 6
\end{pmatrix}
\]
which has the eigenvalues $\frac{7 \pm 3\sqrt{5}}{2}$. The eigenvector corresponding to 
the $H_7/H_8$ substitution is 
\[
\begin{pmatrix}
3\sqrt{5}-5 \\
10
\end{pmatrix}
\]
Therefore, the ratio of clusters $H_7$ and $H_8$ approaches $3\sqrt{5}-5 : 10$.

Each cluster has one reflected tile (posterior-side tile) and six or seven unreflected ones 
(anterior-side tiles), which means the density of the reflected tiles in $\mathscr{T}_h$ is
\[
\frac{(3\sqrt{5}-5)+10}{7(3\sqrt{5}-5)+8\times 10}
\]
which is equal to $\frac{3-\sqrt{5}}{6} \approx 0.127322\ldots$ \cite{ref4}.

Therefore, the following statements hold:

\begin{itemize}
\item In \mbox{$H_8(\text{AH}+7\text{BH}, \infty)$}, the density of AH-tiles is $\frac{3-\sqrt{5}}{6} \approx 0.127322\ldots$.
\item In \mbox{$H_8(7\text{AH}+\text{BH}, \infty)$}, the density of BH-tiles is $\frac{3-\sqrt{5}}{6} \approx 0.127322\ldots$.
\item In \mbox{$H_8(\text{AT}+7\text{BT}, \infty)$}, the density of AT-tiles is $\frac{3-\sqrt{5}}{6} \approx 0.127322\ldots$.
\item In \mbox{$H_8(7\text{AT}+\text{BT}, \infty)$}, the density of BT-tiles is $\frac{3-\sqrt{5}}{6} \approx 0.127322\ldots$.
\end{itemize}

These properties correspond to the statement in Amfirifma's Perspective that the tile set can 
generate $\mathscr{T}_h$ in which the densities of the two types of tiles are mutually reversed. 

A tiling is \textit{dihedral} if it is generated using two types of prototiles in the tile set (set of prototiles). 
More generally, a tiling is \textit{$n$-hedral} if it is generated using $n$ types of tiles in the tile set \cite{ref1}. 
The tilings $\mathscr{T}_h$ generated by hat or turtle tiles are monohedral. On the other hand, 
\mbox{$H_8(\text{AH}+7\text{BH}, \infty)$}, \mbox{$H_8(7\text{AH}+\text{BH}, \infty)$}, 
\mbox{$H_8(7\text{AT}+\text{BT}, \infty)$}, and \mbox{$H_8(\text{AT}+7\text{BT}, \infty)$} are 
dihedral and are regarded as distinct patterns. Therefore, Amfirifma's Perspective can be viewed 
as a method for transforming a monohedral tiling generated by hat or turtle tiles into a dihedral 
tiling generated by a tile set consisting of two types of tiles. However, there is no guarantee 
that the tile set consisting of two types of tiles resulting from this transformation is aperiodic.

\vspace{2\baselineskip}
\section{Relationship between $\mathscr{T}_s$ generated by hat and turtle tiles and the Type 5 tiling}
\label{section3}

\subsection{Relationship between clusters formed by hat and turtle tiles and \mbox{Tile$(1, 1)$}} 
\label{subsection3.1}

\begin{figure}[hb]
 \centering\includegraphics[width=15cm,clip]{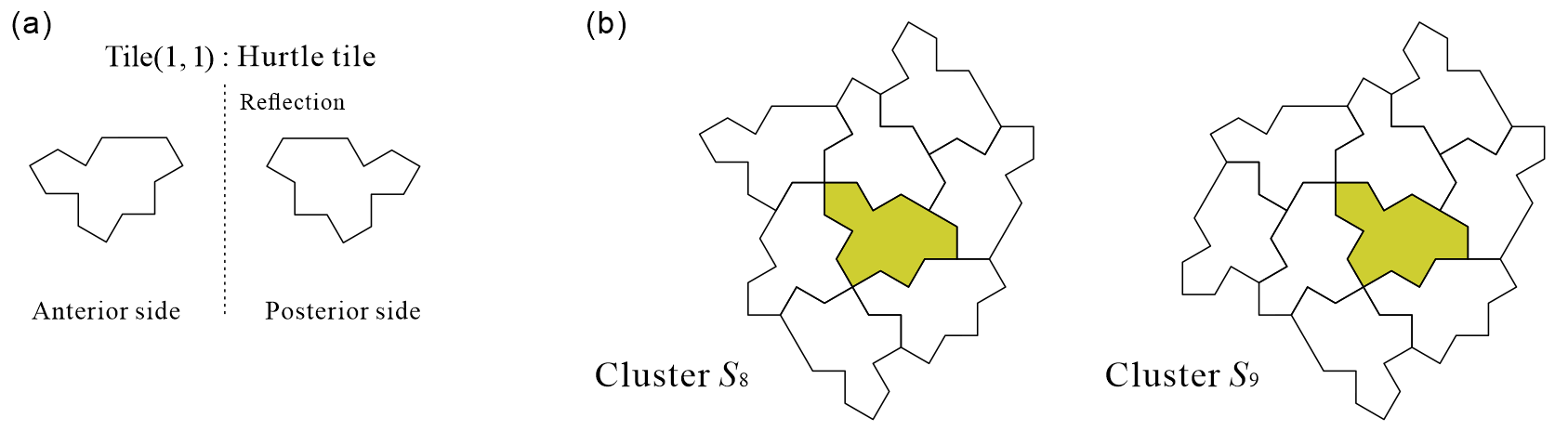} 
  \caption{{\small \mbox{Tile$(1, 1)$} and the clusters $S_8$ and $S_9$.} }
\label{Fig.3-1}
\end{figure}

Smith et al. \cite{ref6} showed that \mbox{Tile$(1, 1)$}, shown in Figure~\ref{Fig.3-1}(a), admits only 
non-periodic tilings if its reflection is prohibited during tiling generation. This \mbox{Tile$(1, 1)$} is 
referred to as chiral aperiodic monotile\footnote{ 
More precisely, \mbox{Tile$(1, 1)$} is classified as a ``weakly chiral aperiodic monotile'' \cite{ref6}. 
Then, \mbox{Tile$(1, 1)$} is also sometimes referred to as the ``hurtle tile,'' as it represents an entity 
that is intermediate between the hat tile and the turtle tile, though it is closer in nature to the turtle tile.
}, and the corresponding substitution method is described in \cite{ref6}.

In this manuscript, the cluster formed by nine \mbox{Tile$(1, 1)$} of Step 2, shown in Figure A.1 
of \cite{ref6}, is denoted by $S_9$. The cluster formed by eight \mbox{Tile$(1, 1)$} is denoted by 
$S_8$ (see Figure~\ref{Fig.3-1}(b)). Furthermore, we refer to the substitution method using $S_8$ 
and $S_9$, described in Appendix ``A. Constructing Spectre Tilings'' of \cite{ref6}, as the $S_8/S_9$ 
substitution. The non-periodic tiling generated using the corresponding clusters and the $S_8/S_9$ 
substitution is denoted by $\mathscr{T}_s$.

Section 3 of \cite{ref6} shows that there is a one-to-one correspondence (bijection) between tilings 
by \mbox{Tile$(1, 1)$} and tilings by combinations of hat tiles and turtle tiles, demonstrating that 
these two types of tilings are combinatorially equivalent. As shown in Figure 3.1 of \cite{ref6}, 
if the lengths of the even edges throughout the tiling $\mathscr{T}_s$ by \mbox{Tile$(1, 1)$} 
are increased from 1 to $\sqrt{3}$, even tiles are transformed into turtle tiles, whereas odd tiles 
into hat tiles. Conversely, if the lengths of the odd edges throughout the tiling $\mathscr{T}_s$ 
by \mbox{Tile$(1, 1)$} are increased to $\sqrt{3}$, even tiles are transformed into hat tiles, 
whereas odd tiles into turtle tiles. Note that the odd tiles correspond to the ocher-colored tiles in 
Figure~\ref{Fig.3-1}(b), whereas the even tiles correspond to the uncolored tiles in 
Figure~\ref{Fig.3-1}(b) (See \cite{ref6} for the definitions of even and odd edges).

In other words, clusters with hat and turtle tiles corresponding to clusters $S_8$ and $S_9$, which are formed 
by \mbox{Tile$(1, 1)$}, can be created by appropriately rescaling the tiles by making the length of the shortest 
edges the same. These clusters with hat and turtle tiles, together with the $S_8/S_9$ substitutions, can be 
used to generate a $\mathscr{T}_s$ containing hat and turtle tiles.

\begin{figure}[H]
 \centering\includegraphics[width=15cm,clip]{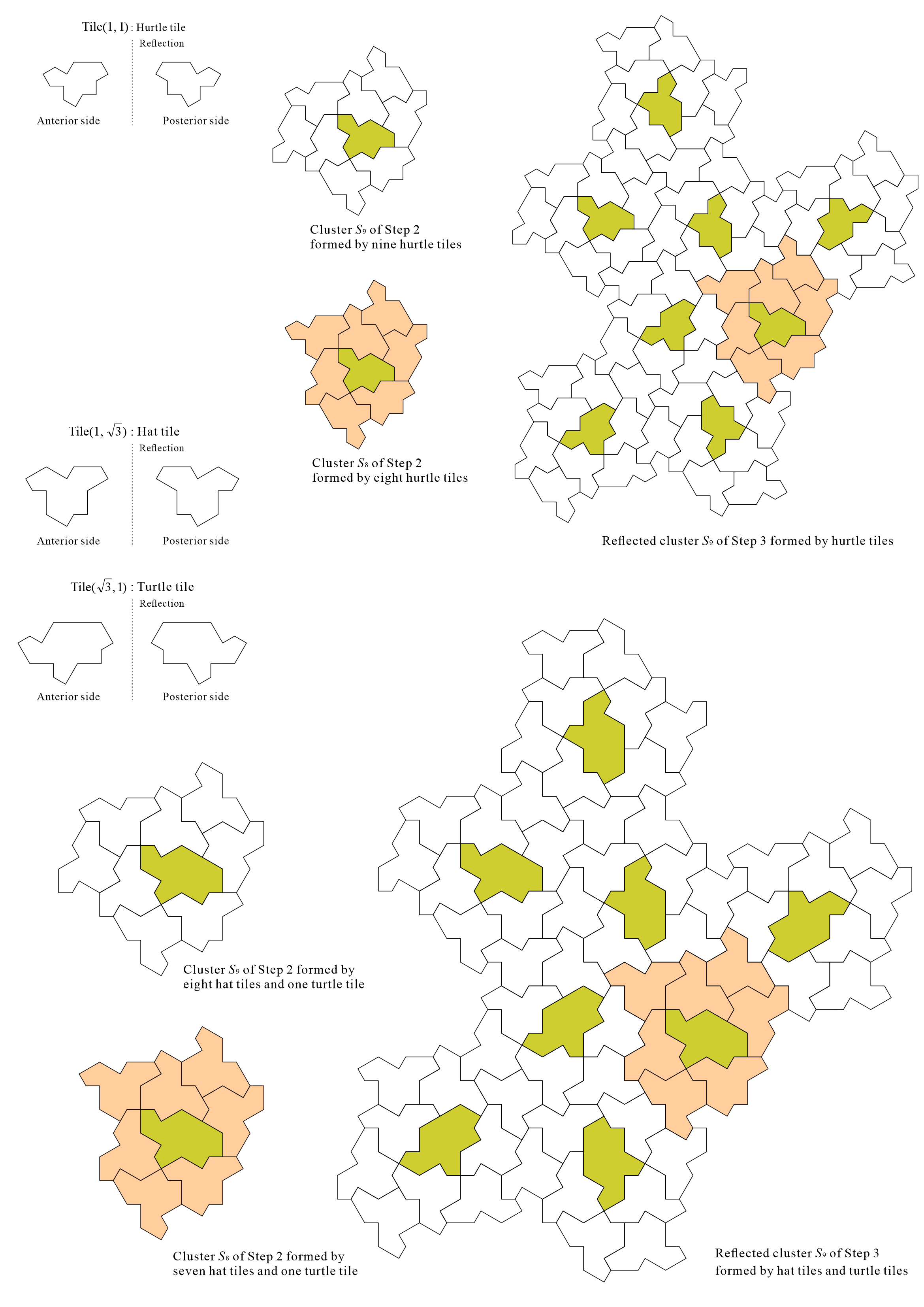} 
  \caption{{\small Clusters $S_8$ and $S_9$ formed by \mbox{Tile$(1, 1)$}, and corresponding clusters 
  $S_8$ and $S_9$ formed by placing turtle tiles in the odd tile positions and using hat tiles for the rest.} }
\label{Fig.3-2}
\end{figure}

\begin{figure}[H]
 \centering\includegraphics[width=15cm,clip]{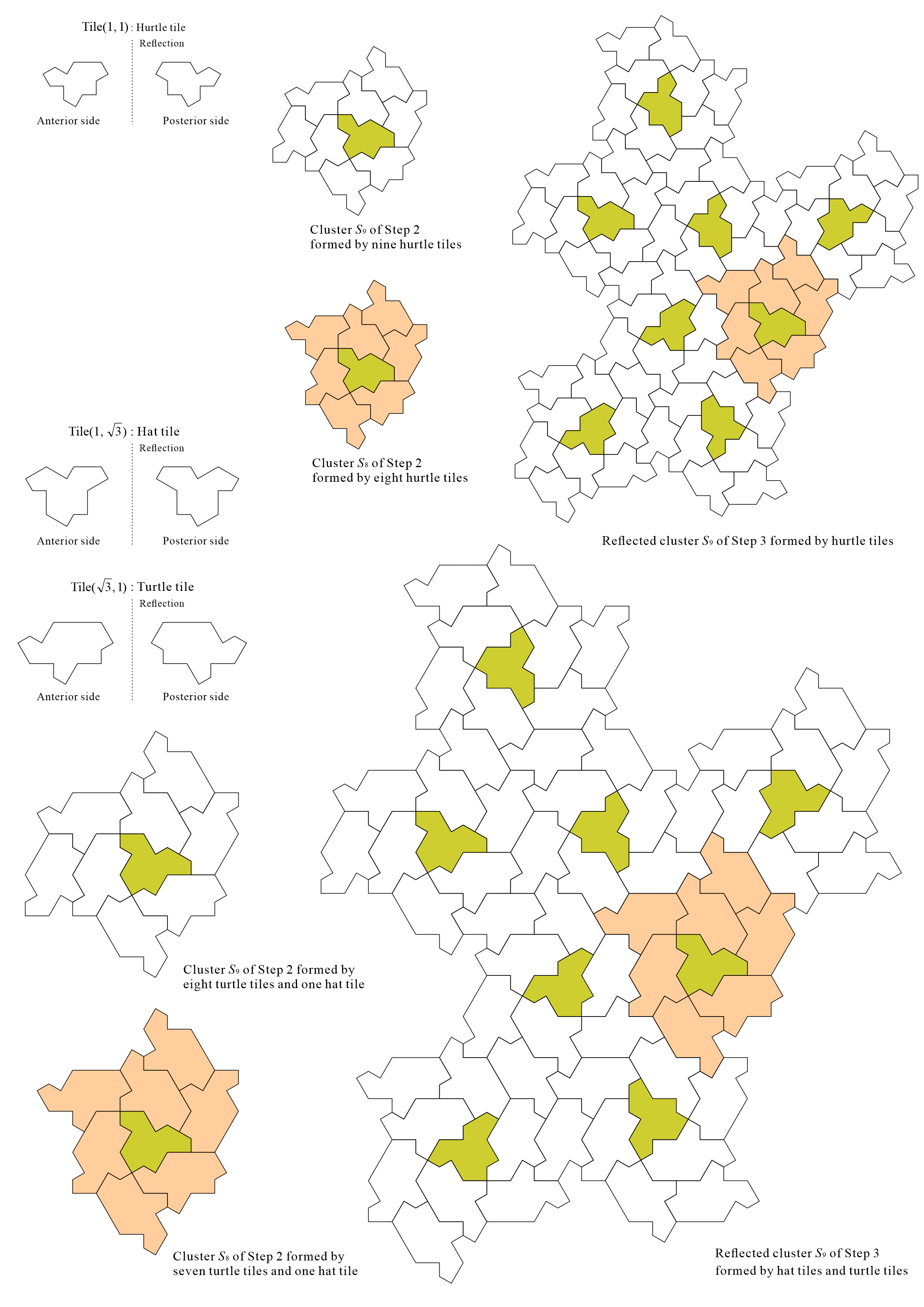} 
  \caption{{\small Clusters $S_8$ and $S_9$ formed by \mbox{Tile$(1, 1)$}, and corresponding clusters 
  $S_8$ and $S_9$ formed by placing hat tiles in the odd tile positions and using turtle tiles for the rest.} }
\label{Fig.3-3}
\end{figure}

Figures~\ref{Fig.3-2} and \ref{Fig.3-3} illustrate this correspondence. The upper parts of Figures~\ref{Fig.3-2} 
and \ref{Fig.3-3} depict the Step-2 clusters $S_8$ and $S_9$ and the reflected Step-3 cluster $S_9$ formed 
by \mbox{Tile$(1, 1)$}\footnote{ 
In this manuscript, the anterior and posterior sides of hurtle tile (i.e., \mbox{Tile$(1, 1)$}) are defined to 
correspond to the anterior and posterior sides of hat and turtle tiles shown in Figure~\ref{Fig.1-2}, 
as illustrated in Figure~\ref{Fig.3-1}. Consistent with the hat and turtle tiles in cluster $H_8$ discussed in 
Section~\ref{section2}, the hurtle tiles forming cluster $S_9$ are designated as being on their anterior sides 
in Figures~\ref{Fig.3-2} and \ref{Fig.3-3}. In the $S_8/S_9$ substitution, if the hurtle tiles in the cluster of 
Step 2 are on their anterior sides, as shown in Figures~\ref{Fig.3-2} and \ref{Fig.3-3}, then the hurtle tiles 
in the cluster of Step 3 are on their posterior sides \cite{ref6}. However, in Figures~\ref{Fig.3-2} and \ref{Fig.3-3}, 
the cluster of Step 3 is reflected and reoriented so that all hurtle tiles are shown on their anterior sides. 
To facilitate visual comparison when tiling with hat and turtle tiles, and when comparing the tilings with 
the corresponding clusters, the $\mathscr{T}_s$ and their associated clusters depicted in this manuscript 
are arranged so that hurtle tiles appear on their anterior sides. \label{footnote6}
}. Note that \mbox{Tile$(1, 1)$} corresponding to odd tiles are colored ochre. The lower part of 
Figure~\ref{Fig.3-2} depicts the corresponding Step-2 clusters $S_8$ and $S_9$ and the reflected Step-3 
cluster $S_9$, which are formed by placing turtle tiles in the odd tile positions and using hat tiles for the rest, 
corresponding to those shown in the upper part. The lower part of Figure~\ref{Fig.3-3} depicts the corresponding 
Step-2 cluster $S_8$ and $S_9$ and the reflected Step-3 cluster $S_9$, which are formed by placing hat tiles 
in the odd tile positions and using turtle tiles for the rest, corresponding to those shown in the upper part.

As noted above, $\mathscr{T}_s$ can be generated using hat and turtle tiles, and cluster $S_9$ can be 
formed by these tiles. This led us to conjecture that Sugimoto's Perspective introduced in 
Section~\ref{section1} might also apply to $\mathscr{T}_s$ generated by hat and turtle tiles. We verified 
that this is indeed the case. The results are presented below.

\subsection{$\mathscr{T}_s$ generated by hat and turtle tiles with Pattern 1}
\label{subsection3.2}

\begin{figure}[H]
 \centering\includegraphics[width=15cm,clip]{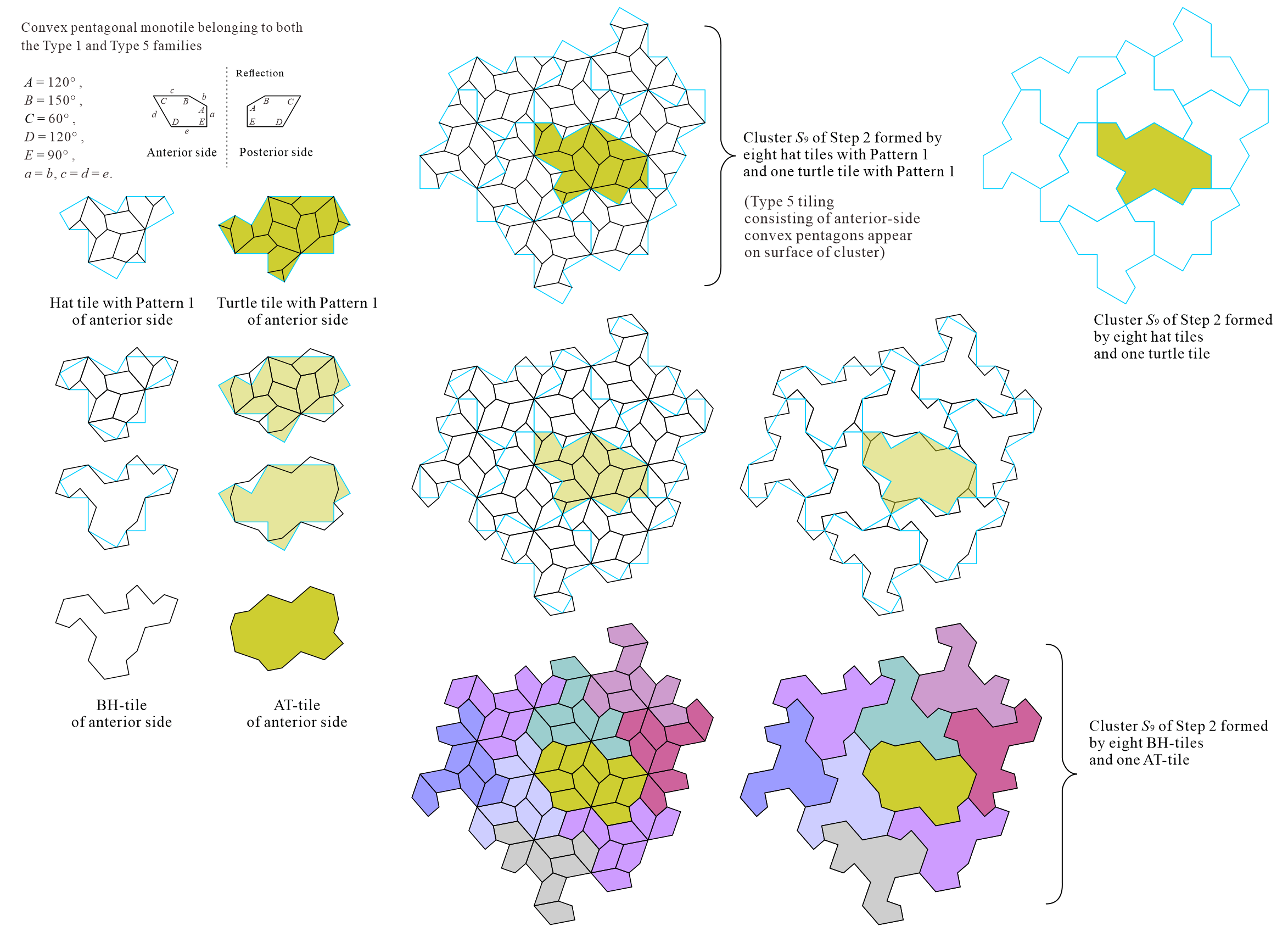} 
  \caption{{\small Step-2 cluster $S_9$ formed by eight hat tiles with Pattern 1 from Figure~\ref{Fig.2-1} 
  and one turtle tile with Pattern 1 from Figure~\ref{Fig.2-5}.} }
\label{Fig.3-4}
\end{figure}

\begin{figure}[H]
 \centering\includegraphics[width=15cm,clip]{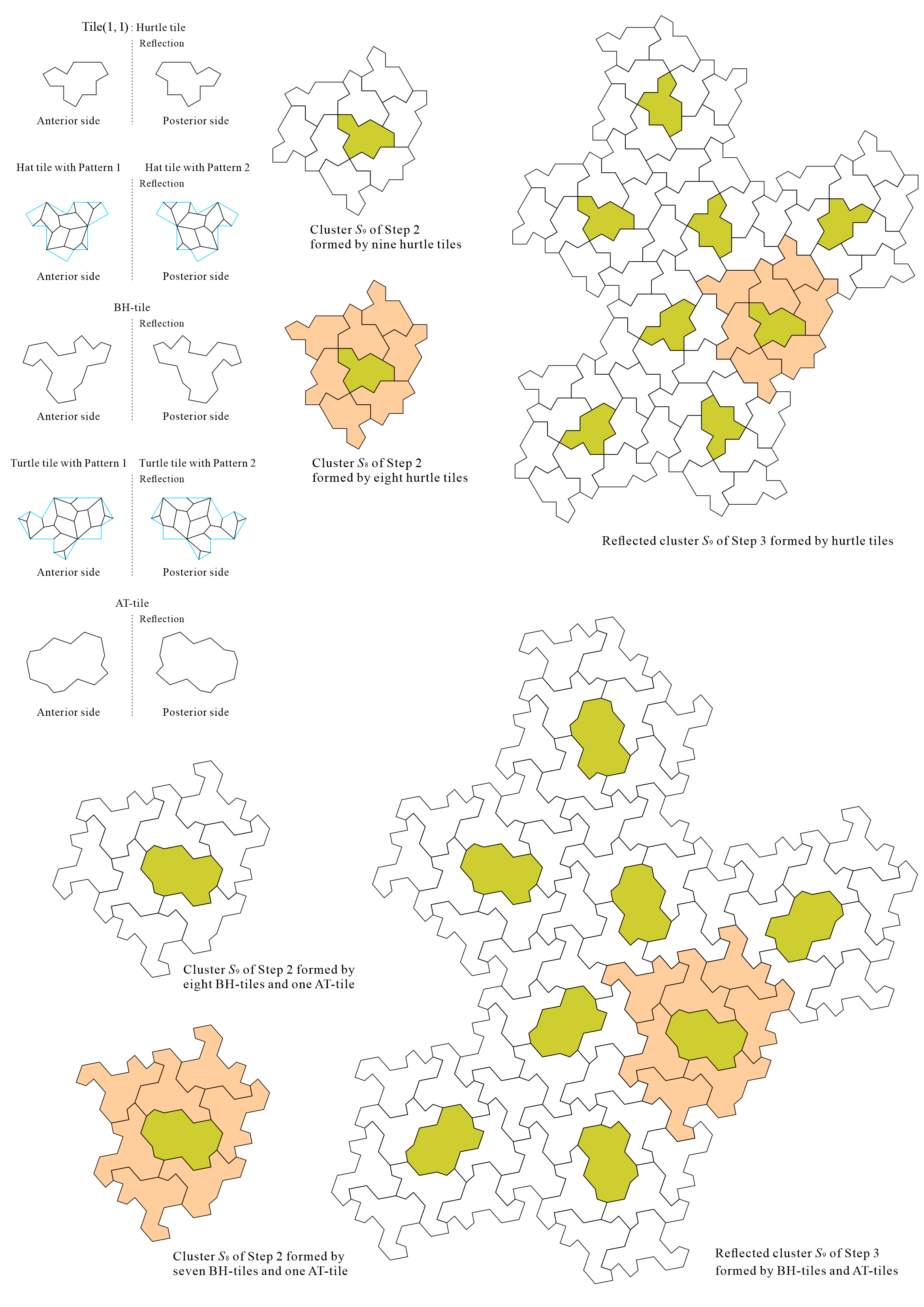} 
  \caption{{\small Application of the $S_8/S_9$ substitution to cluster $S_9$ formed by eight BH-tiles and one AT-tile.} }
\label{Fig.3-5}
\end{figure}

Figure~\ref{Fig.3-4} shows the Step-2 cluster $S_9$, formed by eight hat tiles and one turtle tile 
(corresponding to an odd tile), with the hat tile with Pattern 1 from Figure~\ref{Fig.2-1} and the turtle 
tile with Pattern 1 from Figure~\ref{Fig.2-5} applied. As explained in Footnote~\ref{footnote6}, this manuscript 
adopts representations in which the \mbox{Tile$(1, 1)$} in cluster $S_9$ are shown on their anterior sides. 
Therefore, the cluster $S_9$ shown in Figure~\ref{Fig.3-4} is formed using anterior-side hat and turtle tiles. 
Owing to the properties of Pattern 1, a Type 5 tiling consisting of anterior-side convex pentagons appears 
on the surface of $S_9$. As a result, the cluster $S_9$ formed by hat and turtle tiles with Pattern 1 
corresponds exactly to the cluster $S_9$ formed by eight BH-tiles and one AT-tile. Because the hat and turtle 
tiles forming $S_9$ in Figure~\ref{Fig.3-4} are on the anterior side, the BH-tiles and AT-tiles forming the 
corresponding $S_9$ are also on their anterior side. Therefore, the $S_8/S_9$ substitution can be applied 
to the cluster $S_9$ formed by eight BH-tiles and one AT-tile. This produces the next-step cluster $S_9$ 
shown in Figure~\ref{Fig.3-5}. Repeated application of the substitution generates $\mathscr{T}_s$.

\begin{figure}[H]
 \centering\includegraphics[width=15cm,clip]{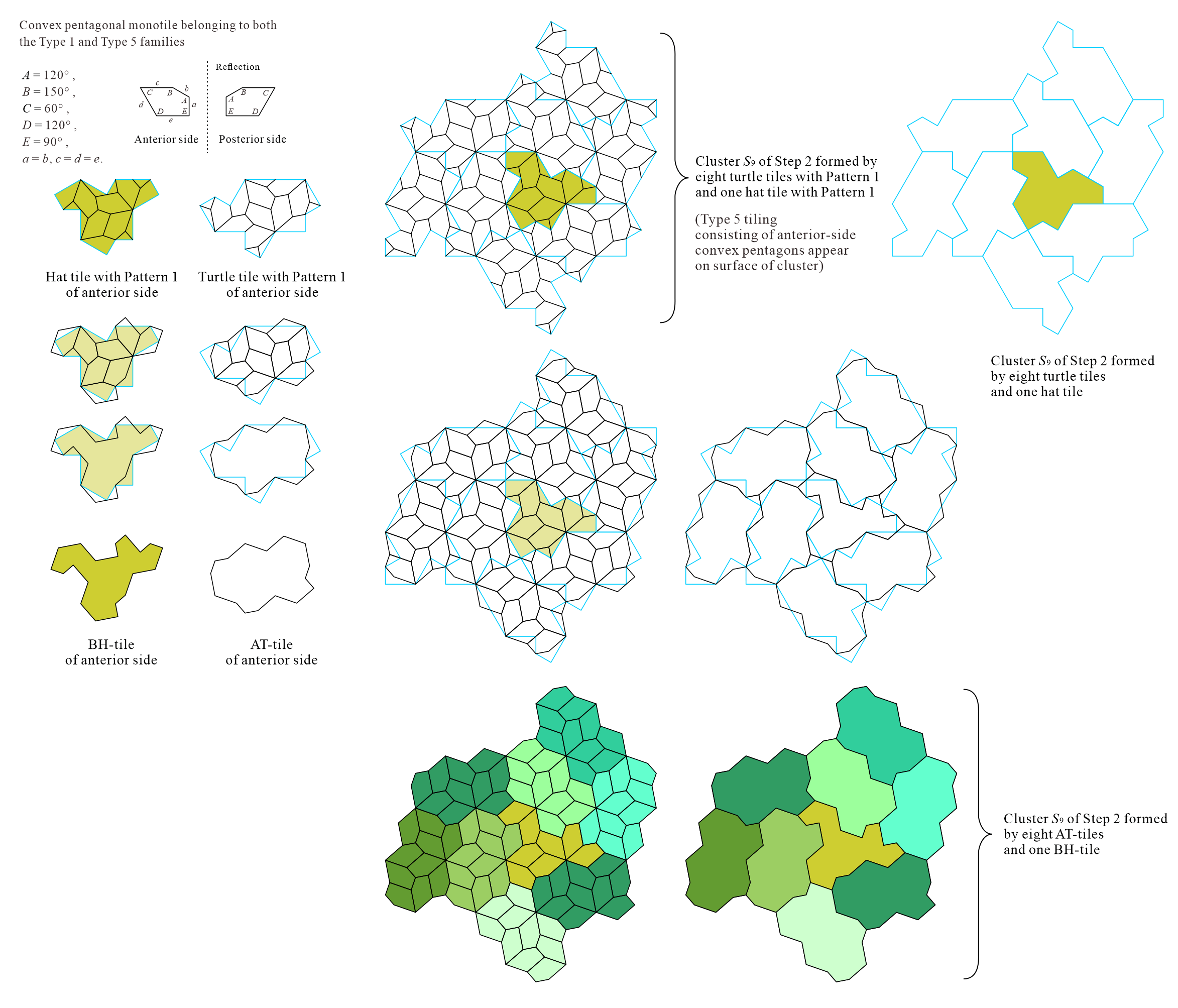} 
  \caption{{\small Step-2 cluster $S_9$ formed by eight turtle tiles with Pattern 1 from Figure~\ref{Fig.2-5} 
  and one hat tile with Pattern 1 from Figure~\ref{Fig.2-1}.} }
\label{Fig.3-6}
\end{figure}

\begin{figure}[H]
 \centering\includegraphics[width=15cm,clip]{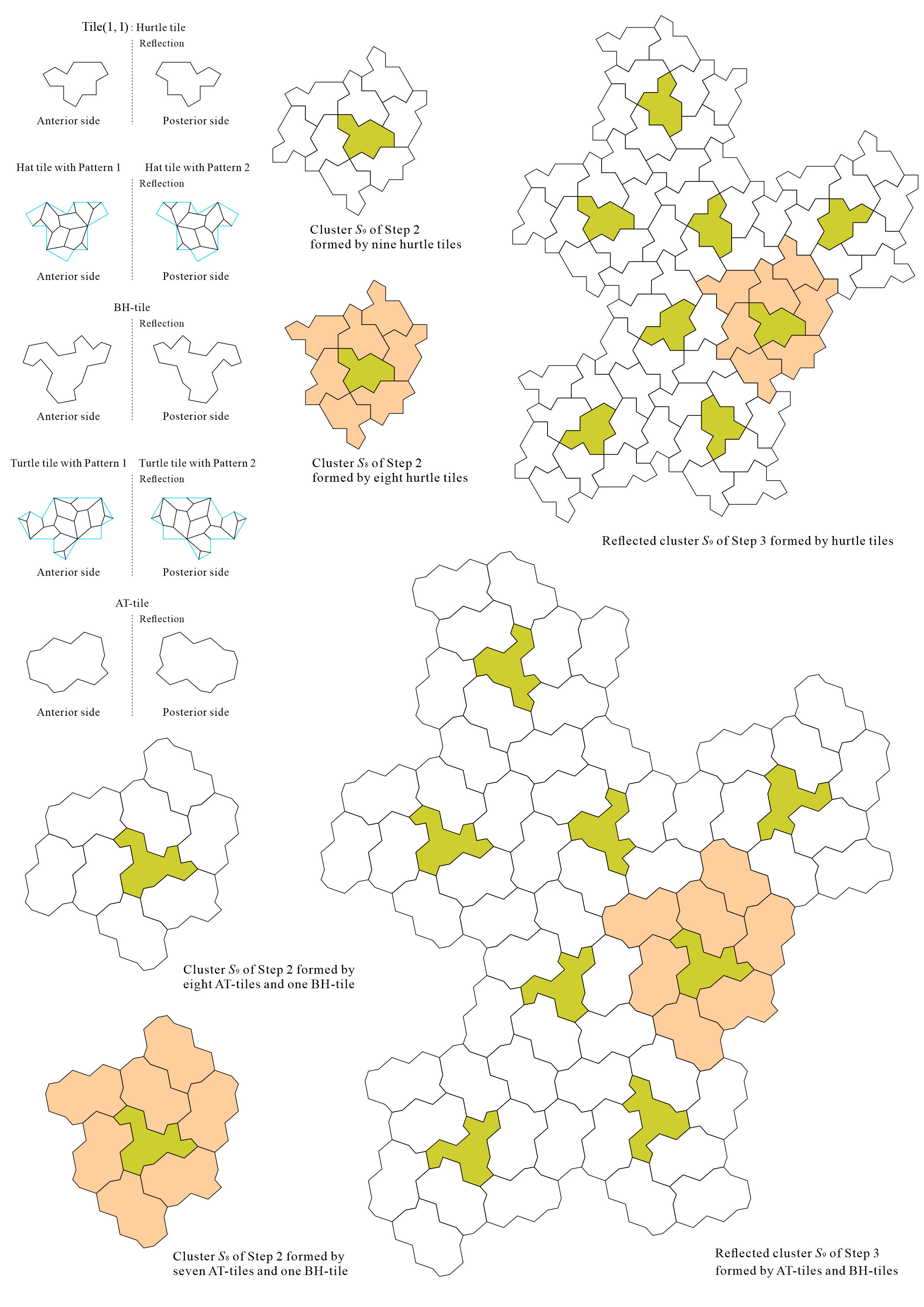} 
  \caption{{\small Application of the $S_8/S_9$ substitution to cluster $S_9$ formed by eight AT-tiles and one BH-tile.} }
\label{Fig.3-7}
\end{figure}

Figure~\ref{Fig.3-6} shows the Step-2 cluster $S_9$, formed by eight turtle tiles and one hat tile 
(corresponding to an odd tile), with the hat tile with Pattern 1 from Figure~\ref{Fig.2-1} and the turtle 
tile with Pattern 1 from Figure~\ref{Fig.2-5} applied. Owing to the properties of Pattern 1, a Type 5 tiling 
consisting of anterior-side convex pentagons appears on the surface of $S_9$. As a result, the cluster 
$S_9$ formed by hat and turtle tiles with Pattern 1 corresponds exactly to the cluster $S_9$ formed by 
eight AT-tiles and one BH-tile. Because the hat and turtle tiles forming $S_9$ in Figure~\ref{Fig.3-6} 
are on the anterior side, the AT-tiles and BH-tiles forming the corresponding $S_9$ are also on 
their anterior side. Therefore, the $S_8/S_9$ substitution can be applied to the cluster $S_9$ 
formed by eight AT-tiles and one BH-tile. This produces the next-step cluster $S_9$ as shown 
in Figure~\ref{Fig.3-7}. Repeated application of the substitution generates $\mathscr{T}_s$.

\subsection{$\mathscr{T}_s$ generated by hat and turtle tiles with Pattern 2}
\label{subsection3.3}

\begin{figure}[H]
 \centering\includegraphics[width=15cm,clip]{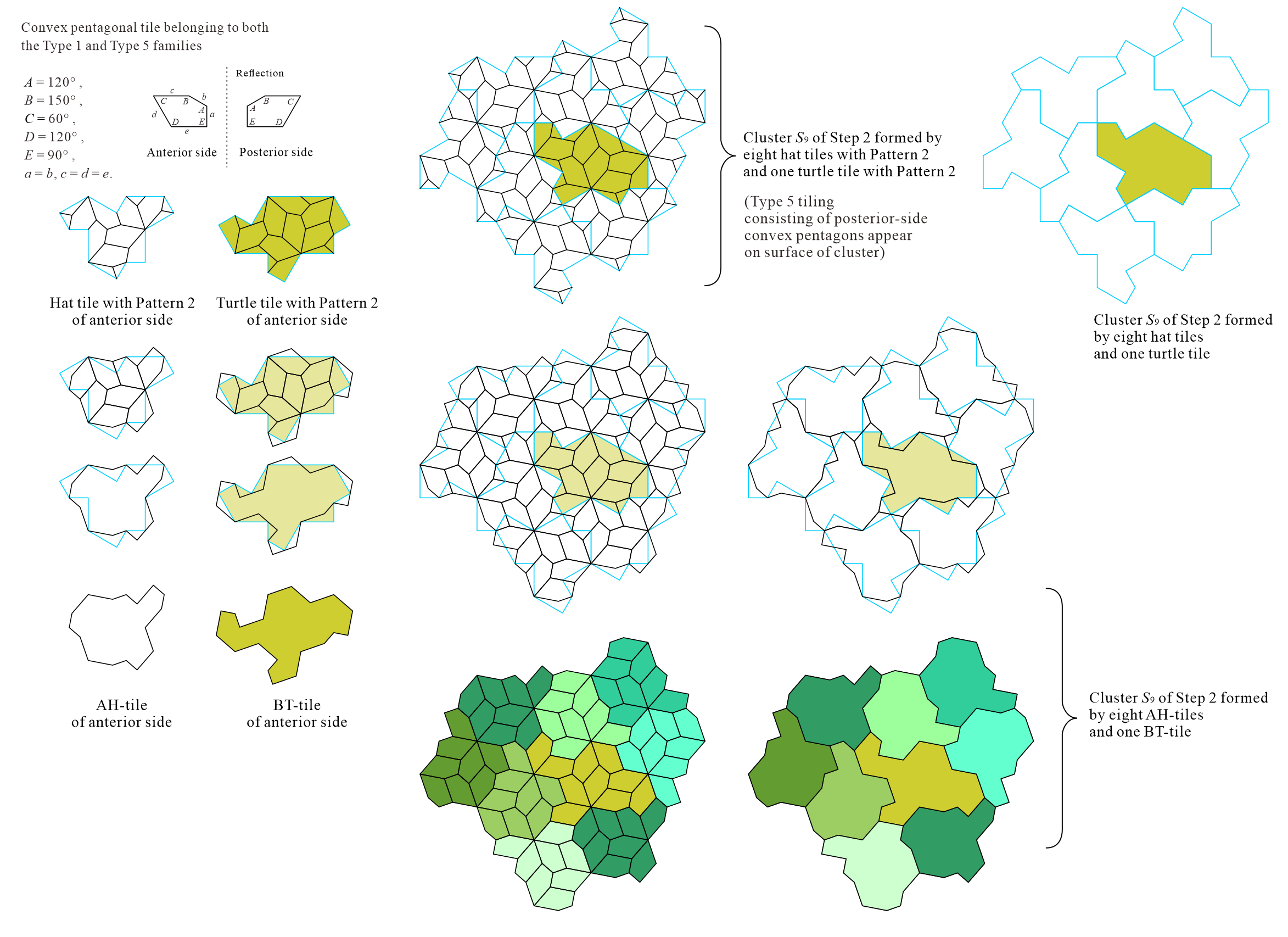} 
  \caption{{\small Step-2 cluster $S_9$ formed by eight hat tiles with Pattern 2 from Figure~\ref{Fig.2-2} 
  and one turtle tile with Pattern 2 from Figure~\ref{Fig.2-6}.} }
\label{Fig.3-8}
\end{figure}

\begin{figure}[H]
 \centering\includegraphics[width=15cm,clip]{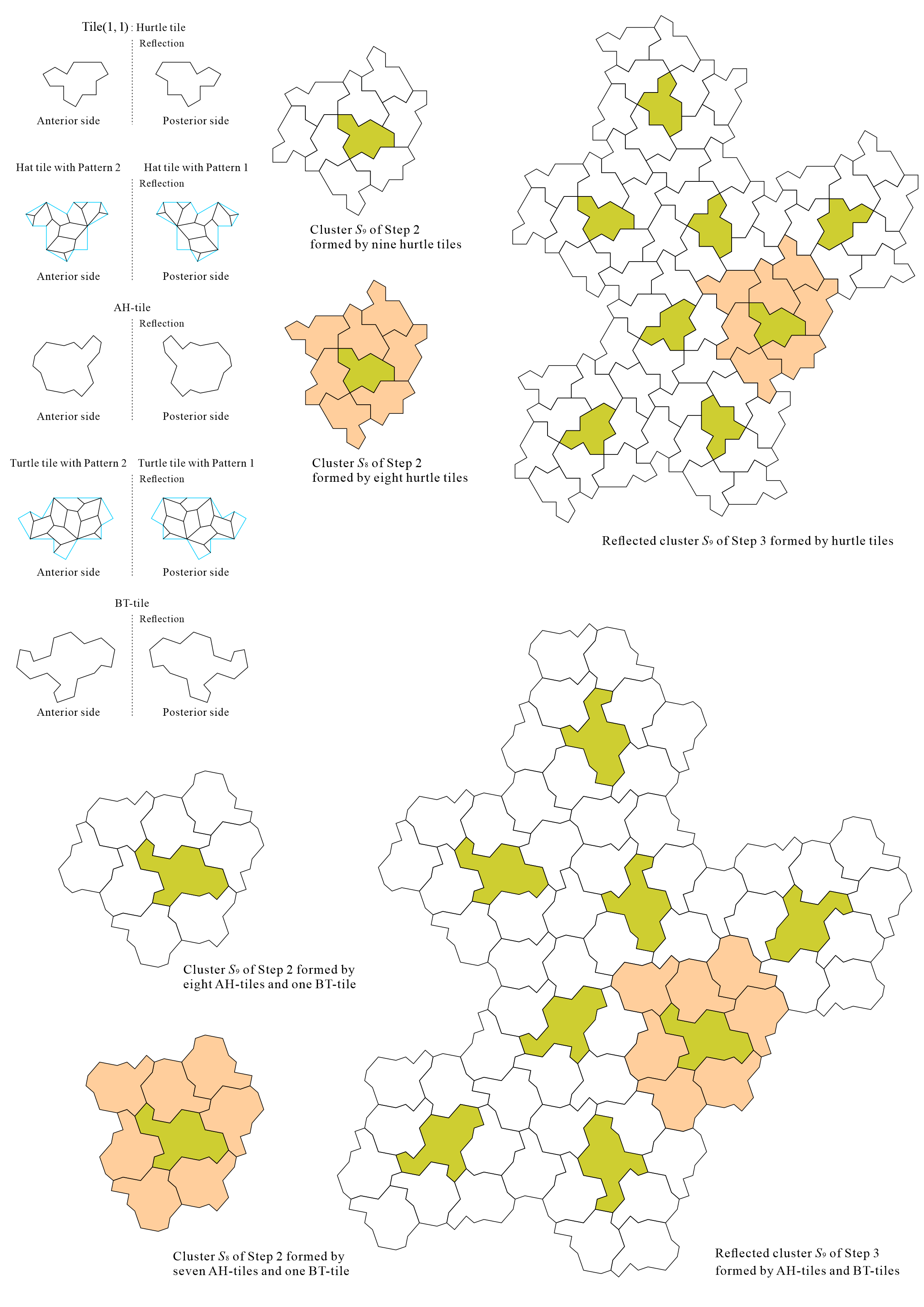} 
  \caption{{\small Application of the $S_8/S_9$ substitution to cluster $S_9$ formed by eight AH-tiles and one BT-tile.} }
\label{Fig.3-9}
\end{figure}

Figure~\ref{Fig.3-8} shows the Step 2 cluster $S_9$, formed by eight hat tiles and one turtle tile 
(corresponding to an odd tile), with the hat tile with Pattern 2 from Figure~\ref{Fig.2-2} and the turtle 
tile with Pattern 2 from Figure~\ref{Fig.2-6} applied. Owing to the properties of Pattern 2, a Type 5 tiling 
consisting of posterior-side convex pentagons appears on the surface of $S_9$. As a result, the cluster 
$S_9$ formed by hat and turtle tiles with Pattern 2 corresponds exactly to the cluster $S_9$ formed by 
eight AH-tiles and one BT-tile. Because the hat and turtle tiles forming $S_9$ in Figure~\ref{Fig.3-8} 
are on the anterior side, the AH-tiles and BT-tiles forming the corresponding $S_9$ are also on 
their anterior side. Therefore, the $S_8/S_9$ substitution can be applied to the cluster 
$S_9$ formed by eight AH-tiles and one BT-tile. This produces the next-step cluster $S_9$ as shown 
in Figure~\ref{Fig.3-9}. Repeated application of the substitution generates $\mathscr{T}_s$.

\begin{figure}[H]
 \centering\includegraphics[width=15cm,clip]{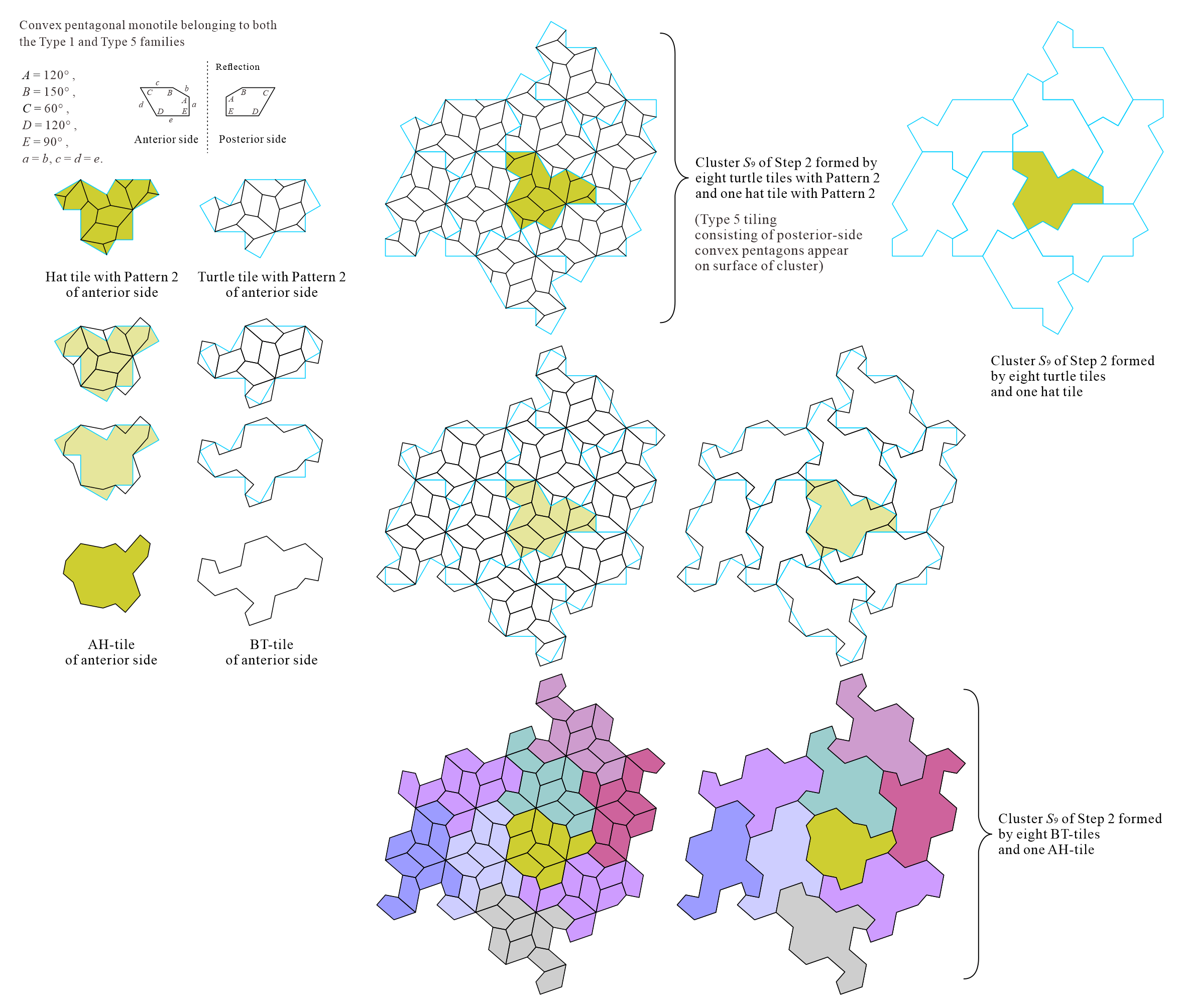} 
  \caption{{\small Step-2 cluster $S_9$ formed by eight turtle tiles with Pattern 2 from Figure~\ref{Fig.2-6} 
  and one hat tile with Pattern 2 from Figure~\ref{Fig.2-2}.} }
\label{Fig.3-10}
\end{figure}

\begin{figure}[H]
 \centering\includegraphics[width=15cm,clip]{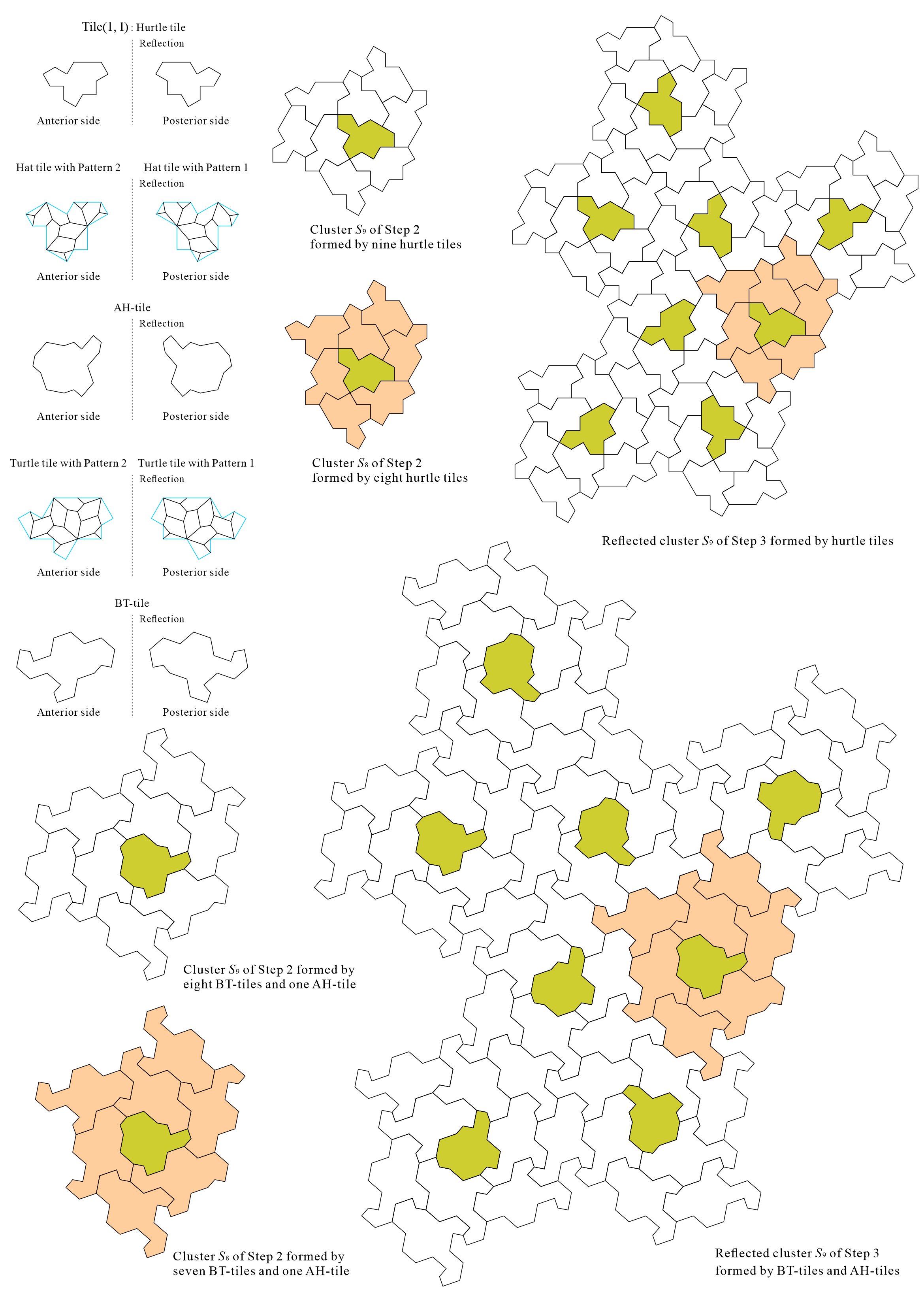} 
  \caption{{\small Application of the $S_8/S_9$ substitution to cluster $S_9$ formed by eight BT-tiles and one AH-tile.} }
\label{Fig.3-11}
\end{figure}

Figure~\ref{Fig.3-10} shows the Step 2 cluster $S_9$, formed by eight turtle tiles and one hat tile 
(corresponding to an odd tile), with the hat tile with Pattern 2 from Figure~\ref{Fig.2-2} and the turtle 
tile with Pattern 2 from Figure~\ref{Fig.2-6} applied. Owing to the properties of Pattern 2, a Type 5 tiling 
consisting of posterior-side convex pentagons appears on the surface of $S_9$. As a result, the cluster
 $S_9$ formed by hat and turtle tiles with Pattern 2 corresponds exactly to the cluster $S_9$ formed by 
 eight BT-tiles and one AH-tile. Because the hat and turtle tiles forming $S_9$ in Figure~\ref{Fig.3-10} 
are on the anterior side, the BT-tiles and AH-tiles forming the corresponding $S_9$ are also on 
their anterior side. Therefore, the $S_8/S_9$ substitution can be applied to the cluster
  $S_9$ formed by eight BT-tiles and one AH-tile. This produces the next-step cluster $S_9$ as shown 
  in Figure~\ref{Fig.3-11}. Repeated application of the substitution generates $\mathscr{T}_s$.

\subsection{Reversal of the densities of two types of tiles in $\mathscr{T}_s$}
\label{subsection3.4}

Let $S_9(n)$ denote the cluster corresponding to $S_9$ obtained after $n$ steps of the $S_8/S_9$ 
substitution. Accordingly, the Step-2 cluster $S_9$ shown in Figure~\ref{Fig.3-1}(b) is denoted by $S_9(2)$.

The cluster $S_9(2)$ formed by eight hat tiles and one turtle tile in Figure~\ref{Fig.3-2} is denoted by 
\mbox{$S_9(8\text{H}+\text{T}, 0)$}. The same notation is used for the other cases. Therefore, the cluster $S_9(2)$ 
formed by  eight turtle tiles and one hat tile in Figure~\ref{Fig.3-3} is denoted by \mbox{$S_9(8\text{T}+\text{H}, 0)$}, 
$S_9(2)$ in Figure~\ref{Fig.3-4} is denoted by \mbox{$S_9(8\text{BH}+\text{AT}, 0)$}, $S_9(2)$ in 
Figure~\ref{Fig.3-6} is denoted by \mbox{$S_9(8\text{AT}+\text{BH}, 0)$, $S_9(2)$} in Figure~\ref{Fig.3-8} is 
denoted by \mbox{$S_9(8\text{AH}+\text{BT}, 0)$}, and $S_9(2)$ in Figure~\ref{Fig.3-10} is denoted by 
\mbox{$S_9(8\text{BT}+\text{AH}, 0)$}. \mbox{$S_9(8\text{H}+\text{T}, m)$} denotes the cluster corresponding 
to $S_9$ obtained after $m$ steps of the $S_8/S_9$ substitution starting from \mbox{$S_9(8\text{H}+\text{T}, 0)$} 
(i.e., \mbox{$S_9(8\text{H}+\text{T}, 1)$} corresponds to $S_9(3)$); the same notation is used for the other cases. 
Accordingly, \mbox{$S_9(8\text{H}+\text{T}, \infty)$} and \mbox{$S_9(8\text{T}+\text{H}, \infty)$} denote 
$\mathscr{T}_s$ generated by hat tiles and turtle tiles, respectively. Similarly, \mbox{$S_9(8\text{BH}+\text{AT}, \infty)$} 
and \mbox{$S_9(8\text{AT}+\text{BH}, \infty)$} denote $\mathscr{T}_s$ generated by BH-tile and AT-tile, 
respectively, and \mbox{$S_9(8\text{AH}+\text{BT}, \infty)$} and \mbox{$S_9(8\text{BT}+\text{AH}, \infty)$} 
denote $\mathscr{T}_s$ generated by AH-tile and BT-tile, respectively.

The $S_8/S_9$ substitution can be represented by the matrix 
\[
\begin{pmatrix}
1 & 1 \\
6 & 7
\end{pmatrix}
\]
which has the eigenvalues $4 \pm \sqrt{15}$. The eigenvector corresponding to the $S_8/S_9$ substitution is 
\[
\begin{pmatrix}
\sqrt{15}-3 \\
6
\end{pmatrix}
\]
Therefore, the ratio of clusters $S_8$ to $S_9$ approaches $\sqrt{15}-3 : 6$.

Each cluster has one odd tile and seven or eight even tiles, which means the density of the odd tiles in $\mathscr{T}_s$ is
\[
\frac{(\sqrt{15}-3)+6}{8(\sqrt{15}-3)+9\times 6}
\]
which is equal to $\frac{5-\sqrt{15}}{10} \approx 0.112701\ldots$.

Therefore, the following statements hold:
\begin{itemize}
\item In \mbox{$S_9(8\text{H}+\text{T}, \infty)$}, the density of turtle tiles is $\frac{5-\sqrt{15}}{10} \approx 0.112701\ldots$.
\item In \mbox{$S_9(8\text{T}+\text{H}, \infty)$}, the density of hat tiles is $\frac{5-\sqrt{15}}{10} \approx 0.112701\ldots$.
\item In \mbox{$S_9(8\text{BH}+\text{AT}, \infty)$}, the density of AT-tiles is $\frac{5-\sqrt{15}}{10} \approx 0.112701\ldots$.
\item In \mbox{$S_9(8\text{AT}+\text{BH}, \infty)$}, the density of BH-tiles is $\frac{5-\sqrt{15}}{10} \approx 0.112701\ldots$.
\item In \mbox{$S_9(8\text{AH}+\text{BT}, \infty)$}, the density of BT-tiles is $\frac{5-\sqrt{15}}{10} \approx 0.112701\ldots$.
\item In \mbox{$S_9(8\text{BT}+\text{AH}, \infty)$}, the density of AH-tiles is $\frac{5-\sqrt{15}}{10} \approx 0.112701\ldots$.
\end{itemize}

The tiling $\mathscr{T}_s$ generated by \mbox{Tile$(1, 1)$} is monohedral. \mbox{$S_9(8\text{H}+\text{T}, \infty)$}, 
\mbox{$S_9(8\text{T}+\text{H}, \infty)$}, \mbox{$S_9(8\text{BH}+\text{AT}, \infty)$}, 
\mbox{$S_9(8\text{AT}+\text{BH}, \infty)$}, \mbox{$S_9(8\text{AH}+\text{BT}, \infty)$}, and 
\mbox{$S_9(8\text{BT}+\text{AH}, \infty)$} are dihedral and are regarded as distinct patterns.

\section{Tile set consisting of hat and turtle tiles and periodic tilings}
\label{section4}

When there are no restrictions on tiling generation\footnote{ 
In the Euclidean plane, a tiling is generated by combining tiles obtained by applying isometries such 
as translations, rotations, and reflections. The phrase ``when there are no restrictions on tiling generation'' 
means that translations, rotations, and reflections are all allowed to be applied to the tiles during tiling 
generation. On the other hand, for example, a case in which the use of reflected tiles is not allowed 
during tiling generation can be considered, and this corresponds to a case in which restrictions are imposed. 
A tiling using a tile set consisting of hat tiles and turtle tiles of the same orientation, with their sizes 
adjusted so that their shortest edges have the same length, can be regarded as equivalent to a tiling 
by \mbox{Tile$(1, 1)$}, and when the use of reflected tiles is not allowed during tiling generation, 
they admit only non-periodic tilings \cite{ref6}. Unless otherwise noted, the tilings considered in this study 
have no restrictions on tiling generation; however, in light of the results in \cite{ref6}, we explicitly 
emphasize this point here. Note that the periodic tilings generated by the tile set consisting of hat and 
turtle tiles are dihedral.
}, a tile set consisting of two types of tiles (concave polygons), namely hat and turtle tiles whose sizes 
are adjusted so that their shortest edges have the same length, can generate periodic tilings. 
Therefore, this section introduces such tilings.

\begin{figure}[H]
 \centering\includegraphics[width=15cm,clip]{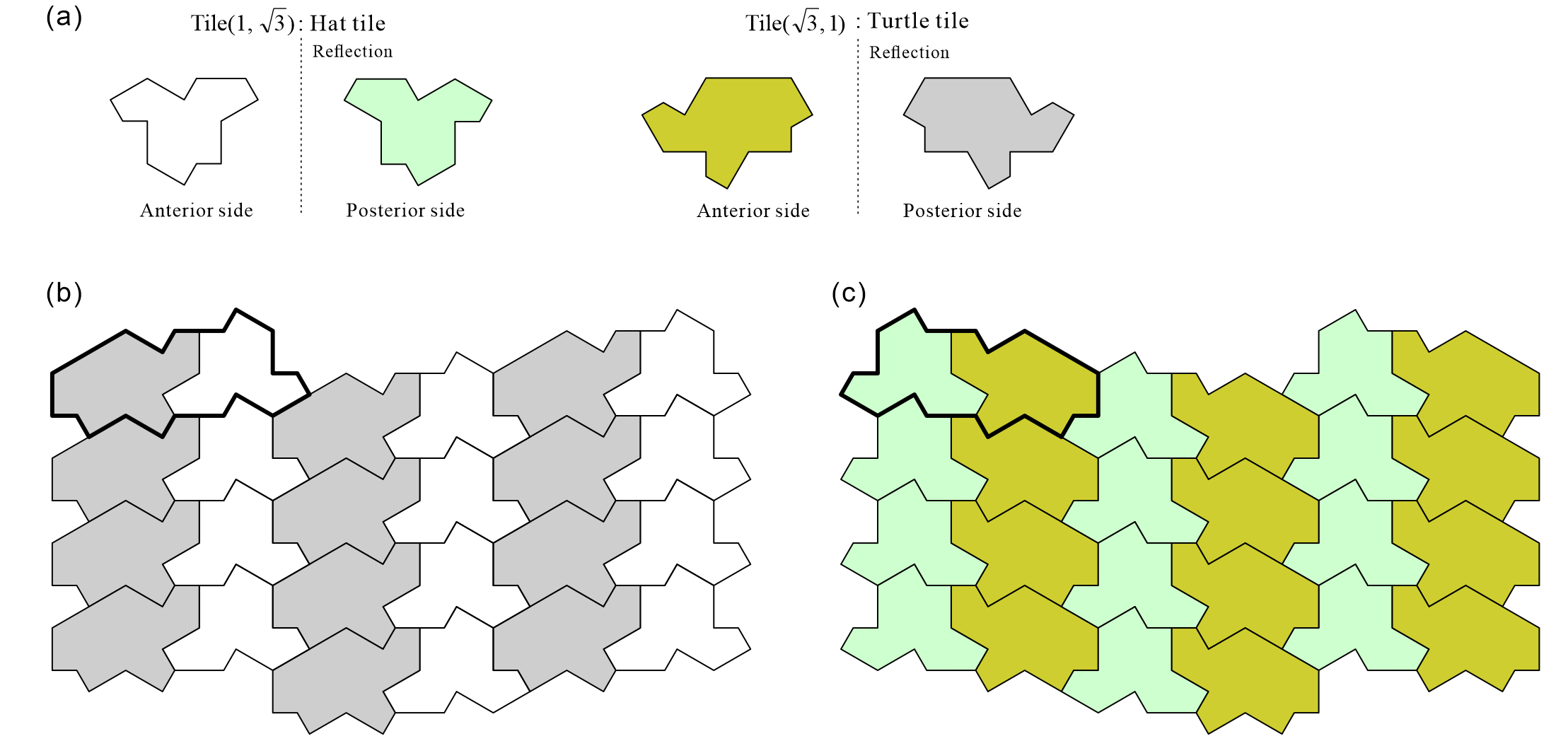} 
  \caption{{\small Periodic tilings in which the ratio of hat tiles to turtle tiles is $1:1$.} }
\label{Fig.4-1}
\end{figure}

Figure~\ref{Fig.4-1} shows periodic tilings in which the ratio of hat tiles to turtle tiles is $1:1$ 
(see also \ref{appA}). When the anterior and posterior sides of the hat and turtle tiles are defined as 
shown in Figure~\ref{Fig.4-1}(a), the periodic tiling shown in Figure~\ref{Fig.4-1}(b) can be generated 
using a translation unit (a unit that generates a periodic tiling only by translations) constructed 
from a pair consisting of an anterior-side hat tile and a posterior-side turtle tile. Similarly, 
the periodic tiling shown in Figure~\ref{Fig.4-1}(c) can be generated using a translation unit 
constructed from a pair consisting of a posterior-side hat tile and an anterior-side turtle tile. 
When two tilings, such as those shown in Figures~\ref{Fig.4-1}(b) and \ref{Fig.4-1}(c), are 
congruent to each other by reflection, they are generally regarded as identical in this study, 
and only one of them is usually shown. However, both tilings are intentionally shown in 
Figure~\ref{Fig.4-1} because this is related to the results presented in Section~\ref{subsection5.2}.

\vspace{2\baselineskip}

\begin{figure}[H]
 \centering\includegraphics[width=15cm,clip]{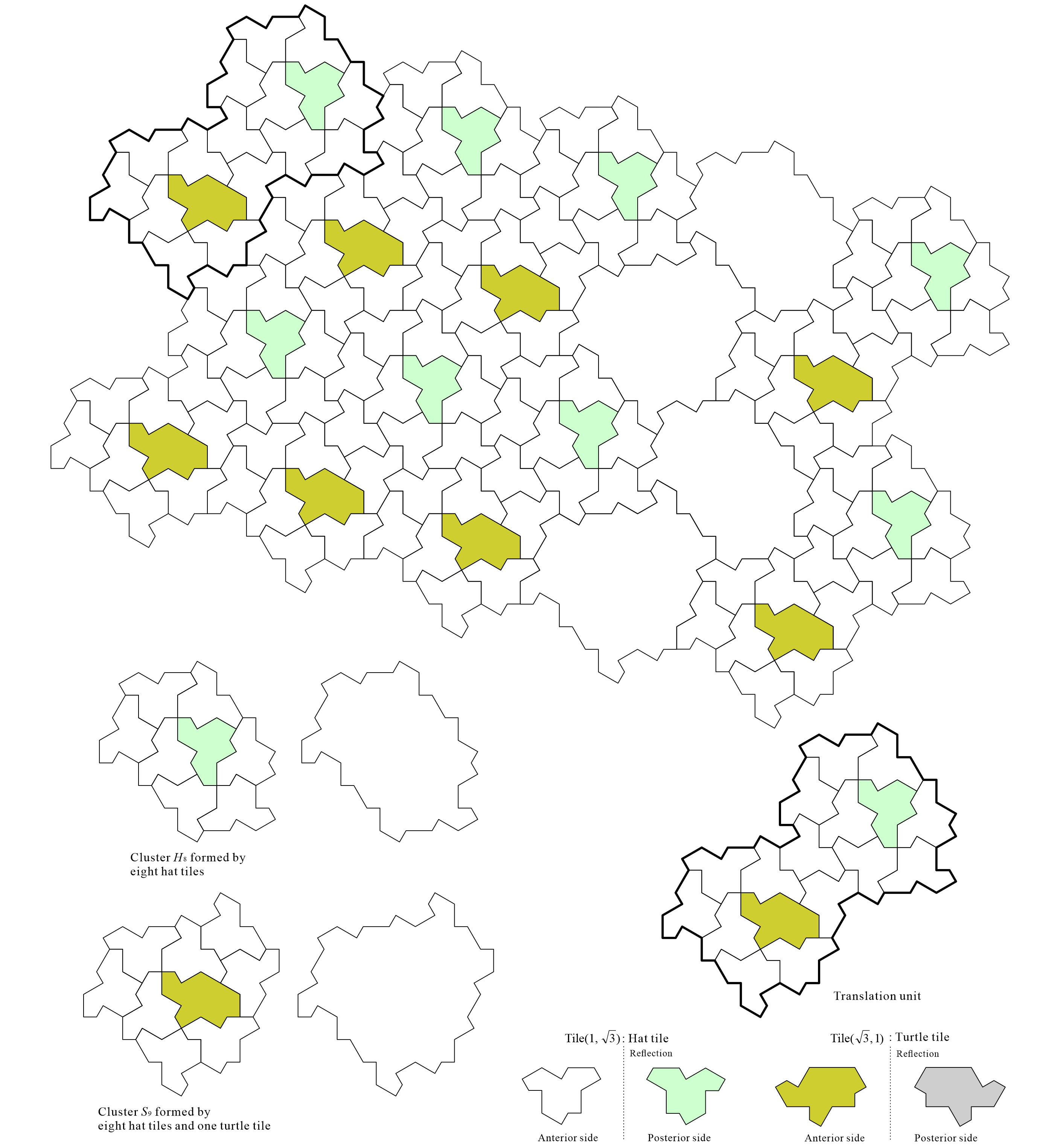} 
  \caption{{\small Periodic tiling generated using a translation unit constructed from cluster $H_8$ formed 
  by eight hat tiles and cluster $S_9$ formed by eight hat tiles and one turtle tile.} }
\label{Fig.4-2}
\end{figure}

Figure~\ref{Fig.4-2} shows a periodic tiling generated using a translation unit constructed from cluster 
$H_8$ formed by eight hat tiles in Figure~\ref{Fig.2-1} and cluster $S_9$ formed by eight hat tiles 
and one turtle tile in Figure~\ref{Fig.3-2}. The ratio of hat tiles to turtle tiles within the translation unit 
is $16:1$ (therefore, the ratio of hat tiles to turtle tiles in this periodic tiling is also $16:1$), and the 
ratio of anterior-side tiles to posterior-side tiles is $16:1$ (there is one posterior-side hat tile in $H_8$, 
while all other tiles are anterior-side hat tiles and turtle tiles).

\vspace{2\baselineskip}

\begin{figure}[H]
 \centering\includegraphics[width=15cm,clip]{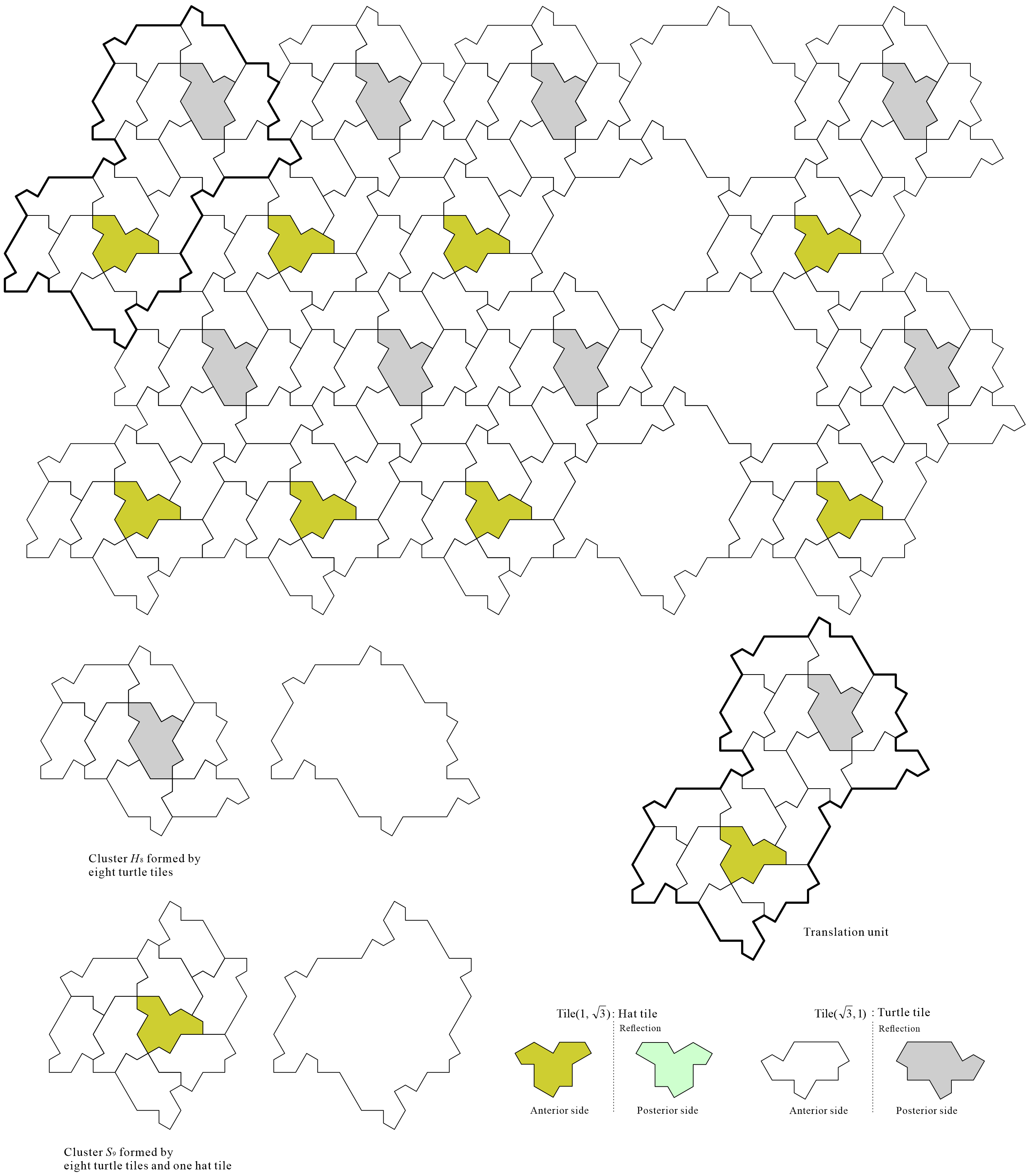} 
  \caption{{\small Periodic tiling generated using a translation unit constructed from cluster $H_8$ formed 
  by eight turtle tiles and cluster $S_9$ formed by eight turtle tiles and one hat tile.} }
\label{Fig.4-3}
\end{figure}

Figure~\ref{Fig.4-3} shows a periodic tiling generated using a translation unit constructed from cluster 
$H_8$ formed by eight turtle tiles in Figure~\ref{Fig.2-5} and cluster $S_9$ formed by eight turtle 
tiles and one hat tile in Figure~\ref{Fig.3-3}. The ratio of turtle tiles to hat tiles within the translation 
unit is $16:1$, and the ratio of anterior-side tiles to posterior-side tiles is $16:1$ (there is one 
posterior-side turtle tile in $H_8$).

\vspace{2\baselineskip}

\begin{figure}[H]
 \centering\includegraphics[width=14cm,clip]{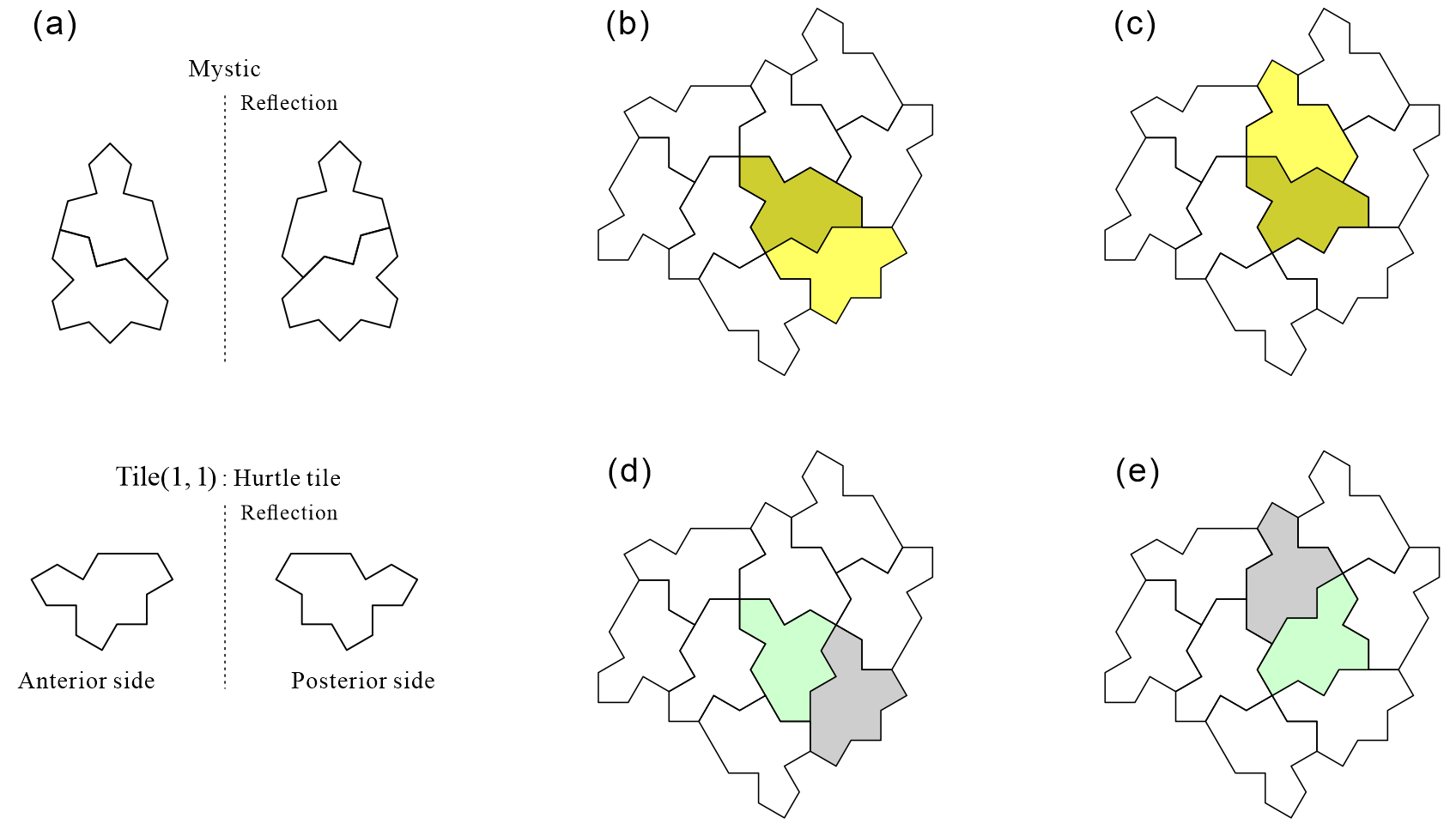} 
  \caption{{\small Cluster $S_9$ formed by \mbox{Tile$(1, 1)$} and the Mystic structure.} }
\label{Fig.4-4}
\end{figure}

\begin{figure}[H]
 \centering\includegraphics[width=14cm,clip]{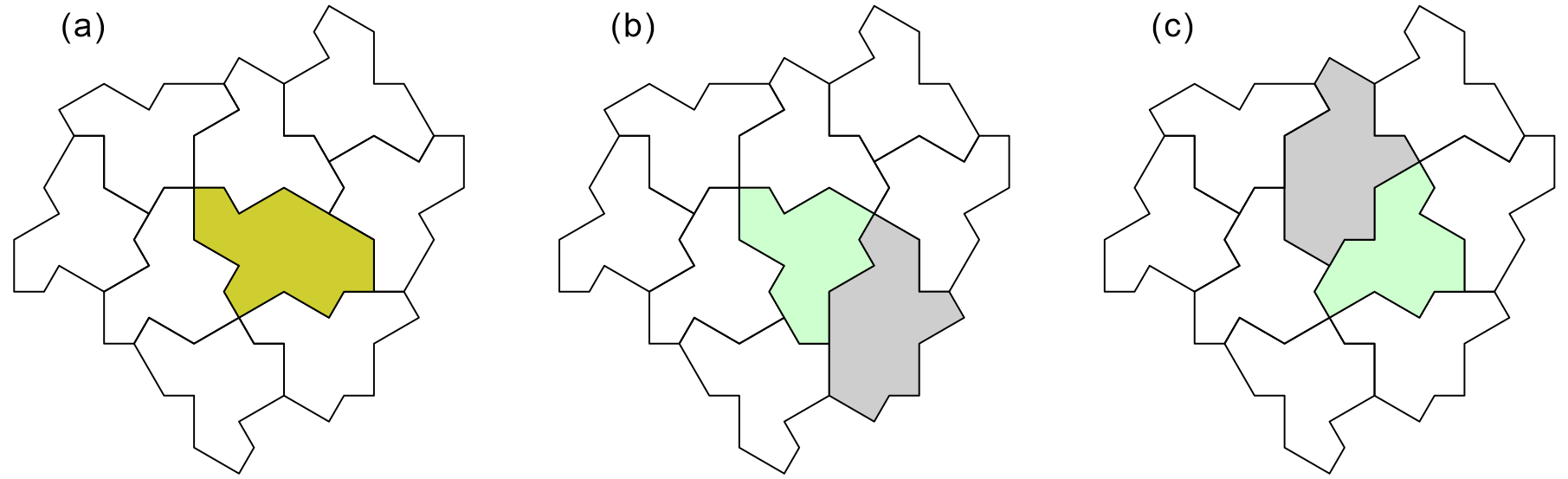} 
  \caption{{\small Original $S_9$ and $S_9$ with the replaced Mystic regions.} }
\label{Fig.4-5}
\end{figure}

Figure~\ref{Fig.4-4}(a) shows a structure formed by two \mbox{Tile$(1, 1)$}, which was named 
``Mystic'' in \cite{ref6}. As shown in Figures~\ref{Fig.4-4}(b) and (c), two regions corresponding 
to the Mystic structure  exist within the Step-2 cluster $S_9$ formed by nine \mbox{Tile$(1, 1)$} 
(i.e., there are two Mystic regions formed by three \mbox{Tile$(1, 1)$}). Because the outer boundary 
of the Mystic has line symmetry, as shown in Figure~\ref{Fig.4-4}(a), the Mystic regions within $S_9$ 
shown in Figures~\ref{Fig.4-4}(b) and (c) can be replaced with reflected versions of the Mystic (i.e., 
a Mystic formed by posterior-side \mbox{Tile$(1, 1)$}) (see Figures~\ref{Fig.4-4}(d) and (e)). Such replacement 
of the Mystic regions within $S_9$ is also possible for a cluster $S_9$ formed by eight hat tiles and one 
turtle tile. Figure~\ref{Fig.4-5} illustrates this property (Figure~\ref{Fig.4-5}(a) shows the original $S_9$, 
while Figures~\ref{Fig.4-5}(b) and (c) show $S_9$ after replacement of the Mystic regions\footnote{  
The Step-2 cluster $S_9$ formed by \mbox{Tile$(1, 1)$} refers to a cluster formed using 
\mbox{Tile$(1, 1)$} that all have the same orientation. Therefore, the clusters shown in 
Figures~\ref{Fig.4-5}(b) and (c), obtained by replacing the Mystic regions with reflected structures, 
are not $S_9$. However, in this manuscript, such clusters are also referred to as $S_9$ unless 
doing so causes confusion.
}). A similar replacement is possible for cluster $S_9$ formed by eight turtle tiles and one hat tile.

\begin{figure}[H]
 \centering\includegraphics[width=15cm,clip]{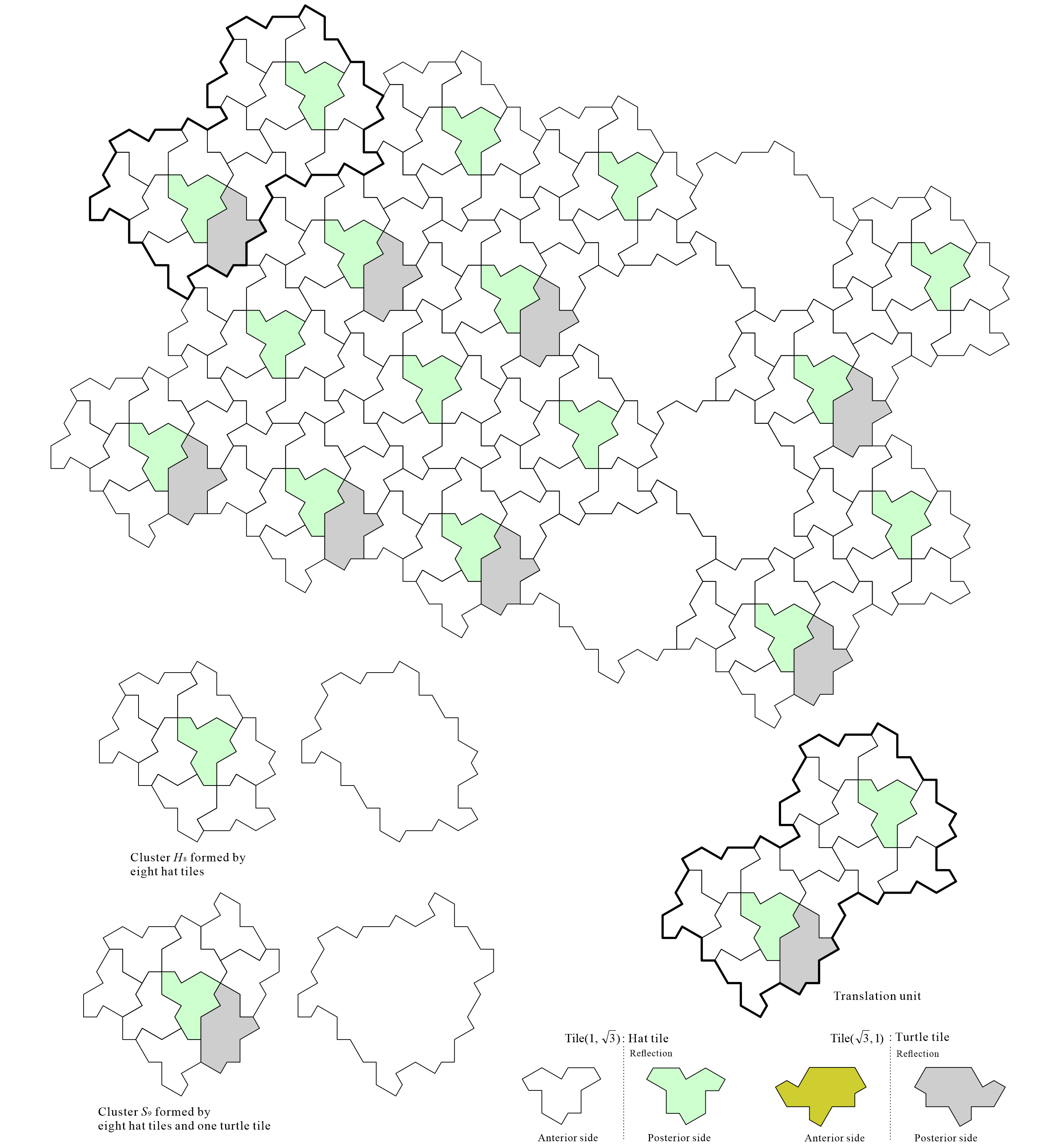} 
  \caption{{\small Periodic tiling generated using a translation unit constructed from $S_9$ in 
  Figure~\ref{Fig.4-5}(b) and $H_8$ formed by eight hat tiles.} }
\label{Fig.4-6}
\end{figure}

The outer boundaries of $S_9$ in Figure~\ref{Fig.4-5}(a) and those of $S_9$ in Figures~\ref{Fig.4-5}(b) 
and (c) are identical. Because the translation unit in Figure~\ref{Fig.4-2} can be regarded as being constructed 
using $S_9$ in Figure~\ref{Fig.4-5}(a) and $H_8$ formed by eight hat tiles, a translation unit with the same 
outer boundary can be constructed using $S_9$ in Figure~\ref{Fig.4-5}(b) and $H_8$ formed by eight hat tiles, 
thereby generating a periodic tiling as shown in Figure~\ref{Fig.4-6}. The ratio of hat tiles to turtle tiles in the 
translation unit of Figure~\ref{Fig.4-6} is $16:1$, and the ratio of anterior-side tiles to posterior-side tiles is $14:3$ 
(there is one posterior-side hat tile in each of $H_8$ and $S_9$, and one posterior-side turtle tile in $S_9$).

\begin{figure}[H]
 \centering\includegraphics[width=15cm,clip]{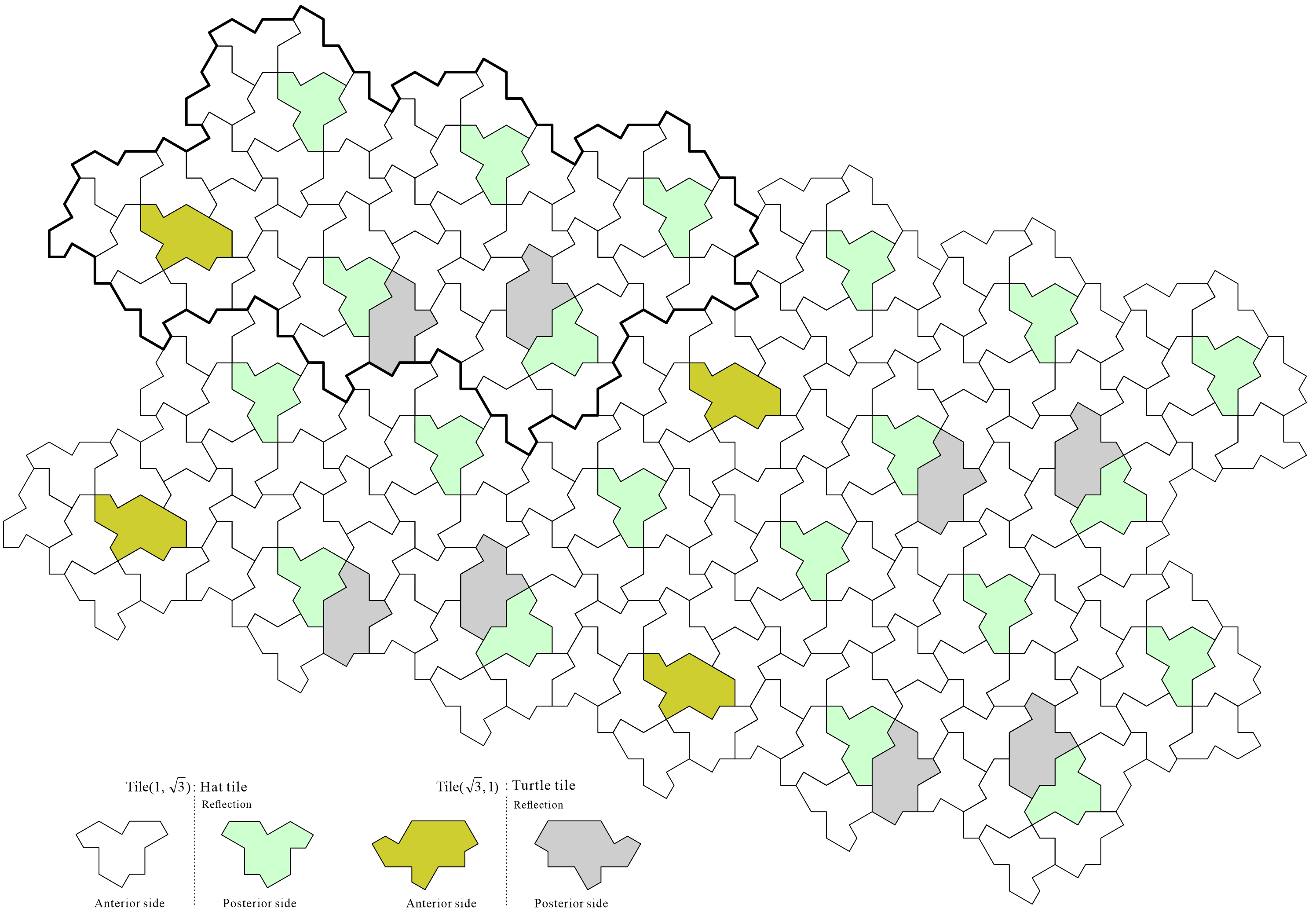} 
  \caption{{\small Periodic tiling using a translation unit constructed from $S_9$ in Figures~\ref{Fig.4-5}(a)--(c) 
  and cluster $H_8$ formed by eight hat tiles.} }
\label{Fig.4-7}
\end{figure}

Similarly, a translation unit using $S_9$ in Figure~\ref{Fig.4-5}(c) and cluster $H_8$ formed by eight 
hat tiles can also be constructed. From the property that the replacement of the Mystic regions 
within $S_9$ described above can produce three types of clusters with different arrangements of 
internal tiles, it can be seen that various translation units can be constructed by combining these 
clusters (furthermore, it is also possible to intentionally avoid generating periodic tilings). For example, 
Figure~\ref{Fig.4-7} shows an example of a periodic tiling using a translation unit constructed from 
$S_9$ in Figures~\ref{Fig.4-5}(a)--(c) and cluster $H_8$ formed by eight hat tiles. It can be observed 
that the tiling contains both anterior-side and posterior-side hat tiles and both anterior-side and 
posterior-side turtle tiles.

Note that there are methods of combining $H_8$ and $S_9$ other than the translation unit shown in 
Figure~\ref{Fig.4-3}. Therefore, more diverse combinations are possible, and various patterns of tilings 
can be generated. See \ref{appA} for details.

Furthermore, the Type 5 tiling patterns can be associated with the periodic tilings introduced in this 
section using the same method as described in Section~\ref{section2}. See \ref{appB} for details.

\vspace{5\baselineskip}
\section{Properties of the AH-tile, BH-tile, AT-tile, and BT-tile}
\label{section5}

\subsection{Heesch number of the AH-tile, BH-tile, AT-tile, and BT-tile}
\label{subsection5.1}

Let $F$ be a tile in the plane. A corona of $F$ is the set consisting of the centrally placed tile $F$ 
and a layer of congruent copies surrounding it. A corona is formed without gaps or overlaps, except 
along tile boundaries. The first corona is the set of all tiles sharing a boundary point with a centrally 
placed tile (including the original tile itself). For an integer $h\ge2$, the $h$-th corona is the set of 
all tiles sharing a boundary point with the $(h-1)$-th corona (including the $(h-1)$-th corona itself). 
The maximum number of corona layers that can be formed around $F$ is called the Heesch number 
of $F$. Therefore, the Heesch number of a polygonal monotile is $\infty$ (infinite), whereas the 
Heesch number of a polygon that can generate neither a tiling nor a corona is 0 (zero) \cite{ref1, ref3, ref7}.

\begin{figure}[H]
 \centering\includegraphics[width=15cm,clip]{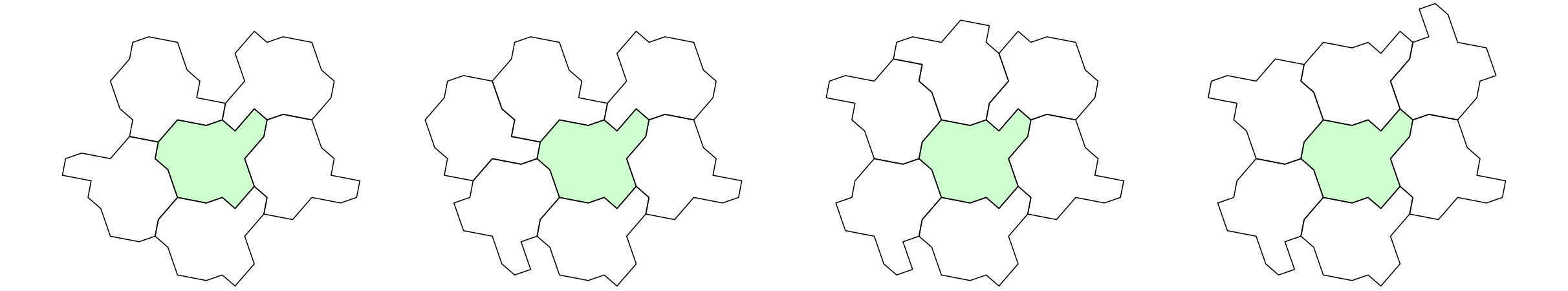} 
  \caption{{\small First corona that can be formed by AH-tile.} }
\label{Fig.5-1}
\end{figure}

\begin{figure}[H]
 \centering\includegraphics[width=15cm,clip]{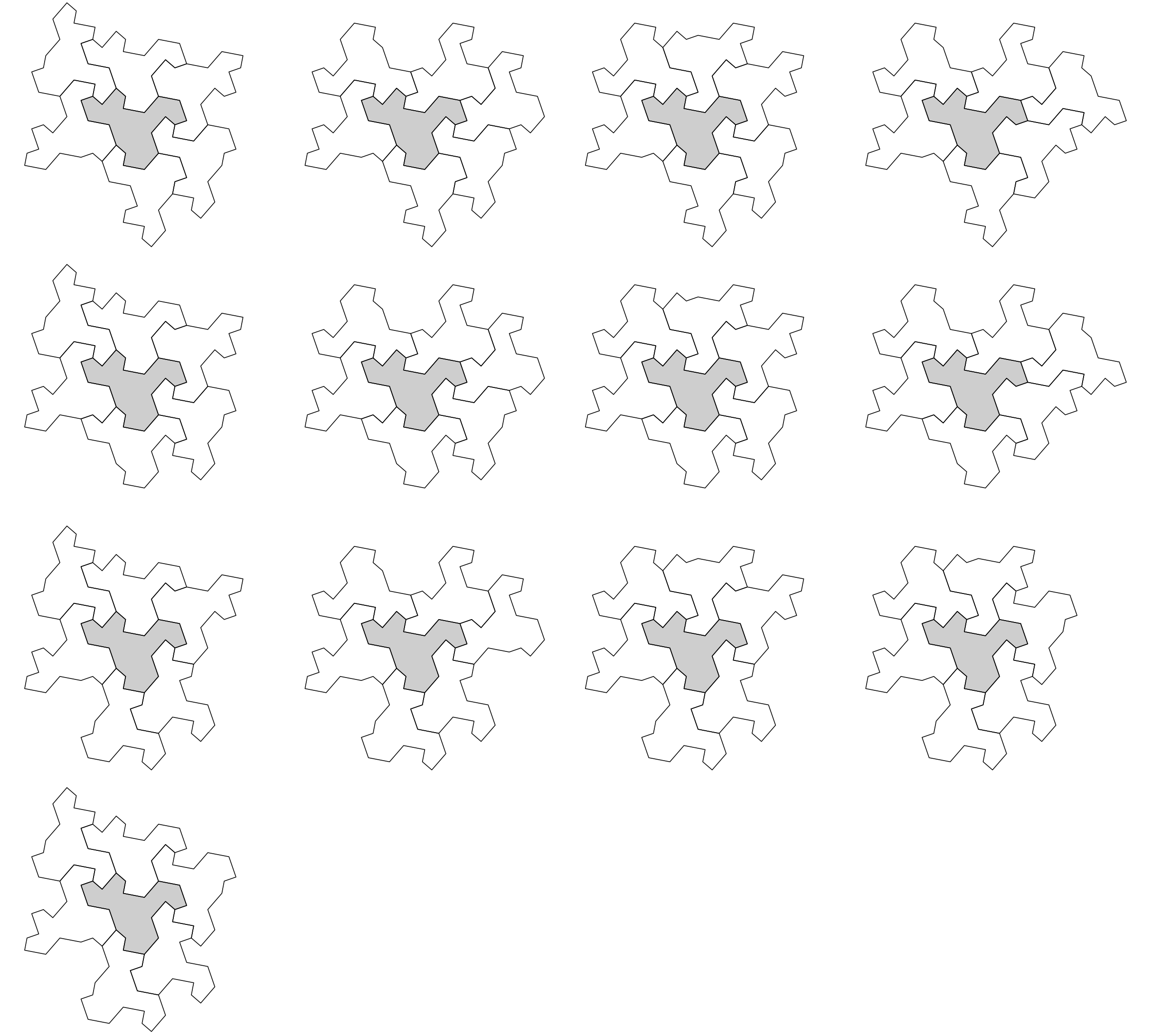} 
  \caption{{\small First corona that can be formed by BH-tile.} }
\label{Fig.5-2}
\end{figure}

\begin{figure}[H]
 \centering\includegraphics[width=15cm,clip]{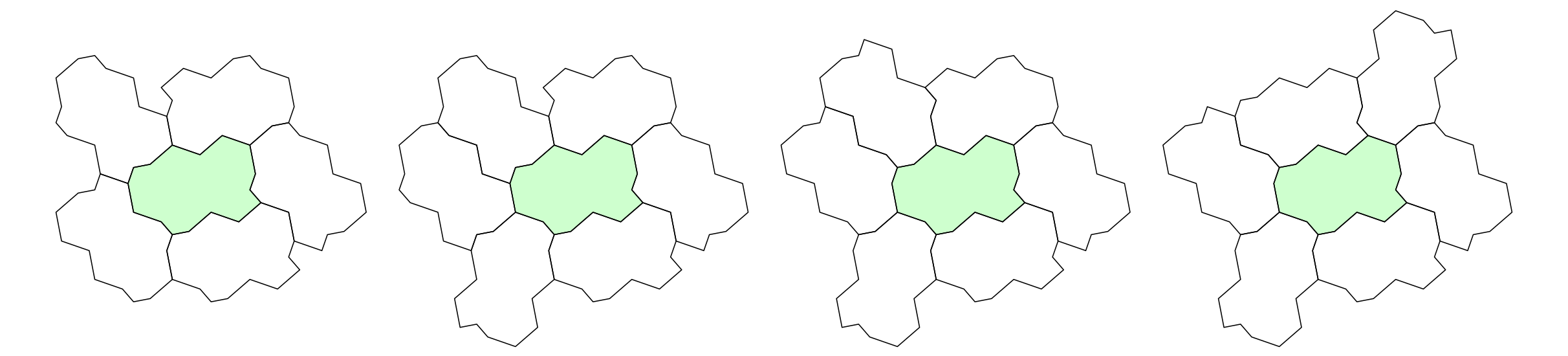} 
  \caption{{\small First corona that can be formed by AT-tile.} }
\label{Fig.5-3}
\end{figure}

\begin{figure}[H]
 \centering\includegraphics[width=15cm,clip]{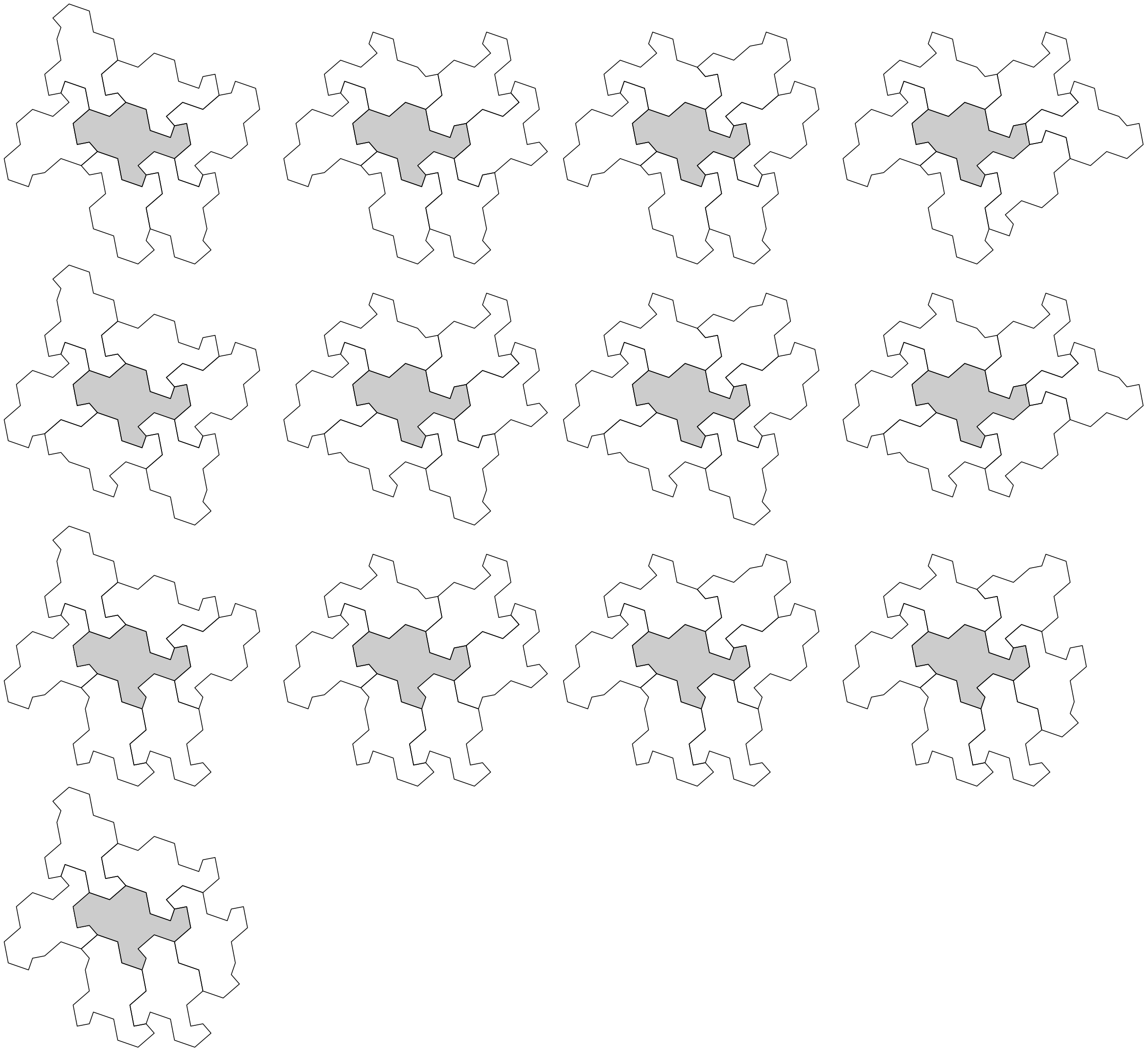} 
  \caption{{\small First corona that can be formed by BT-tile.} }
\label{Fig.5-4}
\end{figure}

We confirmed that the AH-tile and BH-tile shown in Figures~\ref{Fig.2-1} and \ref{Fig.2-2}, and the 
AT-tile and BT-tile shown in Figures~\ref{Fig.2-5} and \ref{Fig.2-6}, all have a Heesch number of 1 
(See Figures~\ref{Fig.5-1}--\ref{Fig.5-4}). That is, each of them can form a first corona but cannot 
form a second corona. Therefore, each of them does not admit monohedral tilings (cannot generate 
a tiling on its own). Note that there is no guarantee that concave polygons created using the same 
method as described in Section~\ref{section2} have the property of not admitting monohedral tilings 
(i.e., having a finite Heesch number)\footnote{ 
At present, we cannot claim that each of AH-tile, BH-tile, AT-tile, and BT-tile corresponding to the 
hat and turtle tiles derived from convex pentagonal monotiles belonging to the Type 5 family does 
not admit monohedral tilings. We have confirmed only that the AH-tile, BH-tile, AT-tile, and BT-tile 
constructed under specific conditions from convex pentagonal monotiles belonging to the Type 5 family 
have a Heesch number of 1. Note that, in this manuscript, AH-tile, BH-tile, AT-tile, and BT-tile generally 
refer to the tiles created based on the convex pentagonal monotiles belonging to both the Type 1 
and Type 5 families shown in Figure~\ref{Fig.1-4}.
}. In fact, Amfirifma and Smith created two types of concave polygons using the same method 
as described in Section~\ref{section2} and presented an example of a tile set consisting of these 
two types of concave polygons that generates $\mathscr{T}_h$, where the tile set contains a 
concave polygon that can generate periodic tilings (see \ref{appC}).

Furthermore, consider, for example, a tile set consisting of the two types of tiles (concave polygons) 
AH-tile and BH-tile. Although it is known that this tile set can generate $\mathscr{T}_h$, it should 
be noted that it may also generate other tilings\footnote{ 
Because the AH-tile and BH-tile have a Heesch number of 1, any other tilings generated by them 
are dihedral.
}. If tilings other than $\mathscr{T}_h$ can be generated, it is also possible that some of these tilings 
are periodic.

\vspace{15\baselineskip}

\vspace{4\baselineskip}

\subsection{Periodic tilings generated by tile sets containing AH-tile, BH-tile, AT-tile, and BT-tile}
\label{subsection5.2}

In this section, we present periodic tilings admitted by tile sets consisting of two or more tile types selected 
from AH-tile, BH-tile, AT-tile, and BT-tile. As described in Section~\ref{subsection5.1}, the AH-tile, BH-tile, 
AT-tile, and BT-tile each have a Heesch number of 1. Therefore, a tile set consisting of any one of these 
tile types alone does not admit tilings. Consequently, we consider tile sets obtained by selecting two or 
more tile types from AH-tile, BH-tile, AT-tile, and BT-tile.

The periodic tilings presented in this section are derived from the periodic tilings generated by the tile set 
consisting of the two types of tiles (concave polygons), namely the hat tile and the turtle tile, introduced in 
Section~\ref{section4}.

\begin{figure}[H]
 \centering\includegraphics[width=15cm,clip]{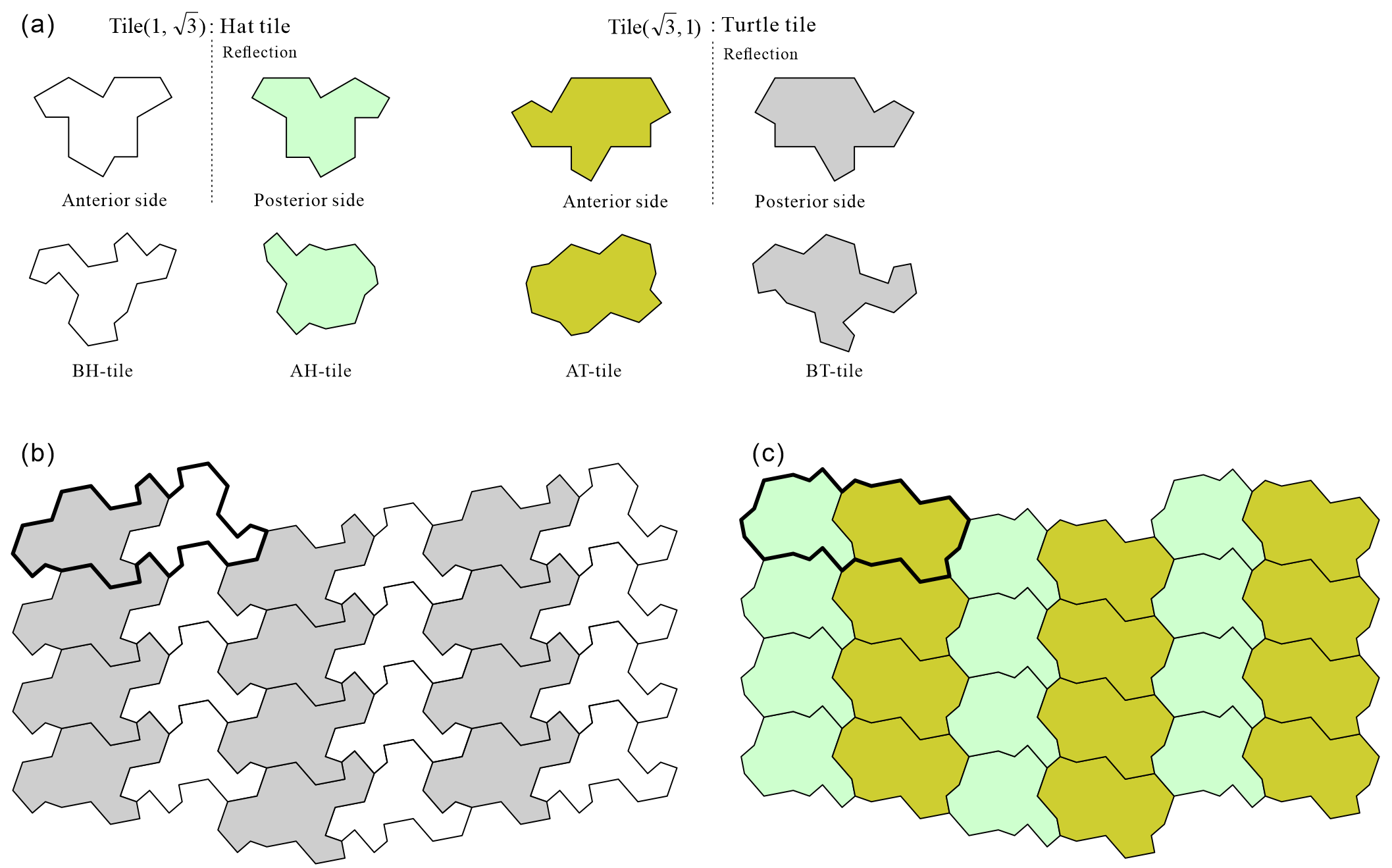} 
  \caption{{\small Periodic tilings admitted by the tile set consisting of BH-tile and BT-tile, and by the 
  tile set consisting of AH-tile and AT-tile.} }
\label{Fig.5-5}
\end{figure}

First, as shown in Figure~\ref{Fig.5-5}(a), the anterior-side hat tile is associated with BH-tile, the posterior-side 
hat tile with AH-tile, the anterior-side turtle tile with AT-tile, and the posterior-side turtle tile with BT-tile. 
As a result, the periodic tilings shown in Figures~\ref{Fig.4-1}(b) and (c) correspond to the tilings shown in 
Figures~\ref{Fig.5-5}(b) and (c), respectively. Therefore, it follows that the tile set consisting of BH-tile and 
BT-tile admits periodic tilings, and that the tile set consisting of AH-tile and AT-tile also admits periodic tilings.

\vspace{5\baselineskip}

\begin{figure}[H]
 \centering\includegraphics[width=15cm,clip]{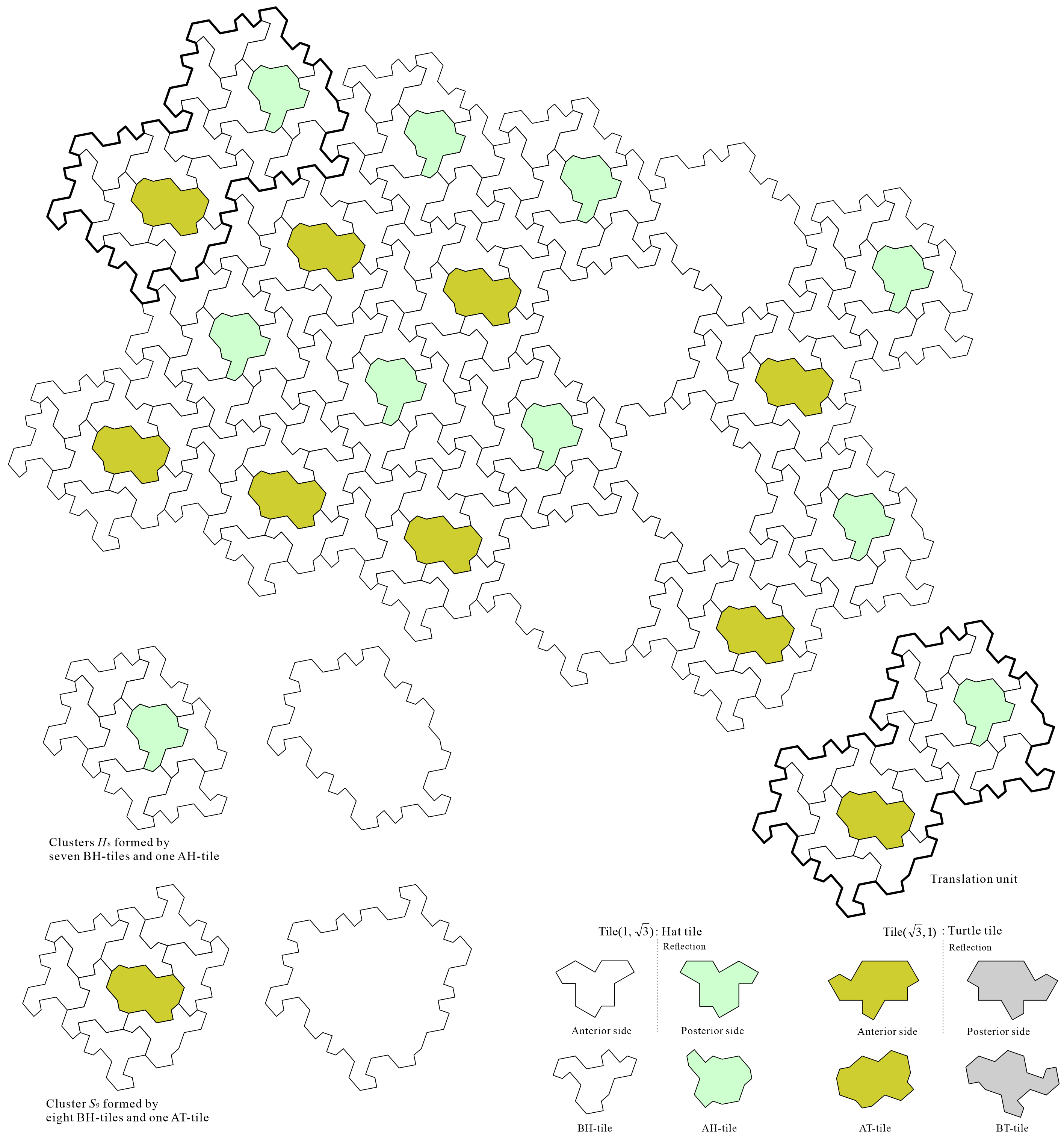} 
  \caption{{\small Periodic tiling generated using a translation unit constructed from AH-tile, BH-tile, and AT-tile.} }
\label{Fig.5-6}
\end{figure}

As mentioned in Section~\ref{section1}, a tile set (set of prototiles) is said to admit tilings if tilings can be 
generated using some of the tiles (prototiles) in the tile set. Because the periodic tiling shown in 
Figure~\ref{Fig.5-5}(c) can be generated using AH-tile and AT-tile, the tile set consisting of AH-tile, BH-tile, 
and AT-tile admits periodic tilings. Furthermore, from cluster $H_8$ in Figure~\ref{Fig.2-1}, cluster $S_9$ 
in Figure~\ref{Fig.3-4}, and the periodic tiling shown in Figure~\ref{Fig.4-2}, it follows that the tile set 
consisting of AH-tile, BH-tile, and AT-tile can also generate the 3-hedral periodic tiling shown in 
Figure~\ref{Fig.5-6}. (In this figure, the anterior-side hat tile, the posterior-side hat tile, the anterior-side 
turtle tile, and the posterior-side turtle tile correspond to BH-tile, AH-tile, AT-tile, and BT-tile, 
respectively.)

\begin{figure}[H]
 \centering\includegraphics[width=15cm,clip]{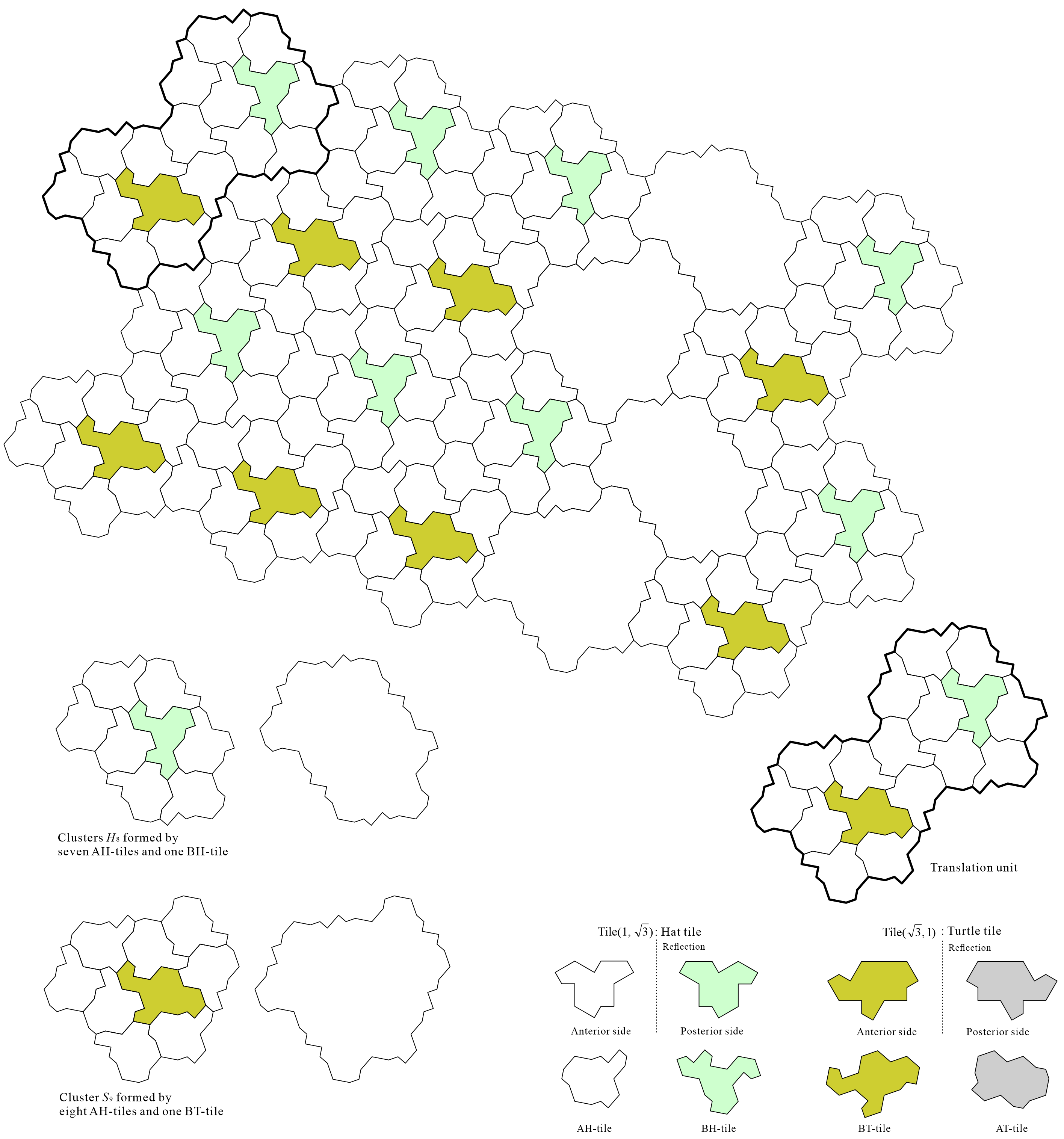} 
  \caption{{\small Periodic tiling generated using a translation unit constructed from AH-tile, BH-tile, and BT-tile.} }
\label{Fig.5-7}
\end{figure}

Because the periodic tiling shown in Figure~\ref{Fig.5-5}(b) can be generated using BH-tile and BT-tile, the tile 
set consisting of AH-tile, BH-tile, and BT-tile admits periodic tilings. Furthermore, from cluster $H_8$ in 
Figure~\ref{Fig.2-2}, cluster $S_9$ in Figure~\ref{Fig.3-8}, and the periodic tiling shown in Figure~\ref{Fig.4-2}, 
it follows that the tile set consisting of AH-tile, BH-tile, and BT-tile can also generate the 3-hedral periodic 
tiling shown in Figure~\ref{Fig.5-7}. (In this figure, the anterior-side hat tile, the posterior-side hat tile, 
the anterior-side turtle tile, and the posterior-side turtle tile correspond to AH-tile, BH-tile, BT-tile, 
and AT-tile, respectively.)

\vspace{2\baselineskip}

\begin{figure}[H]
 \centering\includegraphics[width=15cm,clip]{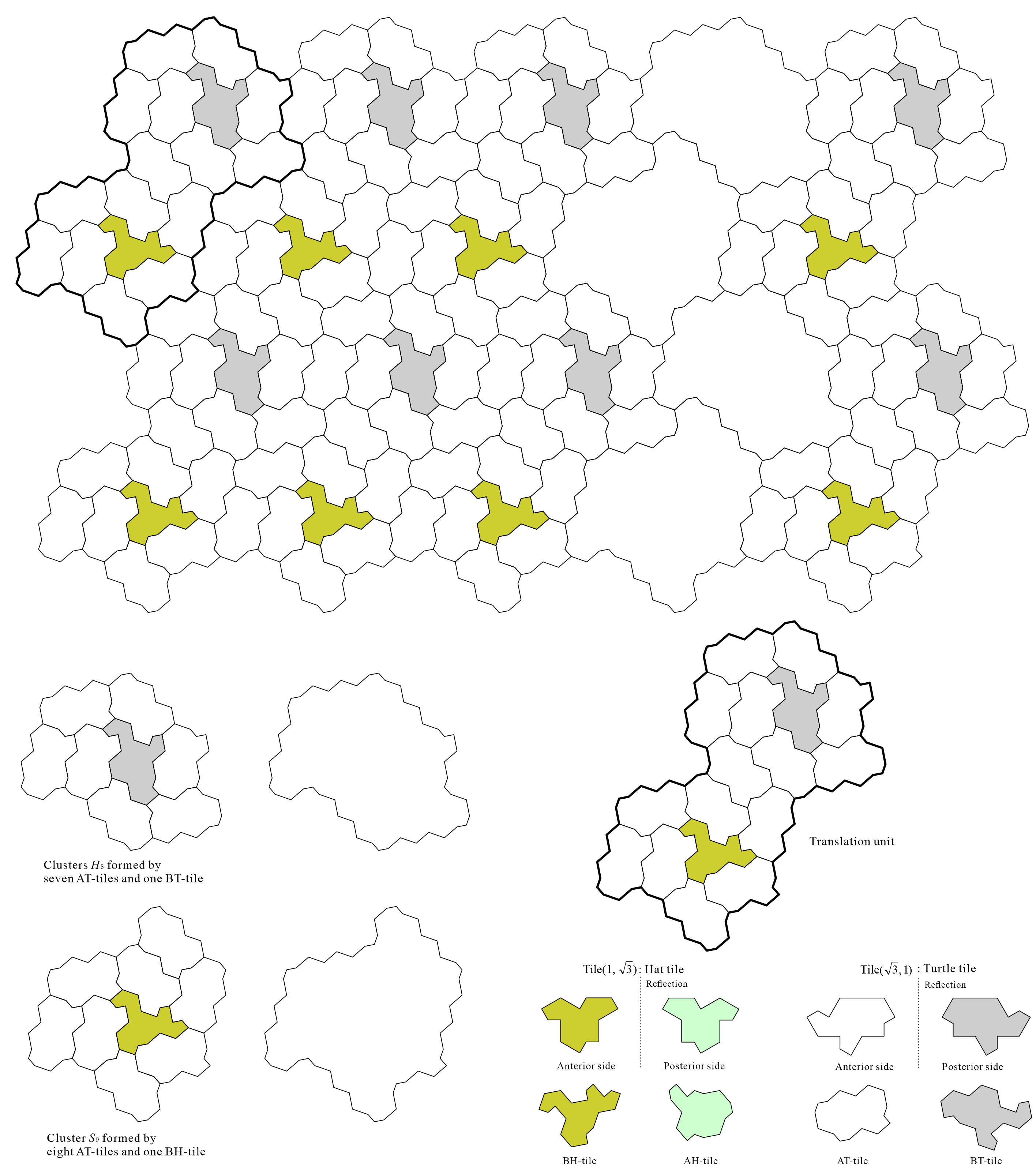} 
  \caption{{\small Periodic tiling generated using a translation unit constructed from BH-tile, AT-tile, and BT-tile.} }
\label{Fig.5-8}
\end{figure}

Because the periodic tiling shown in Figure~\ref{Fig.5-5}(b) can be generated using BH-tile and BT-tile, 
the tile set consisting of BH-tile, AT-tile, and BT-tile admits periodic tilings. Furthermore, from cluster $H_8$ in 
Figure~\ref{Fig.2-5}, cluster $S_9$ in Figure~\ref{Fig.3-6}, and the periodic tiling shown in Figure~\ref{Fig.4-3}, 
it follows that the tile set consisting of BH-tile, AT-tile, and BT-tile can also generate the 3-hedral periodic 
tiling shown in Figure~\ref{Fig.5-8}. (In this figure, the anterior-side hat tile, the posterior-side hat tile, 
the anterior-side turtle tile, and the posterior-side turtle tile correspond to BH-tile, AH-tile, AT-tile, 
and BT-tile, respectively.)

\begin{figure}[H]
 \centering\includegraphics[width=15cm,clip]{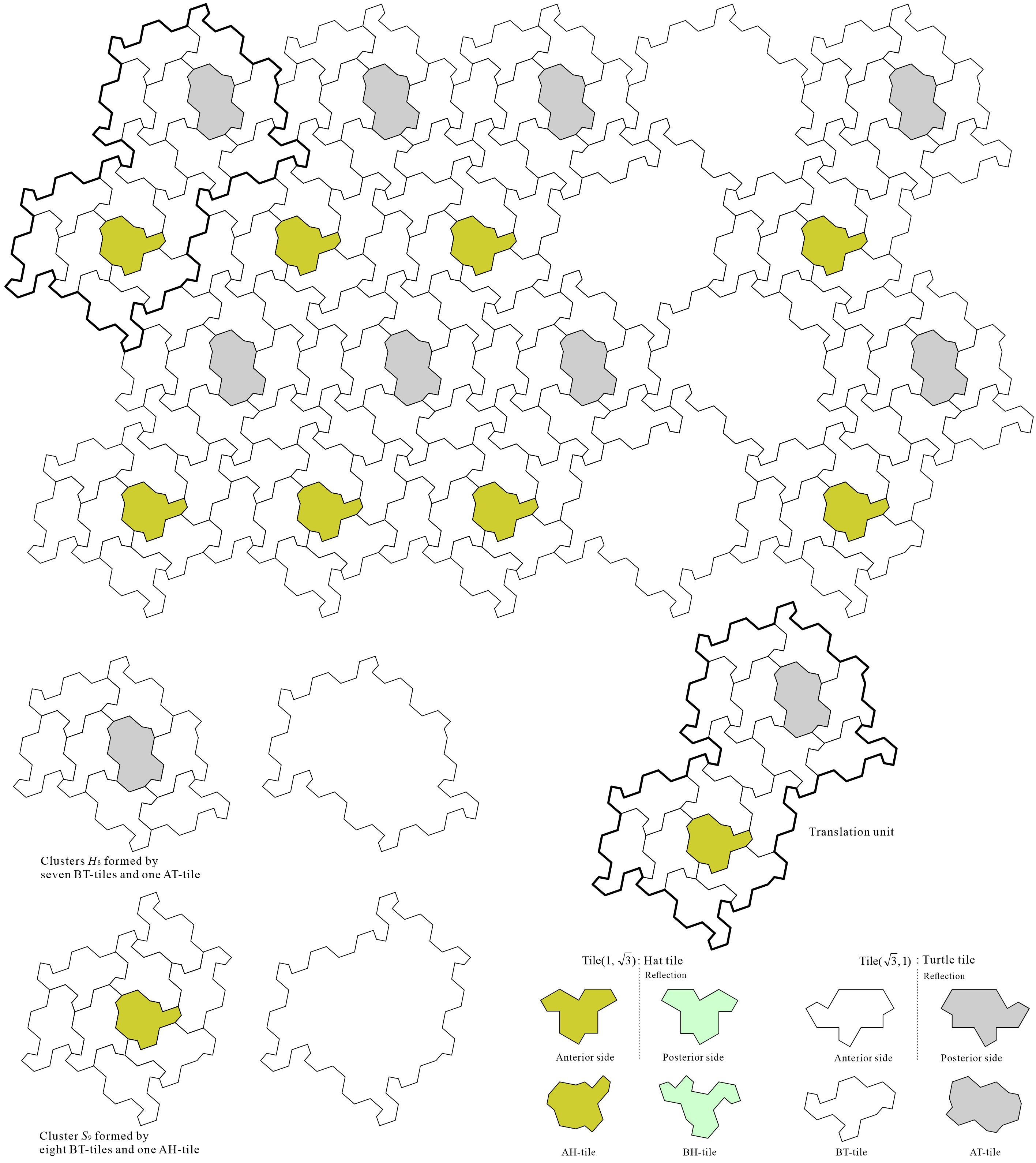} 
  \caption{{\small Periodic tiling generated using a translation unit constructed from AH-tile, AT-tile, and BT-tile.} }
\label{Fig.5-9}
\end{figure}

Because the periodic tiling shown in Figure~\ref{Fig.5-5}(c) can be generated using AH-tile and AT-tile, the tile 
set consisting of AH-tile, AT-tile, and BT-tile admits periodic tilings. Furthermore, from cluster $H_8$ in 
Figure~\ref{Fig.2-6}, cluster $S_9$ in Figure~\ref{Fig.3-10}, and the periodic tiling shown in Figure~\ref{Fig.4-3}, 
it follows that the tile set consisting of AH-tile, AT-tile, and BT-tile can also generate the 3-hedral periodic 
tiling shown in Figure~\ref{Fig.5-9}. (In this figure, the anterior-side hat tile, the posterior-side hat tile, 
the anterior-side turtle tile, and the posterior-side turtle tile correspond to AH-tile, BH-tile, BT-tile, 
and AT-tile, respectively.)

From the above results, it follows that all tile sets consisting of tiles (concave polygons) obtained by 
selecting three or four tile types from AH-tile, BH-tile, AT-tile, and BT-tile admit periodic tilings.

\begin{figure}[H]
 \centering\includegraphics[width=15cm,clip]{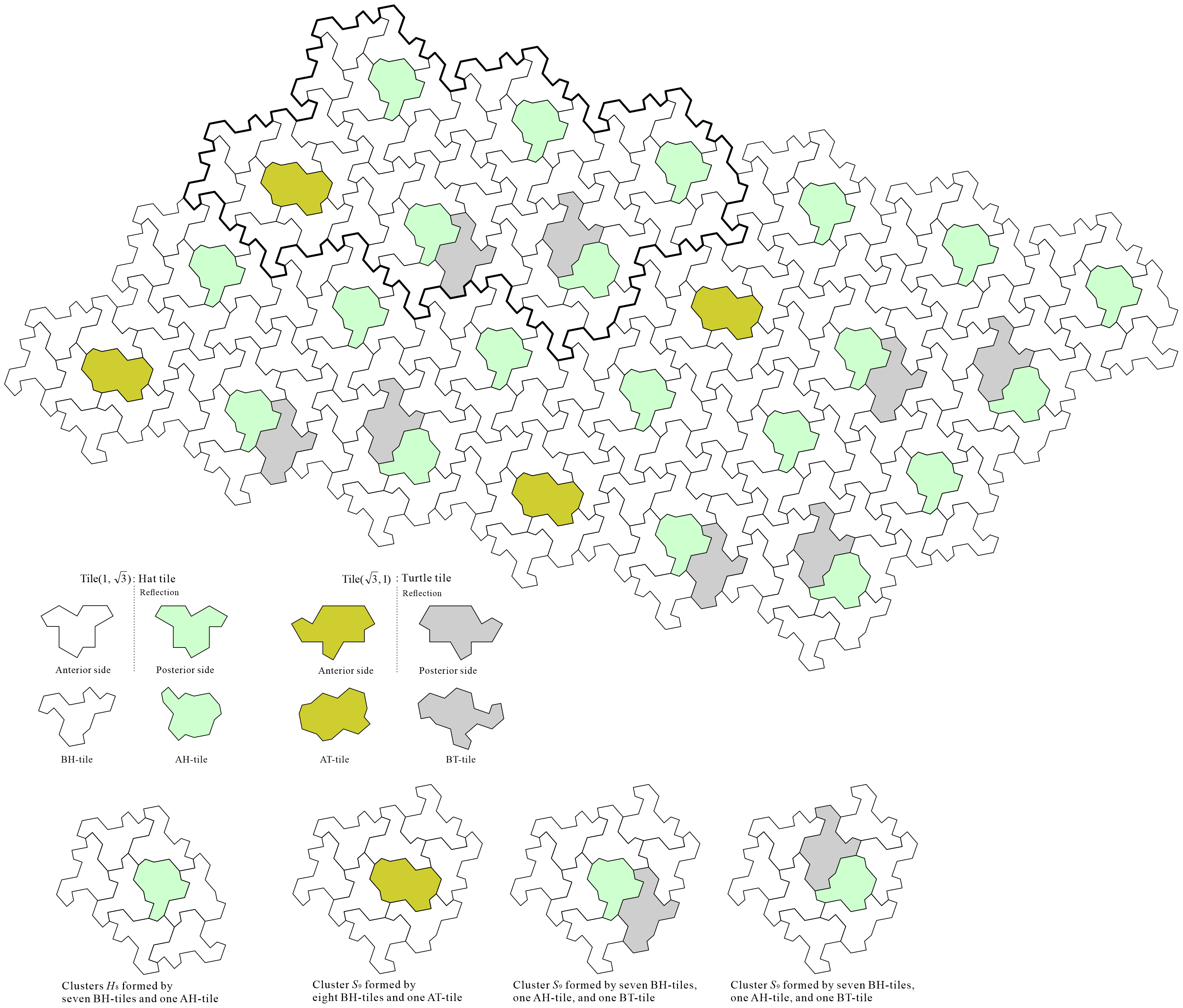} 
  \caption{{\small Periodic tiling generated using a translation unit constructed from AH-tile, BH-tile, AT-tile, and BT-tile.} }
\label{Fig.5-10}
\end{figure}

Furthermore, periodic tilings can also be generated using all four tile types, namely AH-tile, BH-tile, AT-tile, 
and BT-tile. For example, by associating the anterior-side hat tile with BH-tile, the posterior-side hat tile with 
AH-tile, the anterior-side turtle tile with AT-tile, and the posterior-side turtle tile with BT-tile in the periodic 
tiling shown in Figure~\ref{Fig.4-7}, the 4-hedral periodic tiling shown in Figure~\ref{Fig.5-10} is obtained. 
As in the case of the hat tile and the turtle tile, the tile sets consisting of tiles obtained by selecting three 
or four tile types from AH-tile, BH-tile, AT-tile, and BT-tile, can generate a variety of tiling patterns 
(see \ref{appA}).

\vspace{3\baselineskip}

\subsection{Properties of tile sets containing AH-tile, BH-tile, AT-tile, and BT-tile}
\label{subsection5.3}

The following summarizes the current knowledge and the remaining open questions regarding the 
AH-tile, BH-tile, AT-tile, and BT-tile derived from convex pentagonal monotiles belonging to both 
the Type 1 and Type 5 families shown in Figure~\ref{Fig.1-4} (see also Table~\ref{Table.5-1}).

\begin{itemize}
\item The AH-tile is a concave polygon with Heesch number 1 and does not admit monohedral tilings.
\item The BH-tile is a concave polygon with Heesch number 1 and does not admit monohedral tilings.
\item The AT-tile is a concave polygon with Heesch number 1 and does not admit monohedral tilings.
\item The BT-tile is a concave polygon with Heesch number 1 and does not admit monohedral tilings.
\item A tile set consisting of two types of concave polygons, AH-tile and BH-tile, is known to 
generate $\mathscr{T}_h$, but whether it admits other tilings remains unknown.
\item A tile set consisting of two types of concave polygons, AT-tile and BT-tile, is known to 
generate $\mathscr{T}_h$, but whether it admits other tilings remains unknown.
\item A tile set consisting of two types of concave polygons, AT-tile and BH-tile, is known to 
generate $\mathscr{T}_s$, but whether it admits other tilings remains unknown.
\item A tile set consisting of two types of concave polygons, AH-tile and BT-tile, is known to 
generate $\mathscr{T}_s$, but whether it admits other tilings remains unknown.
\item A tile set consisting of two types of concave polygons, AH-tile and AT-tile, admits periodic tilings.
\item A tile set consisting of two types of concave polygons, BH-tile and BT-tile, admits periodic tilings.
\item All tile sets consisting of concave polygons obtained by selecting three or four tile types from 
AH-tile, BH-tile, AT-tile, and BT-tile admit periodic tilings.
\end{itemize}

\vspace{\baselineskip}

Based on the summary above, the following question arises: ``Is there a tile set consisting of two 
types of concave polygons, obtained by selecting two tile types from AH-tile, BH-tile, AT-tile, 
and BT-tile, all of which have Heesch number 1, that is an aperiodic tile set (i.e., an aperiodic set 
of tiles)?'' For example, suppose that a tile set consisting of the two types of concave polygons 
AT-tile and BT-tile generates only $\mathscr{T}_h$, or admits other tilings in addition to $\mathscr{T}_h$, 
and that all such tilings are non-periodic. In that case, the tile set is an aperiodic tile set.

\setlength{\tabcolsep}{8pt}       
\renewcommand{\arraystretch}{1.4} 

\begin{table}[H]
\caption{Relationship among the AH-tile, BH-tile, AT-tile, and BT-tile}
\label{Table.5-1}
\centering
{\small
\begin{tabular}{|l|p{2.9cm}|p{2.9cm}|p{2.9cm}|p{2.9cm}|}
\hline
 & \multicolumn{1}{c|}{AH-tile} & \multicolumn{1}{c|}{BH-tile} & \multicolumn{1}{c|}{AT-tile} & \multicolumn{1}{c|}{BT-tile} \\ \hline

\multicolumn{1}{|c|}{\multirow{8}{*}{AH-tile}} & \raggedright AH-tile has Heesch number 1 (i.e., AH-tile does not admit monohedral tilings). & \raggedright A tile set consisting of two types of concave polygons, AH-tile and BH-tile, is known to generate $\mathscr{T}_h$, but whether it admits other tilings remains unknown. & \raggedright A tile set consisting of AH-tile and AT-tile admits periodic tilings. & \raggedright A tile set consisting of two types of concave polygons, AH-tile and BT-tile, is known to generate $\mathscr{T}_s$, but whether it admits other tilings remains unknown. \tabularnewline \hline

\multicolumn{1}{|c|}{\multirow{8}{*}{BH-tile}} & & \raggedright BH-tile has Heesch number 1 (i.e., BH-tile does not admit monohedral tilings). & \raggedright A tile set consisting of two types of concave polygons, AT-tile and BH-tile, is known to generate $\mathscr{T}_s$, but whether it admits other tilings remains unknown. & \raggedright A tile set consisting of BH-tile and BT-tile admits periodic tilings. \tabularnewline \hline

\multicolumn{1}{|c|}{\multirow{8}{*}{AT-tile}} & & & \raggedright AT-tile has Heesch number 1 (i.e., AT-tile does not admit monohedral tilings). & \raggedright A tile set consisting of two types of concave polygons, AT-tile and BT-tile, is known to generate $\mathscr{T}_h$, but whether it admits other tilings remains unknown. \tabularnewline \hline

\multicolumn{1}{|c|}{\multirow{4}{*}{BT-tile}} & & & & \raggedright BT-tile has Heesch number 1 (i.e., BT-tile does not admit monohedral tilings). \tabularnewline \hline
\end{tabular}
}
\end{table}

\renewcommand{\arraystretch}{1.0}

\vspace{2\baselineskip}
\section{Tile sets generating non-periodic tilings in which the densities of two types of tiles are mutually reversed}
\label{section6}

This section focuses on the properties of the tile sets consisting of the two types of tiles discussed 
in Sections~\ref{section2} and \ref{section3} that generate $\mathscr{T}_h$ and $\mathscr{T}_s$, 
namely, their ability to generate non-periodic tilings in which the densities of the two types of tiles are 
mutually reversed.

\vspace{\baselineskip}

Let \mbox{$P\{\text{A-tile}, \text{B-tile}\}$} denote a tile set consisting of the two types of tiles, A-tile and 
B-tile, without matching rules, where A-tile and B-tile are non-congruent (i.e., not mirror images of each other), 
and assume that \mbox{$P\{\text{A-tile}, \text{B-tile}\}$} can generate tiling using both A-tile and B-tile. 
The aperiodic tile set (i.e., aperiodic set of tiles) corresponding to \mbox{$P\{\text{A-tile}, \text{B-tile}\}$} is 
denoted by \mbox{$\textit{ASP}\{\text{A-tile}, \text{B-tile}\}$}. An \mbox{$\textit{ASP}\{\text{A-tile}, \text{B-tile}\}$} 
that generates non-periodic tilings in which the densities of A-tile and B-tile are mutually reversed is denoted by 
\mbox{$\textit{ASPmr}\{\text{A-tile}, \text{B-tile}\}$}\footnote{ 
It is assumed that the tilings generated by \mbox{$P\{\text{A-tile}, \text{B-tile}\}$} can be obtained 
without restrictions on tiling generation, although they may also be obtained when restrictions are imposed. 
Furthermore, the property of generating non-periodic tilings in which the relative densities of A-tile and 
B-tile are reversed suggests that \mbox{$P\{\text{A-tile}, \text{B-tile}\}$} generates at least two 
patterns of dihedral tiling. Note that, as described in Section~\ref{section1}, in this study, the 
anterior-side and posterior-side of a tile are regarded as the same type of tile; therefore, 
A-tile and B-tile are naturally not congruent. However, for clarity, we have explicitly added the 
explanation that ``A-tile and B-tile are not congruent.'' Consequently, while the hat tile or turtle 
tile can be regarded as generating non-periodic tilings in which the relative densities of 
anterior-side and posterior-side tiles are mutually reversed, in those cases the 
two tiles (concave polygons) are congruent, and thus they do not correspond to 
\mbox{$\textit{ASP}\{\text{A-tile}, \text{B-tile}\}$} or \mbox{$\textit{ASPmr}\{\text{A-tile}, \text{B-tile}\}$}.
}.

The combinations of the Heesch numbers and aperiodicity exhibited by the tiles in 
\mbox{$\textit{ASP}\{\text{A-tile}, \text{B-tile}\}$} and \mbox{$\textit{ASPmr}\{\text{A-tile}, \text{B-tile}\}$} can 
be classified into the following six categories:

\begin{enumerate}[label={(Pair-\arabic*)}, leftmargin=*]
\item Both tiles have Heesch number 0.
\item One tile has Heesch number 0, and the other has a positive Heesch number.
\item Both tiles have positive Heesch numbers.
\item One tile is an aperiodic monotile, and the other has Heesch number 0.
\item One tile is an aperiodic monotile, and the other has a positive Heesch number.
\item Both tiles are aperiodic monotiles.
\end{enumerate}

\vspace{\baselineskip}

It should be noted that tile sets corresponding to \mbox{$\textit{ASP}\{\text{A-tile}, \text{B-tile}\}$}, 
that is, aperiodic tile sets consisting of two types of concave polygons without matching rules, already exist. 
For example, the pairs of concave polygons obtained from the outer boundaries of clusters $H_7$ and $H_8$ 
formed at each step of the $H_7/H_8$ substitution may be regarded as aperiodic tile sets consisting of 
two types of concave polygons without matching rules. This observation suggests that classifying 
tile sets solely according to the shapes of their constituent tiles is not particularly meaningful. 
To the best of our knowledge, no summary or clear information is available that identifies which 
known tile sets corresponding to \mbox{$\textit{ASP}\{\text{A-tile}, \text{B-tile}\}$} fall into categories 
(Pair-1)--(Pair-6), and compiling such a summary would be of interest \cite{ref14}. However, we focus on the 
possibility that tile sets based on the AH-tile, BH-tile, AT-tile, and BT-tile introduced in Sections~\ref{section2},  
\ref{section3}, and \ref{section5} are \mbox{$\textit{ASPmr}\{\text{A-tile}, \text{B-tile}\}$}. 
If no restrictions are imposed on the tiles or tile sets, a tile set capable of generating non-periodic tilings 
in which the densities of the two types of tiles are mutually reversed is not particularly remarkable. 
However, the situation changes significantly if the tile set is aperiodic. Therefore, we believe that whether 
an \mbox{$\textit{ASPmr}\{\text{A-tile}, \text{B-tile}\}$} exists is a very interesting question.

As mentioned in Footnote~\ref{footnote2}, a set of prototiles is said to be aperiodic if it admits tilings of 
the plane and all tilings generated from those prototiles are non-periodic \cite{ref1, ref2, ref3}. 
For example, in an aperiodic tile set consisting of two types of tiles, if either tile is considered 
individually, it cannot generate periodic tilings (i.e., neither tile can generate periodic tilings on its own). 
Conversely, if all tilings generated using the tiles in the tile set are non-periodic, then the tile set is 
aperiodic, and the tiles themselves may be aperiodic monotiles. For example, suppose there exists a 
tile set \mbox{$P\{\text{A-tile}, \text{B-tile}\}$} in which A-tile is an aperiodic monotile and B-tile has a 
positive Heesch number. Clearly, A-tile admits only non-periodic tilings on its own, whereas B-tile does 
not admit monohedral tilings. If all tilings generated using both A-tile and B-tile are non-periodic, then 
this tile set is an \mbox{$\textit{ASP}\{\text{A-tile}, \text{B-tile}\}$} corresponding to (Pair-5).

\mbox{$P\{\text{hat tile}, \text{turtle tile}\}$} is a tile set consisting of two types of tiles (concave polygons) 
without matching rules that has the property of generating non-periodic tilings in which the densities of 
the two types of tiles are mutually reversed. However, because this tile set also admits the periodic 
tilings introduced in Section~\ref{section4}, it is not an \mbox{$\textit{ASPmr}\{\text{A-tile}, \text{B-tile}\}$}.

\vspace{\baselineskip}

\begin{ques}
Does an \mbox{$\textit{ASPmr}\{\text{A-tile}, \text{B-tile}\}$} exist?
\end{ques}

Regarding this question, the following possibilities can currently be identified:
\begin{itemize}
\item \mbox{$P\{\text{AH-tile}, \text{BH-tile}\}$} can generate $\mathscr{T}_h$, in which both 
AH-tile and BH-tile have Heesch number 1 and the densities of the two types of polygons are 
mutually reversed. Accordingly, \mbox{$P\{\text{AH-tile}, \text{BH-tile}\}$} may be an 
\mbox{$\textit{ASPmr}\{\text{AH-tile}, \text{BH-tile}\}$} corresponding to (Pair-3).
\item \mbox{$P\{\text{AT-tile}, \text{BT-tile}\}$} can generate $\mathscr{T}_h$, in which both 
AT-tile and BT-tile have Heesch number 1 and the densities of the two types of polygons are 
mutually reversed. Accordingly, \mbox{$P\{\text{AT-tile}, \text{BT-tile}\}$} may be an 
\mbox{$\textit{ASPmr}\{\text{AT-tile}, \text{BT-tile}\}$} corresponding to (Pair-3).
\item \mbox{$P\{\text{AT-tile}, \text{BH-tile}\}$} can generate $\mathscr{T}_s$, in which both 
AT-tile and BH-tile have Heesch number 1 and the densities of the two types of polygons are 
mutually reversed. Accordingly, \mbox{$P\{\text{AT-tile}, \text{BH-tile}\}$} may be an 
\mbox{$\textit{ASPmr}\{\text{AT-tile}, \text{BH-tile}\}$} corresponding to (Pair-3).
\item \mbox{$P\{\text{AH-tile}, \text{BT-tile}\}$} can generate $\mathscr{T}_s$, in which both 
AH-tile and BT-tile have Heesch number 1 and the densities of the two types of polygons are 
mutually reversed. Accordingly, \mbox{$P\{\text{AH-tile}, \text{BT-tile}\}$} may be an 
\mbox{$\textit{ASPmr}\{\text{AH-tile}, \text{BT-tile}\}$} corresponding to (Pair-3).
\end{itemize}

\vspace{1\baselineskip}
\section{Conclusion}
\label{section7}

At present, we are primarily investigating whether the tile set \mbox{$P\{\text{AT-tile}, \text{BT-tile}\}$} 
is an \mbox{$\textit{ASPmr}\{\text{AT-tile}, \text{BT-tile}\}$} corresponding to (Pair-3). We also intend 
to investigate the other tile sets. However, because the number of combinations to be examined is 
enormous, it remains unclear how far our current approach can be extended or whether it will 
ultimately yield the desired results.

\section*{Acknowledgments}
The author would like to express sincere gratitude to David Smith, Amfirifma, and Ortwin Schenker 
for their valuable information and helpful comments during the course of this research. In particular, 
the information they provided on the periodic tilings generated by the tile set consisting of the two 
types of tiles (concave polygons), hat tiles and turtle tiles, described in Section~\ref{section4} and 
\ref{appA}, as well as the example of a concave polygon admitting periodic tilings presented 
in \ref{appC}, was extremely valuable to this study. 
The author would like to thank Professor Emeritus Yoshio Agaoka of Hiroshima University for 
providing appropriate advice during the preparation of this manuscript.

\vspace{1\baselineskip}

\vspace{1\baselineskip}
\appendix
\def\thesection{Appendix \Alph{section}}
\renewcommand{\thefigure}{\Alph{section}.\arabic{figure}}
\renewcommand{\thetable}{\Alph{section}.\arabic{table}}

\section{Examples of Various Tilings}
\label{appA}

\begin{figure}[H]
 \centering\includegraphics[width=12cm,clip]{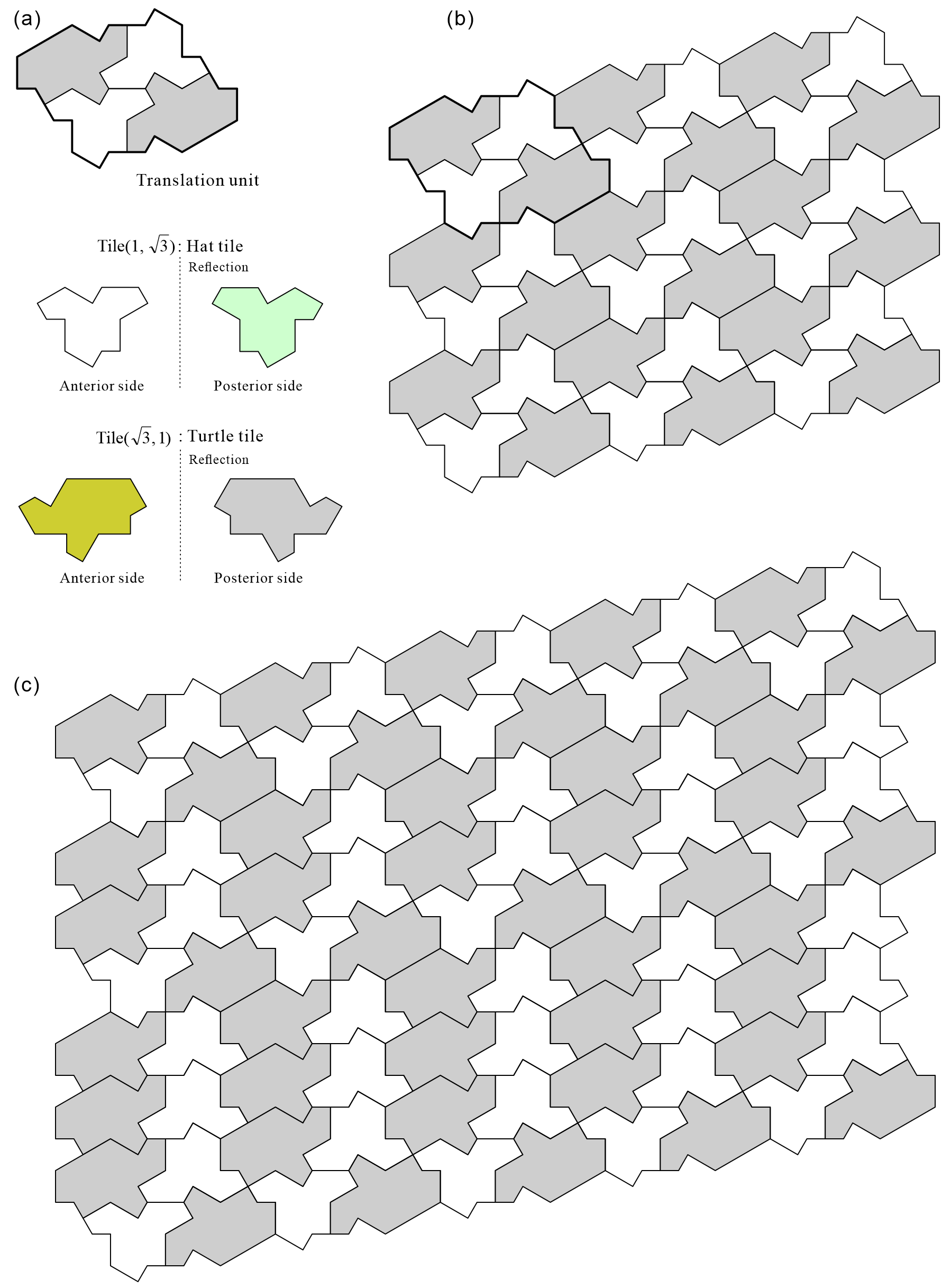} 
  \caption{{\small Tilings different from that shown in Figure~\ref{Fig.4-1}, with a $1:1$ ratio of hat tiles to turtle tiles.} }
\label{Fig.A-1}
\end{figure}

A periodic tiling such as that shown in Figure~\ref{Fig.A-1}(b) can be generated using a translation 
unit constructed from two anterior-side hat tiles and two posterior-side turtle tiles, as shown in 
Figure~\ref{Fig.A-1}(a). Furthermore, by focusing on the tilings in Figures~\ref{Fig.4-1}(b) and 
\ref{Fig.A-1}(b) and their translation units, we observe that tilings can also be generated by freely 
connecting in the vertical direction belts formed by extending these translation units through 
horizontal translation in the same direction, as shown in Figure~\ref{Fig.A-1}(c). 
The tilings in Figures~\ref{Fig.4-1} and \ref{Fig.A-1} have a $1:1$ ratio of hat tiles to turtle tiles.

\begin{figure}[H]
 \centering\includegraphics[width=15cm,clip]{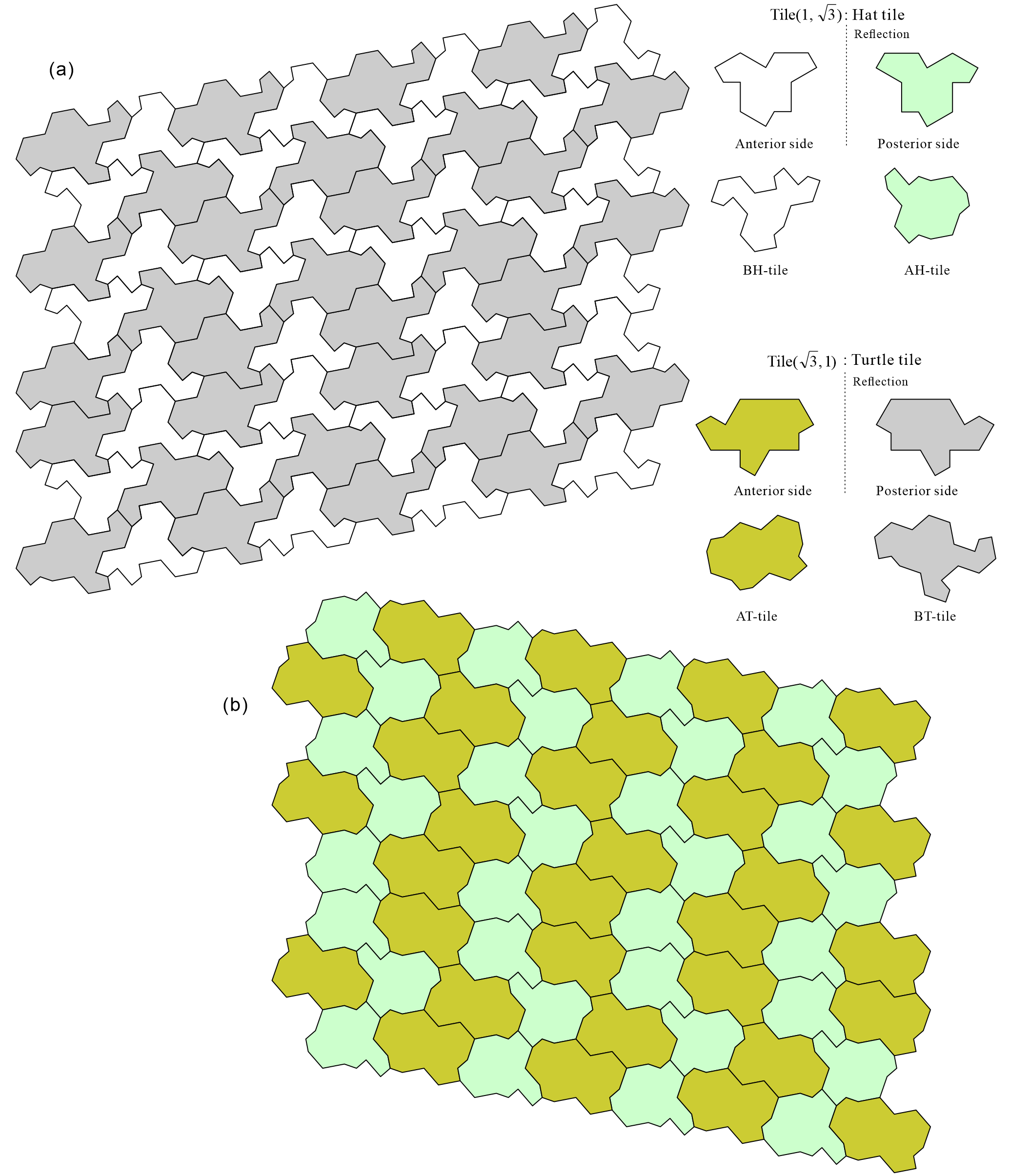} 
  \caption{{\small Tiling generated using BH-tiles and BT-tiles and tiling generated using AH-tiles and AT-tiles.} }
\label{Fig.A-2}
\end{figure}

The tiling generated using BH-tiles and BT-tiles, based on the properties and relationships of 
Figures~\ref{Fig.4-1}, \ref{Fig.5-5}, and \ref{Fig.A-1}, is shown in Figure~\ref{Fig.A-2}(a), while the tiling 
generated using AH-tiles and AT-tiles is shown in Figure~\ref{Fig.A-2}(b).

\vspace{3\baselineskip}

\begin{figure}[H]
 \centering\includegraphics[width=15cm,clip]{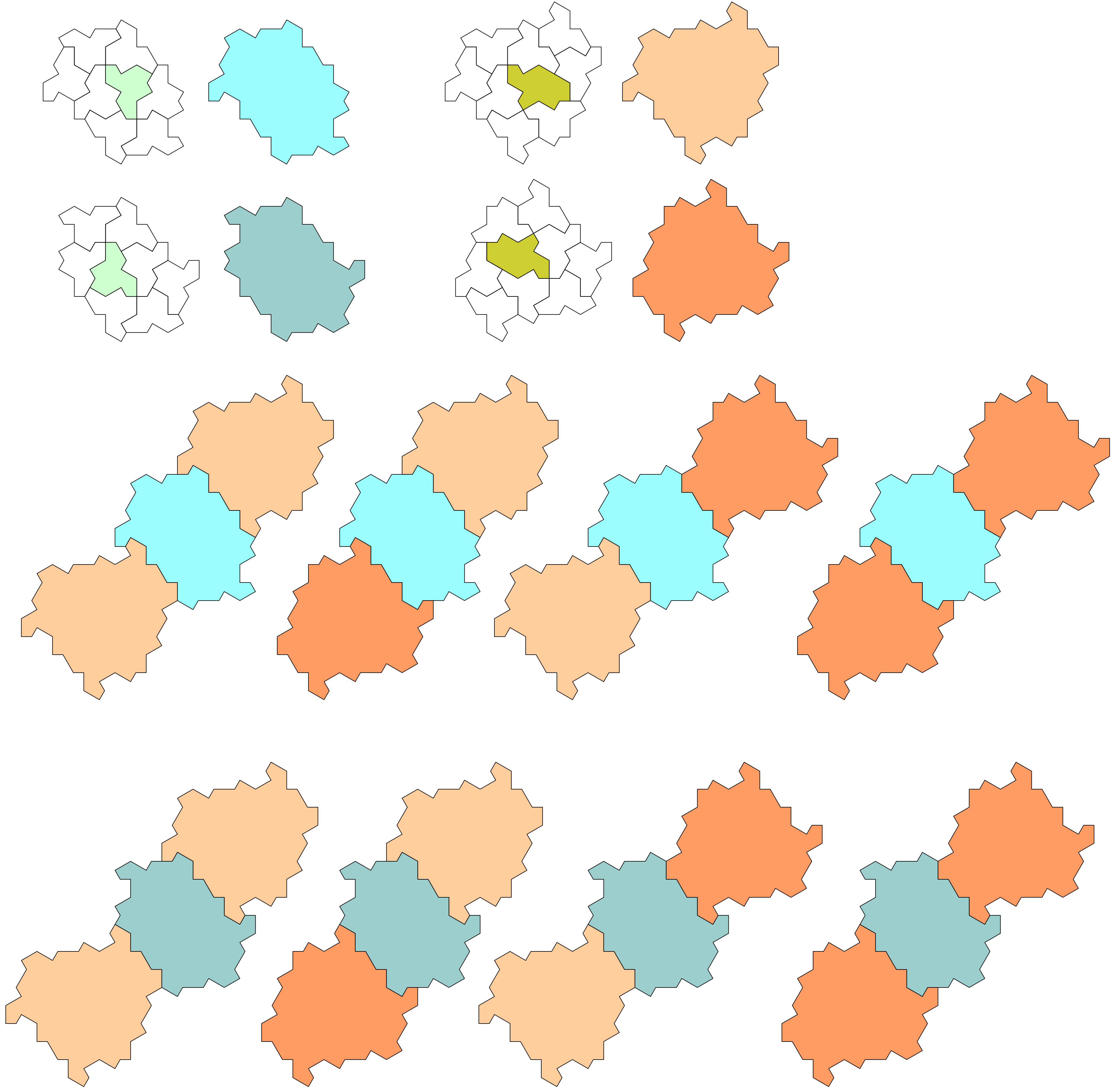} 
  \caption{{\small Variations in the ways of combining $H_8$ and $S_9$.} }
\label{Fig.A-3}
\end{figure}

Figure~\ref{Fig.A-3} shows variations in the ways of combining cluster $H_8$, formed by 
eight hat tiles shown in Figure~\ref{Fig.2-1}, and cluster $S_9$, formed by eight hat 
tiles and one turtle tile shown in Figure~\ref{Fig.3-2}. Specifically, $H_8$ and $S_9$ can 
each be freely combined either in their original orientation or after a $180^\circ$ rotation.

\vspace{10\baselineskip}

\begin{figure}[H]
 \centering\includegraphics[width=15cm,clip]{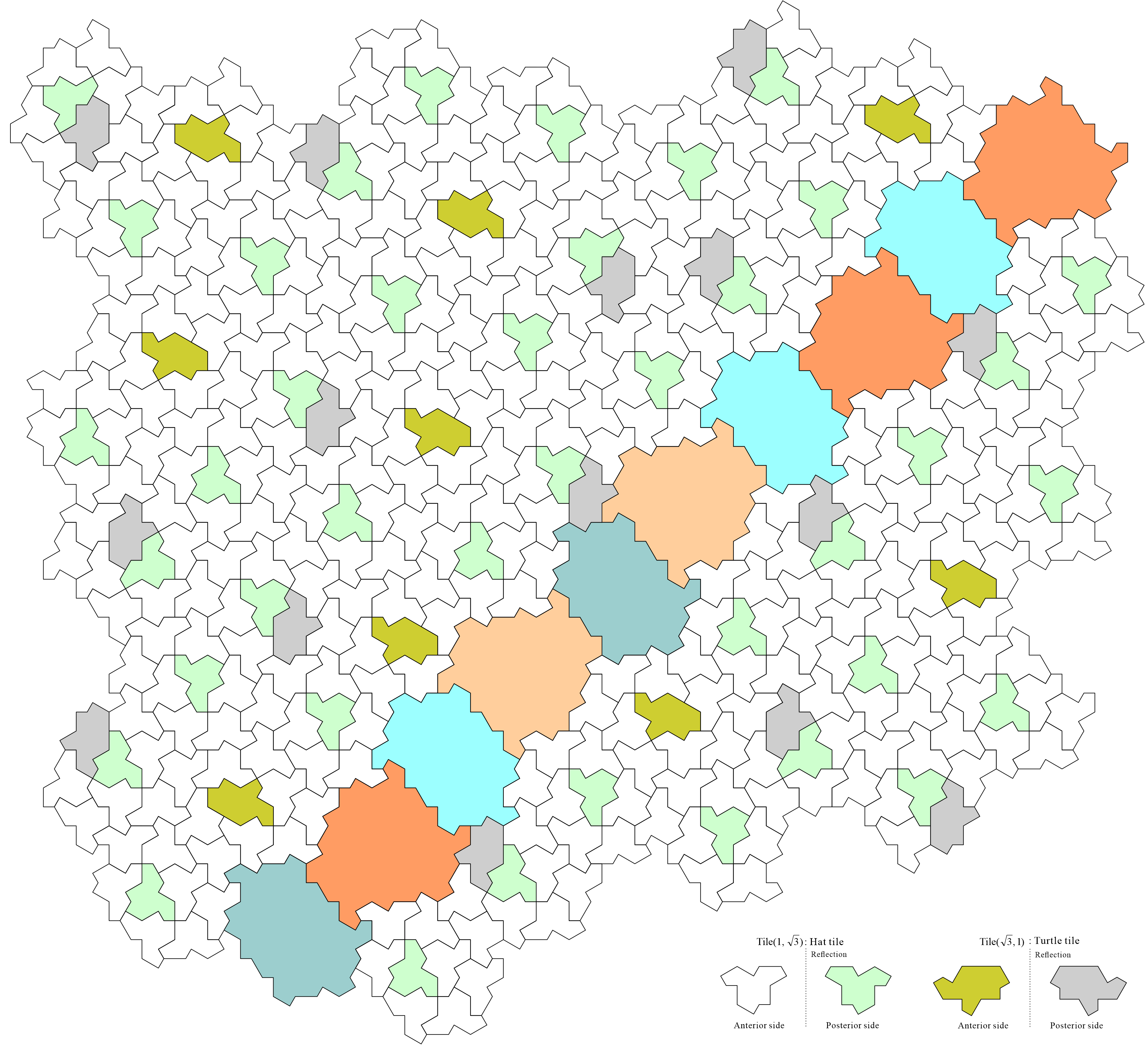} 
  \caption{{\small A tiling using hat tiles and turtle tiles generated by combining the property 
  that $H_8$ and $S_9$ can each be rotated by $180^\circ$ with the property that the 
  Mystic regions within $S_9$ can be replaced, as shown in Figure~\ref{Fig.4-5}.} }
\label{Fig.A-4}
\end{figure}

Figure~\ref{Fig.A-4} shows an example of a tiling using hat tiles and turtle tiles, generated 
by combining the property that clusters $H_8$ and $S_9$ can each be rotated by $180^\circ$ 
with the property that the Mystic regions within $S_9$ can be replaced, as shown in 
Figure~\ref{Fig.4-5}.

\vspace{20\baselineskip}

\begin{figure}[H]
 \centering\includegraphics[width=15cm,clip]{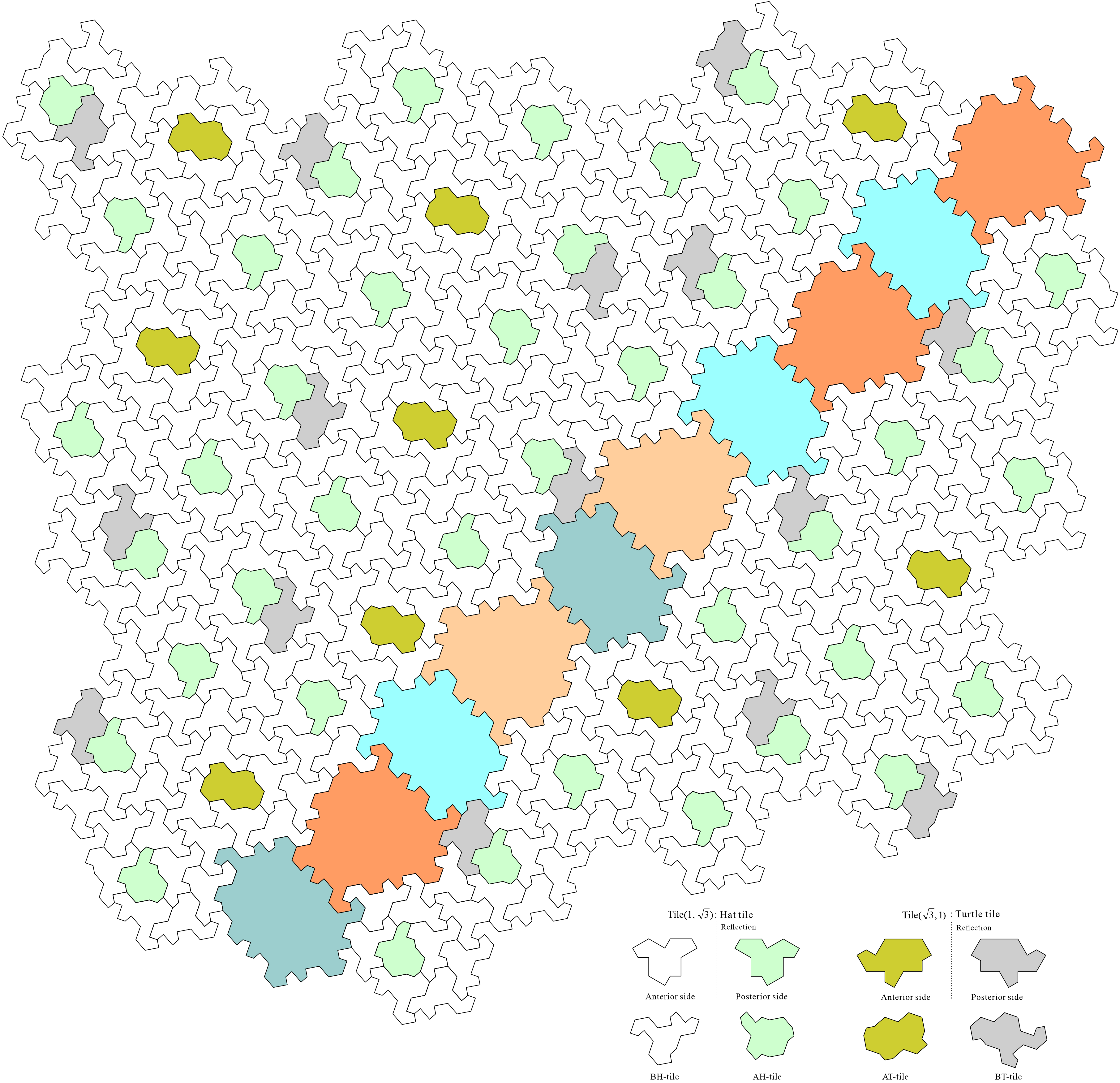} 
  \caption{{\small A tiling using AH-tile, BH-tile, AT-tile, and BT-tile generated by combining 
  the property that $H_8$ and $S_9$ can each be rotated by $180^\circ$ with the property 
  that the Mystic regions within $S_9$ can be replaced.} }
\label{Fig.A-5}
\end{figure}

Figure~\ref{Fig.A-5} shows an example of a tiling using all four tile types, AH-tile, BH-tile, 
AT-tile, and BT-tile, generated by utilizing both the property that clusters $H_8$ and $S_9$ 
can each be rotated by $180^\circ$ and the property that the Mystic regions within $S_9$ 
can be replaced.

\vspace{8\baselineskip}
\section{Correspondence with Type 5 tiling patterns}
\label{appB}

Section~\ref{section4} introduced various periodic tilings using hat tiles and turtle tiles. 
The Type 5 tiling patterns can likewise be associated with these tilings. It should be noted 
that, in order to establish the correspondence with a Type 5 tiling (i.e., to make a Type 5 tiling 
pattern appear on the surface of the tiling), the translation unit must be constructed 
using hat tiles and turtle tiles that have the same pattern. For example, if a translation unit 
is constructed using a hat tile with Pattern 1 and a turtle tile with Pattern 2, the resulting 
dihedral tiling using hat tiles and turtle tiles is valid, but the corresponding monohedral tiling 
pattern of convex pentagons (i.e., the Type 5 tiling pattern) does not appear on its surface.

\begin{figure}[H]
 \centering\includegraphics[width=15cm,clip]{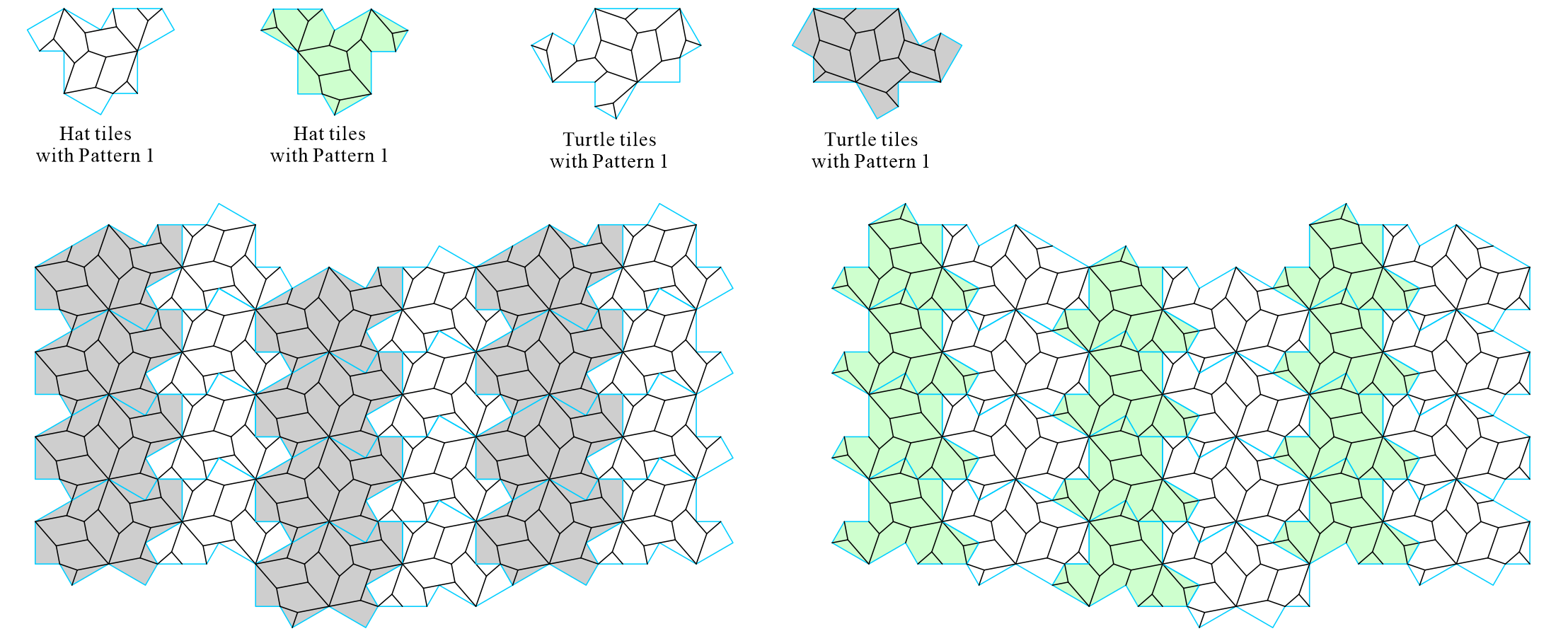} 
  \caption{{\small Periodic tiling in Figure~\ref{Fig.4-1} reconstructed using hat tiles and turtle tiles with Pattern 1.} }
\label{Fig.B-1}
\end{figure}

\begin{figure}[H]
 \centering\includegraphics[width=15cm,clip]{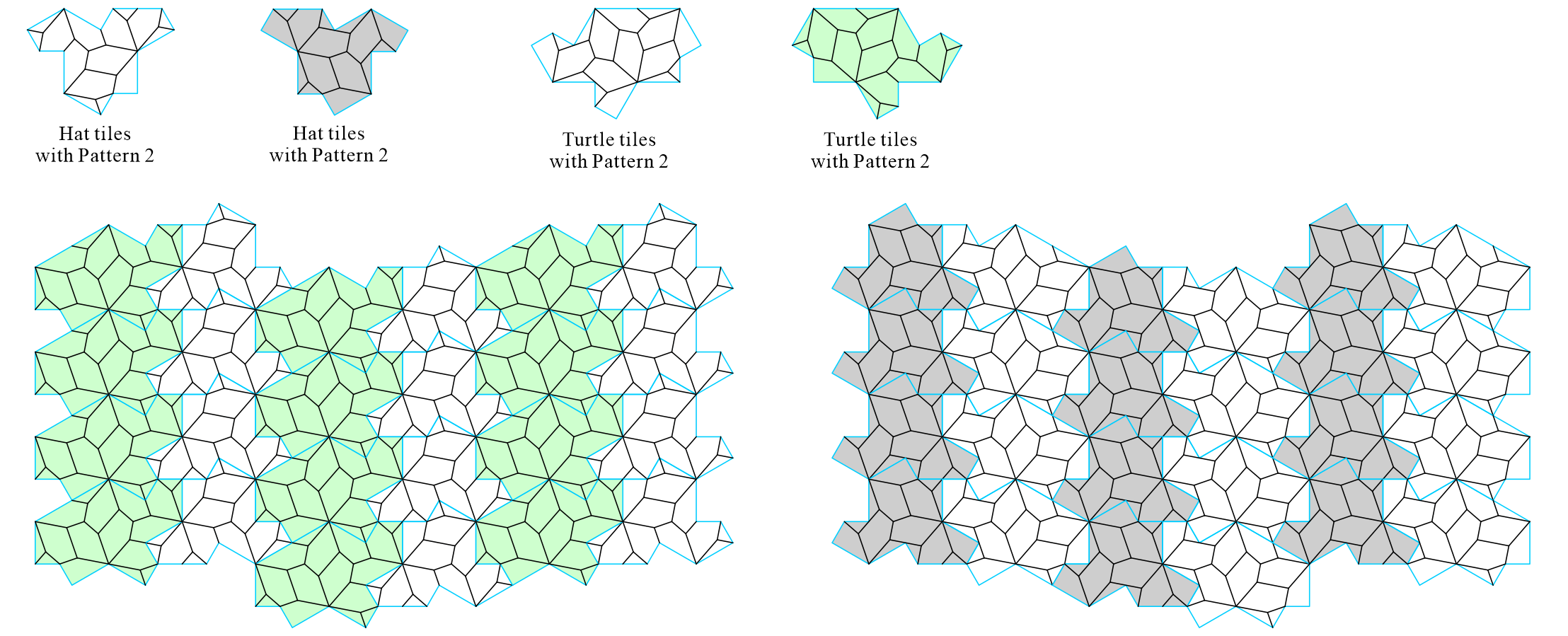} 
  \caption{{\small Periodic tiling in Figure~\ref{Fig.4-1} reconstructed using hat tiles and turtle tiles with Pattern 2.} }
\label{Fig.B-2}
\end{figure}

Figure~\ref{Fig.B-1} shows the periodic tiling in Figure~\ref{Fig.4-1} reconstructed using hat 
tiles and turtle tiles with Pattern 1. Figure~\ref{Fig.B-2} shows the periodic tiling in 
Figure~\ref{Fig.4-1} reconstructed using hat tiles and turtle tiles with Pattern 2.

\begin{figure}[H]
 \centering\includegraphics[width=14.5cm,clip]{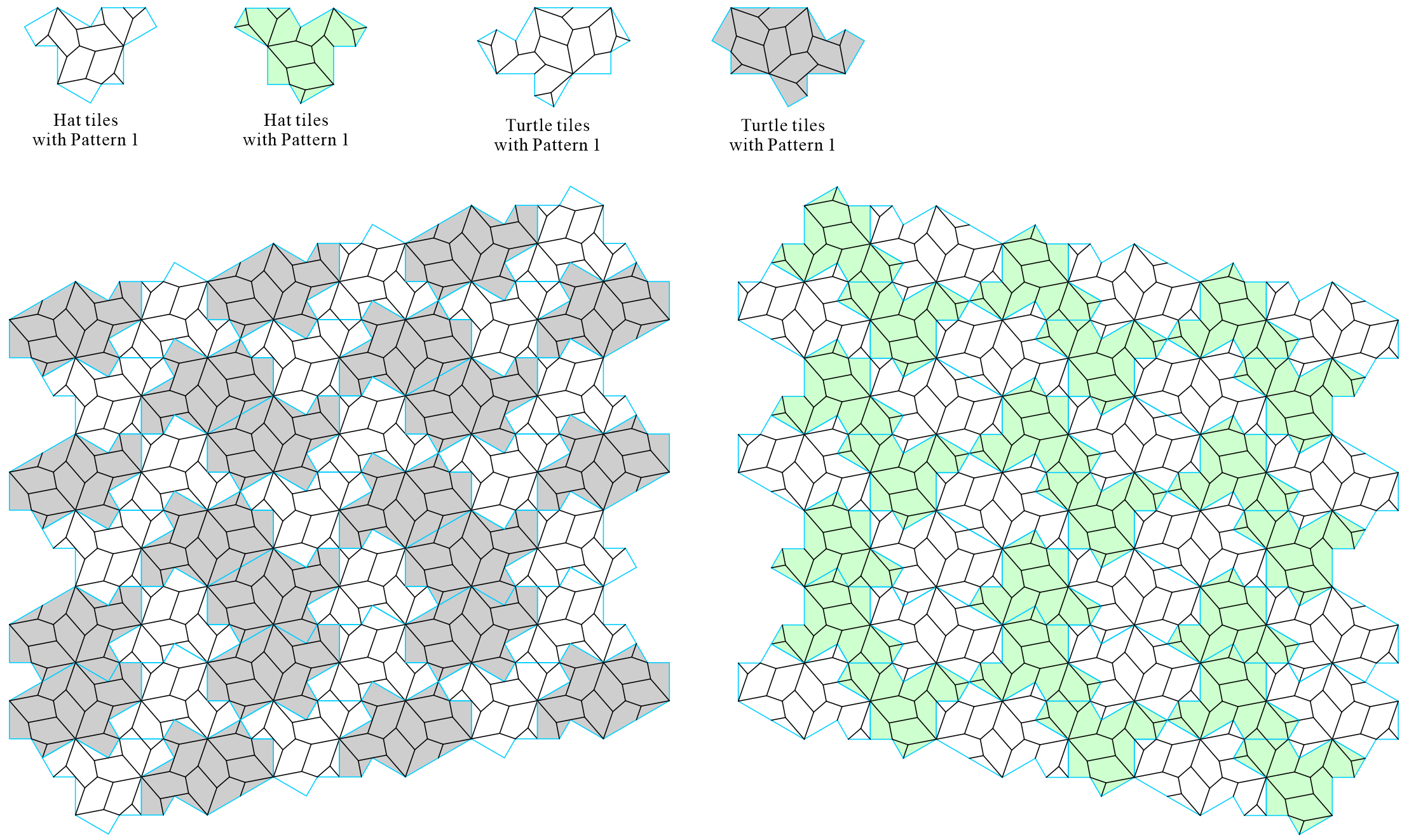} 
  \caption{{\small Tilings with the properties described in Figure~\ref{Fig.A-1}, using hat tiles and turtle tiles with Pattern 1.} }
\label{Fig.B-3}
\end{figure}

\begin{figure}[H]
 \centering\includegraphics[width=14.5cm,clip]{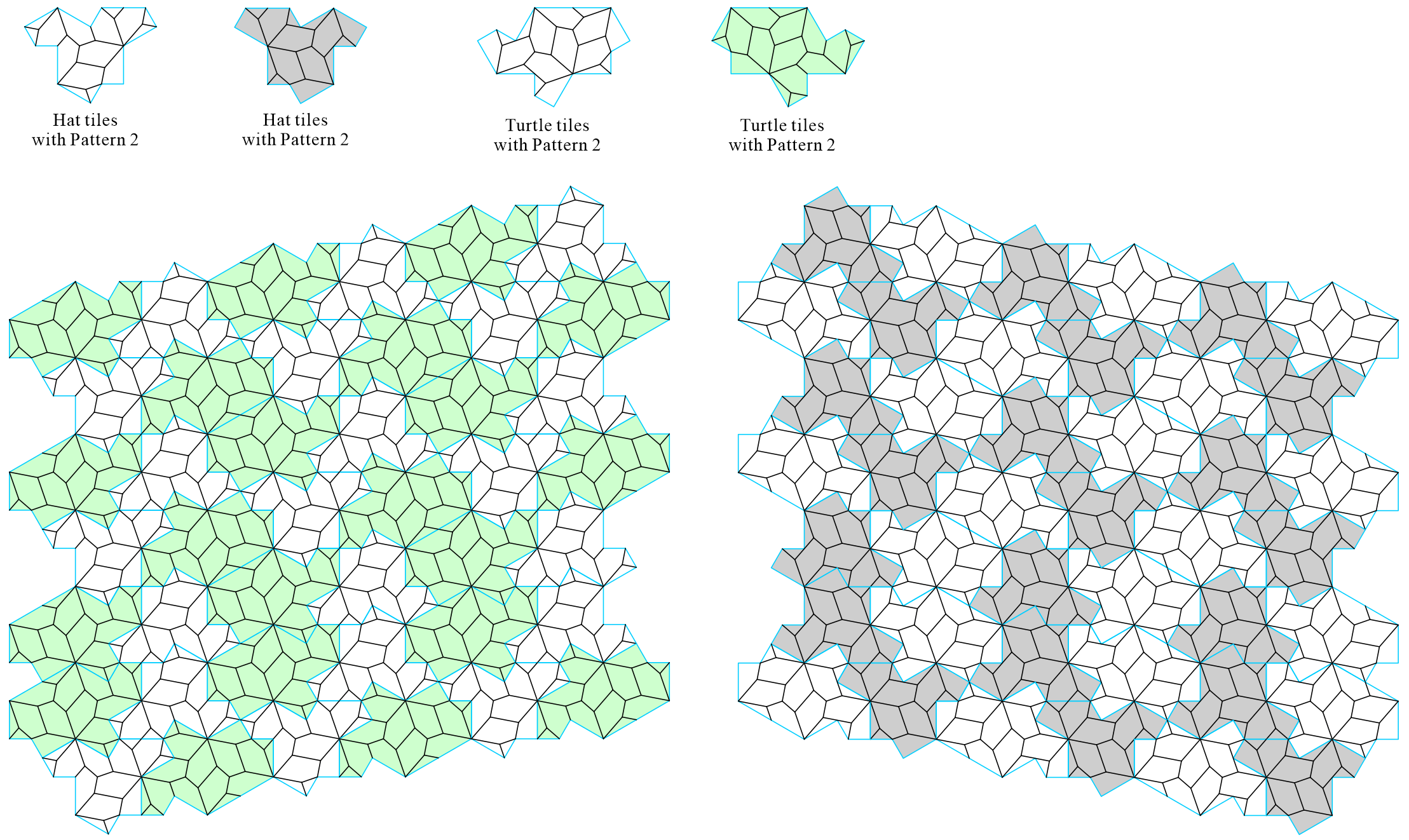} 
  \caption{{\small Tilings with the properties described in Figure~\ref{Fig.A-1}, using hat tiles and turtle tiles with Pattern 2.} }
\label{Fig.B-4}
\end{figure}

\vspace{2\baselineskip}
Figure~\ref{Fig.B-3} shows the tiling with the properties described in Figure~\ref{Fig.A-1}, reconstructed 
using hat tiles and turtle tiles with Pattern 1. Figure~\ref{Fig.B-4} shows the tiling with the properties 
described in Figure~\ref{Fig.A-1}, reconstructed using hat tiles and turtle tiles with Pattern 2.

\vspace{2\baselineskip}

\begin{figure}[H]
 \centering\includegraphics[width=15cm,clip]{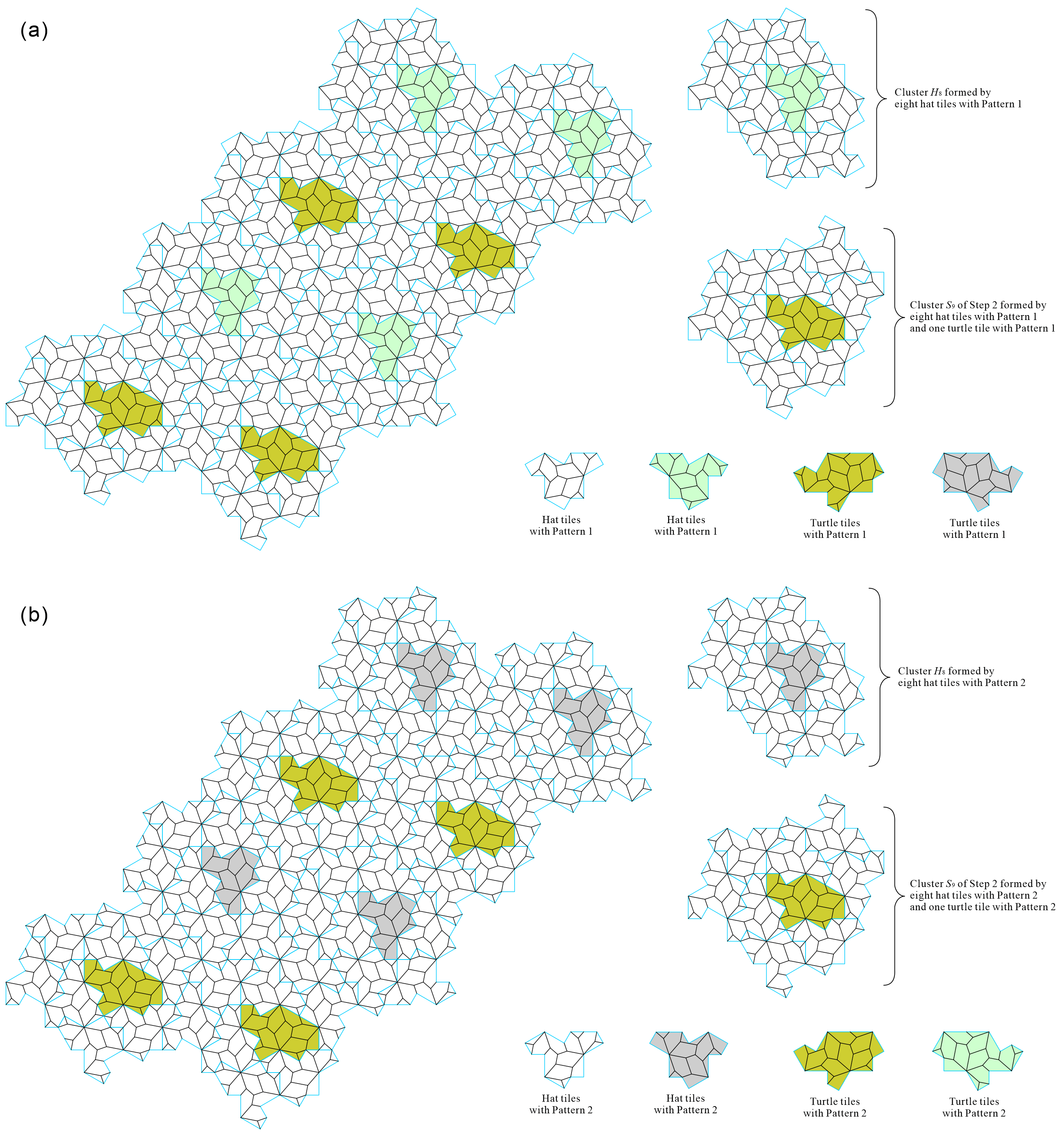} 
  \caption{{\small Periodic tilings in Figure~\ref{Fig.4-2} reconstructed using hat tiles and turtle tiles 
  with Pattern 1 and with Pattern 2.} }
\label{Fig.B-5}
\end{figure}

Figure~\ref{Fig.B-5}(a) shows the periodic tiling in Figure~\ref{Fig.4-2} reconstructed using hat 
tiles and turtle tiles with Pattern 1. Figure~\ref{Fig.B-5}(b) shows the periodic tiling in 
Figure~\ref{Fig.4-2} reconstructed using hat tiles and turtle tiles with Pattern 2.

\vspace{6\baselineskip}

\begin{figure}[H]
 \centering\includegraphics[width=15cm,clip]{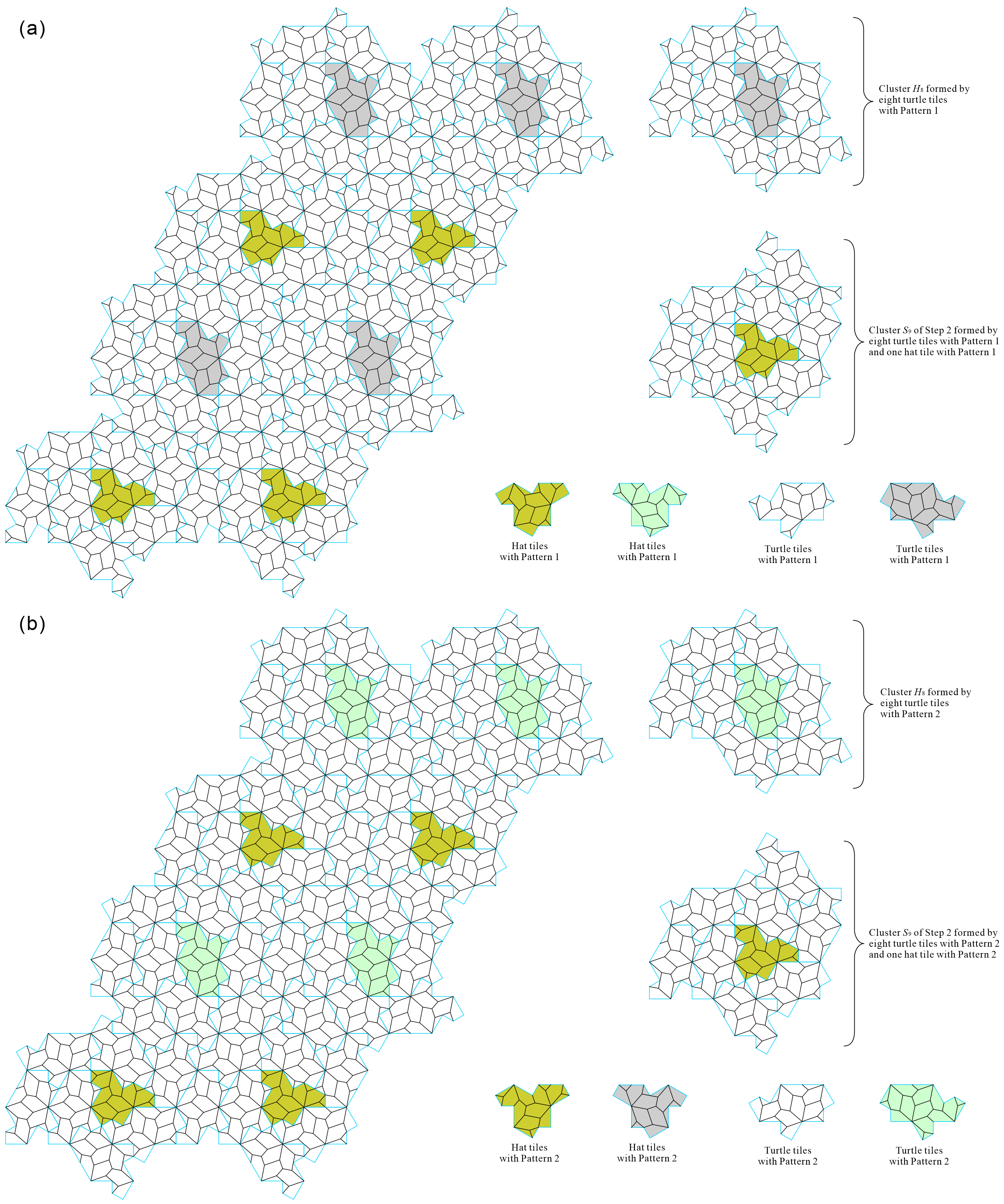} 
  \caption{{\small Periodic tilings in Figure~\ref{Fig.4-3} reconstructed using hat tiles and turtle tiles 
  with Pattern 1 and with Pattern 2.} }
\label{Fig.B-6}
\end{figure}

Figure~\ref{Fig.B-6}(a) shows the periodic tiling in Figure~\ref{Fig.4-3} reconstructed using hat 
tiles and turtle tiles with Pattern 1. Figure~\ref{Fig.B-6}(b) shows the periodic tiling in 
Figure~\ref{Fig.4-3} reconstructed using hat tiles and turtle tiles with Pattern 2.

\vspace{2\baselineskip}

\begin{figure}[H]
 \centering\includegraphics[width=11.5cm,clip]{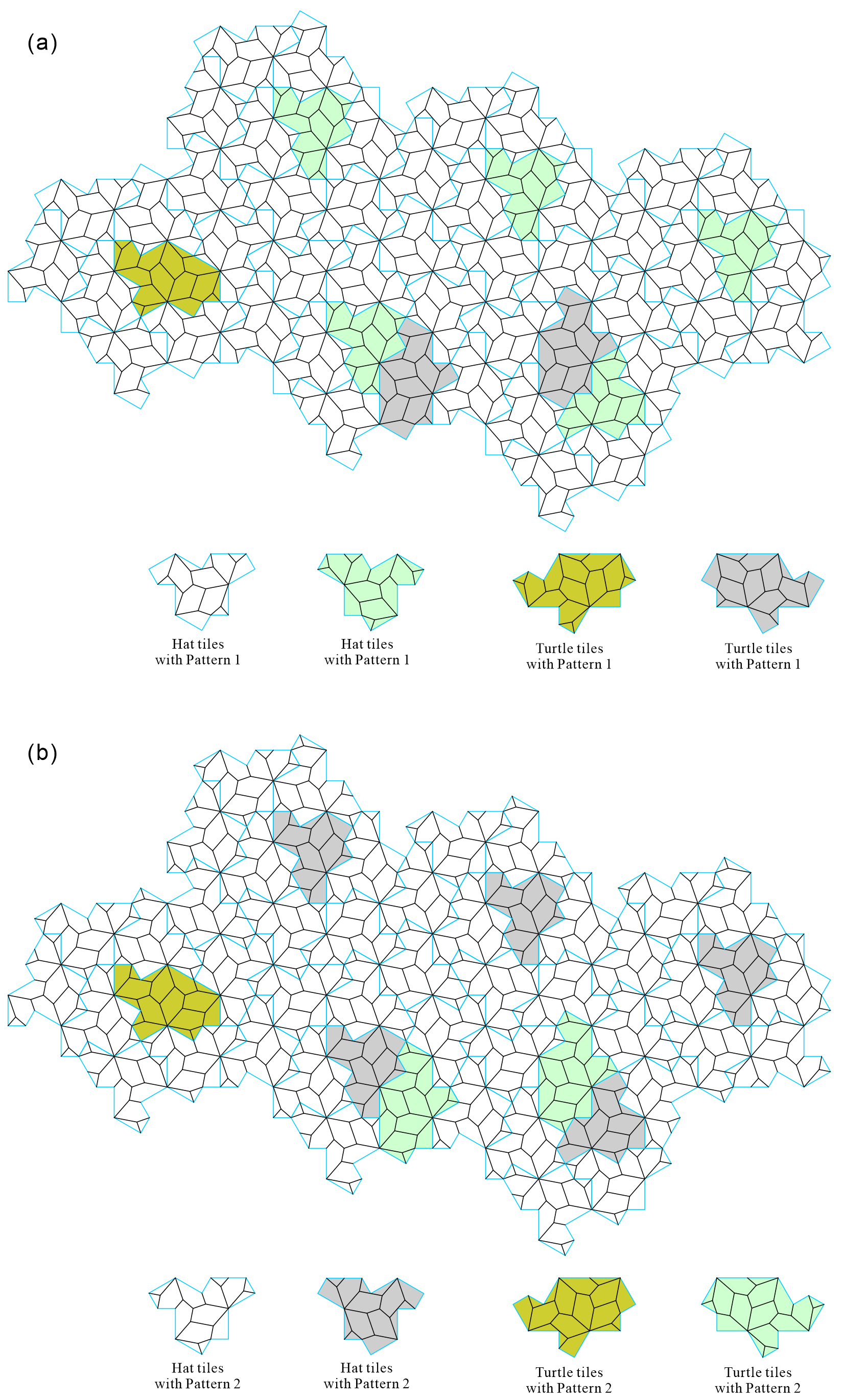} 
  \caption{{\small Translation units of the periodic tiling in Figure~\ref{Fig.4-7} reconstructed 
  using hat tiles and turtle tiles with Pattern 1 and with Pattern 2.} }
\label{Fig.B-7}
\end{figure}

Figure~\ref{Fig.B-7}(a) shows the translation unit of the periodic tiling in Figure~\ref{Fig.4-7} 
reconstructed using hat tiles and turtle tiles with Pattern 1. Figure~\ref{Fig.B-7}(b) shows the 
translation unit of the periodic tiling in Figure~\ref{Fig.4-7} reconstructed using hat tiles and 
turtle tiles with Pattern 2.

\begin{figure}[H]
 \centering\includegraphics[width=15cm,clip]{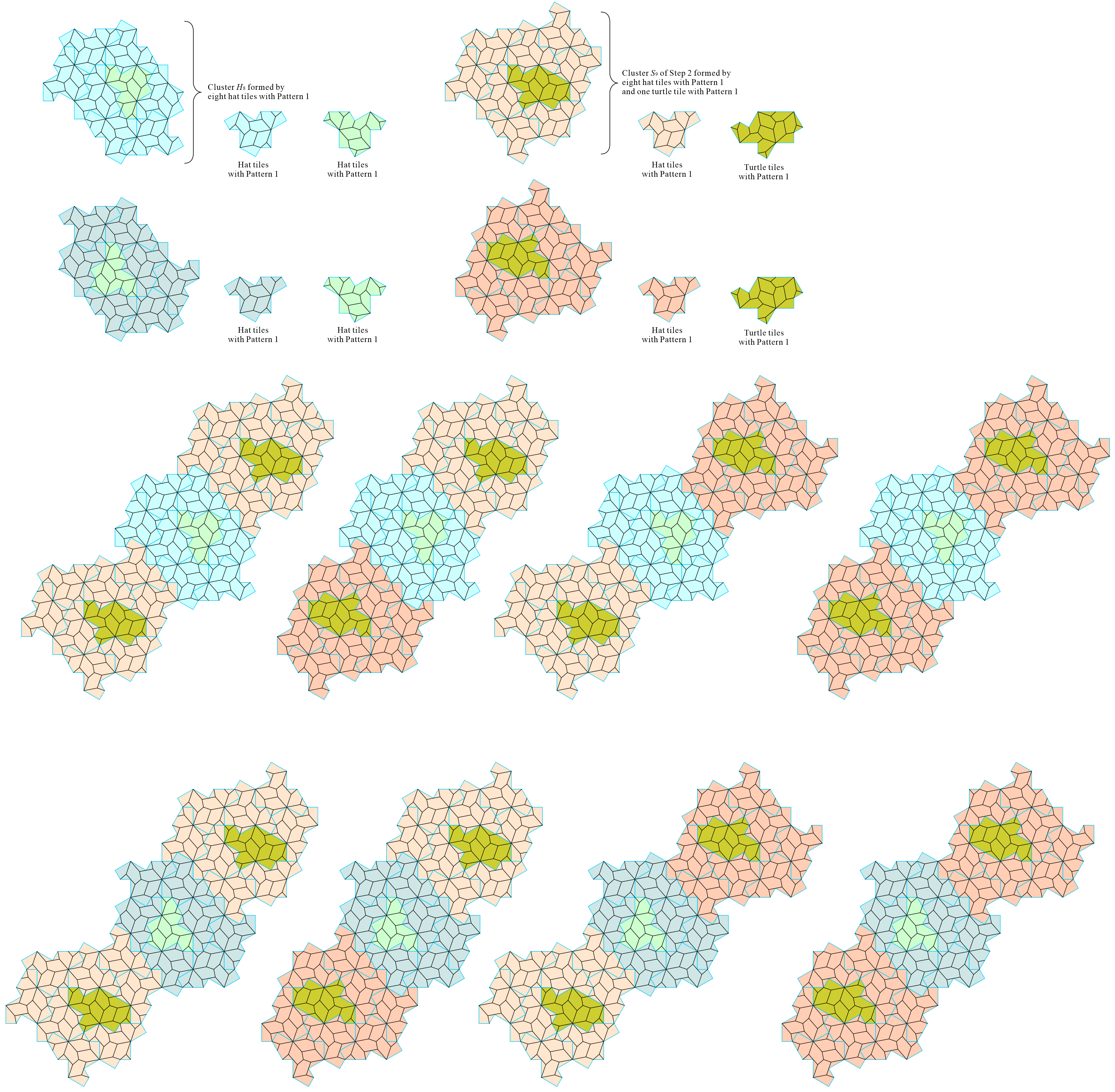} 
  \caption{{\small Variations in the combinations of $H_8$ and $S_9$ in Figure~\ref{Fig.A-3} drawn 
  using hat tiles and turtle tiles with Pattern 1.} }
\label{Fig.B-8}
\end{figure}

\begin{figure}[H]
 \centering\includegraphics[width=15cm,clip]{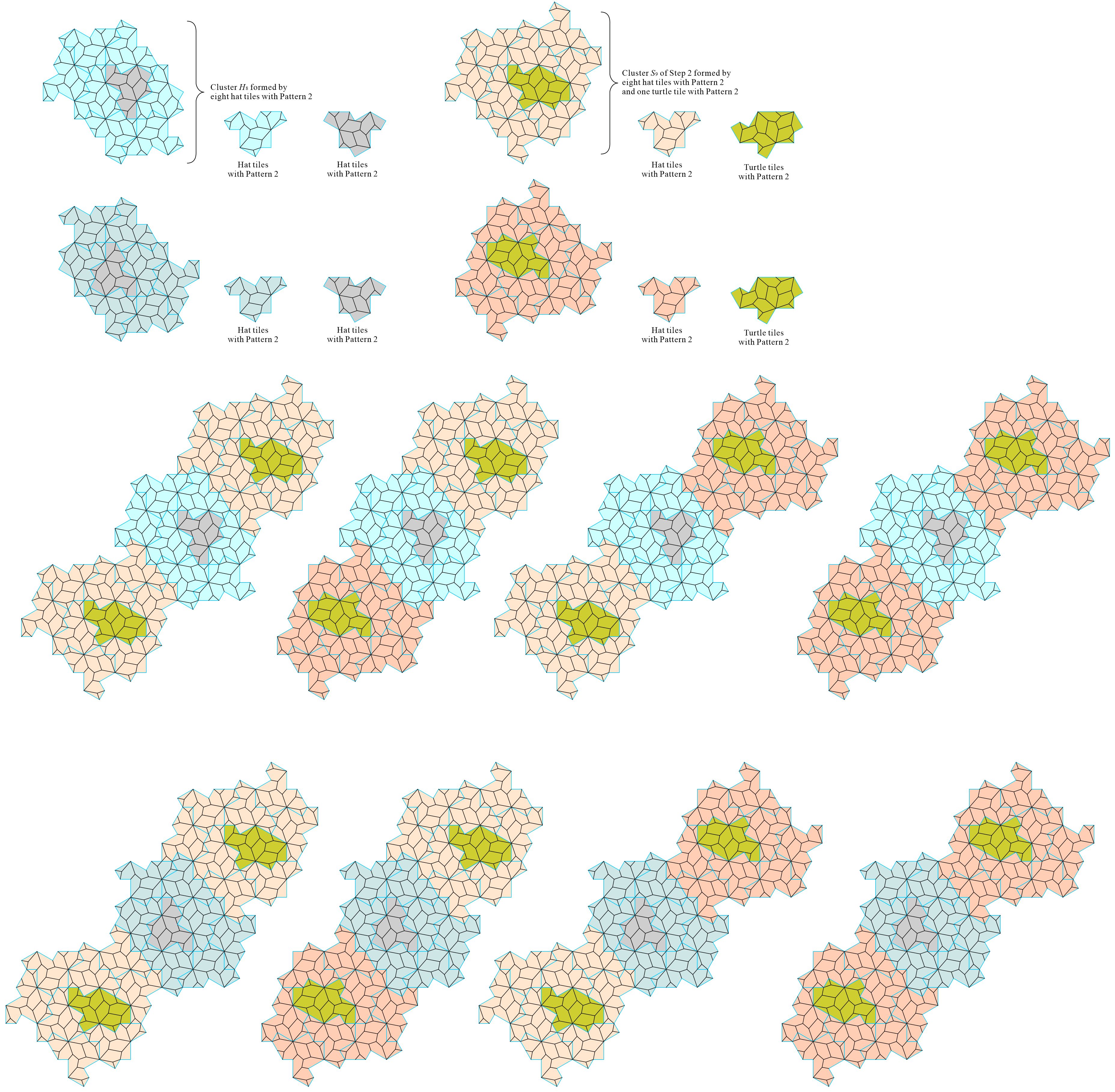} 
  \caption{{\small Variations in the combinations of $H_8$ and $S_9$ in Figure~\ref{Fig.A-3} drawn 
  using hat tiles and turtle tiles with Pattern 2.} }
\label{Fig.B-9}
\end{figure}

Figure~\ref{Fig.B-8} shows the variations in the combinations of $H_8$ and $S_9$ in Figure~\ref{Fig.A-3} 
drawn using hat tiles and turtle tiles with Pattern 1. Figure~\ref{Fig.B-9} shows the variations in 
the combinations of $H_8$ and $S_9$ in Figure~\ref{Fig.A-3} drawn using hat tiles and turtle tiles 
with Pattern 2.

\vspace{8\baselineskip}
\section{An example in which one of the concave polygons admits periodic tilings}
\label{appC}

\begin{figure}[H]
 \centering\includegraphics[width=10.5cm,clip]{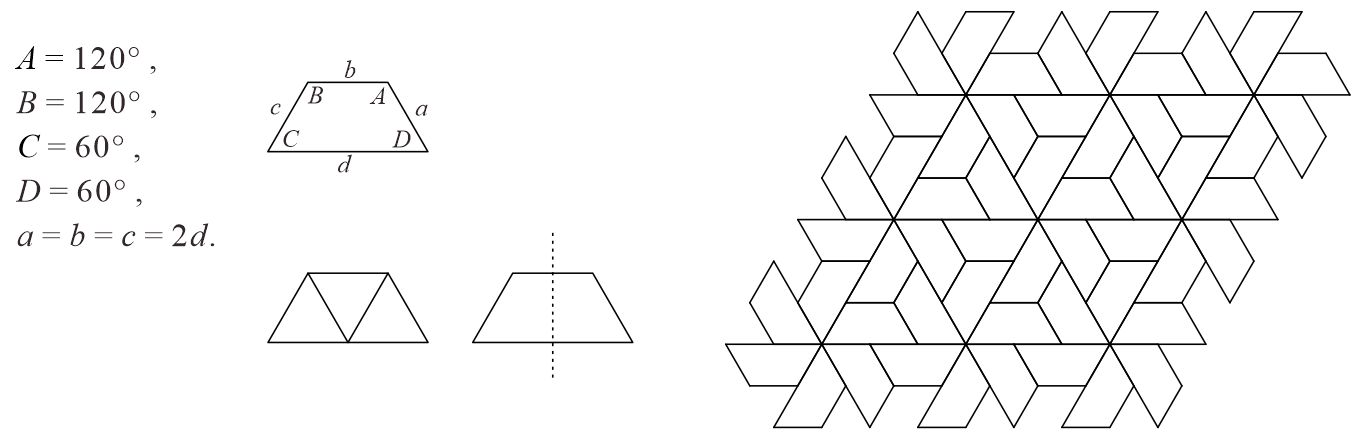} 
  \caption{{\small Periodic tiling with valence-six vertices, formed using a trapezoidal monotile 
  that has line symmetry and can be divided into three congruent equilateral triangles.} }
\label{Fig.C-1}
\end{figure}

\begin{figure}[H]
 \centering\includegraphics[width=14.5cm,clip]{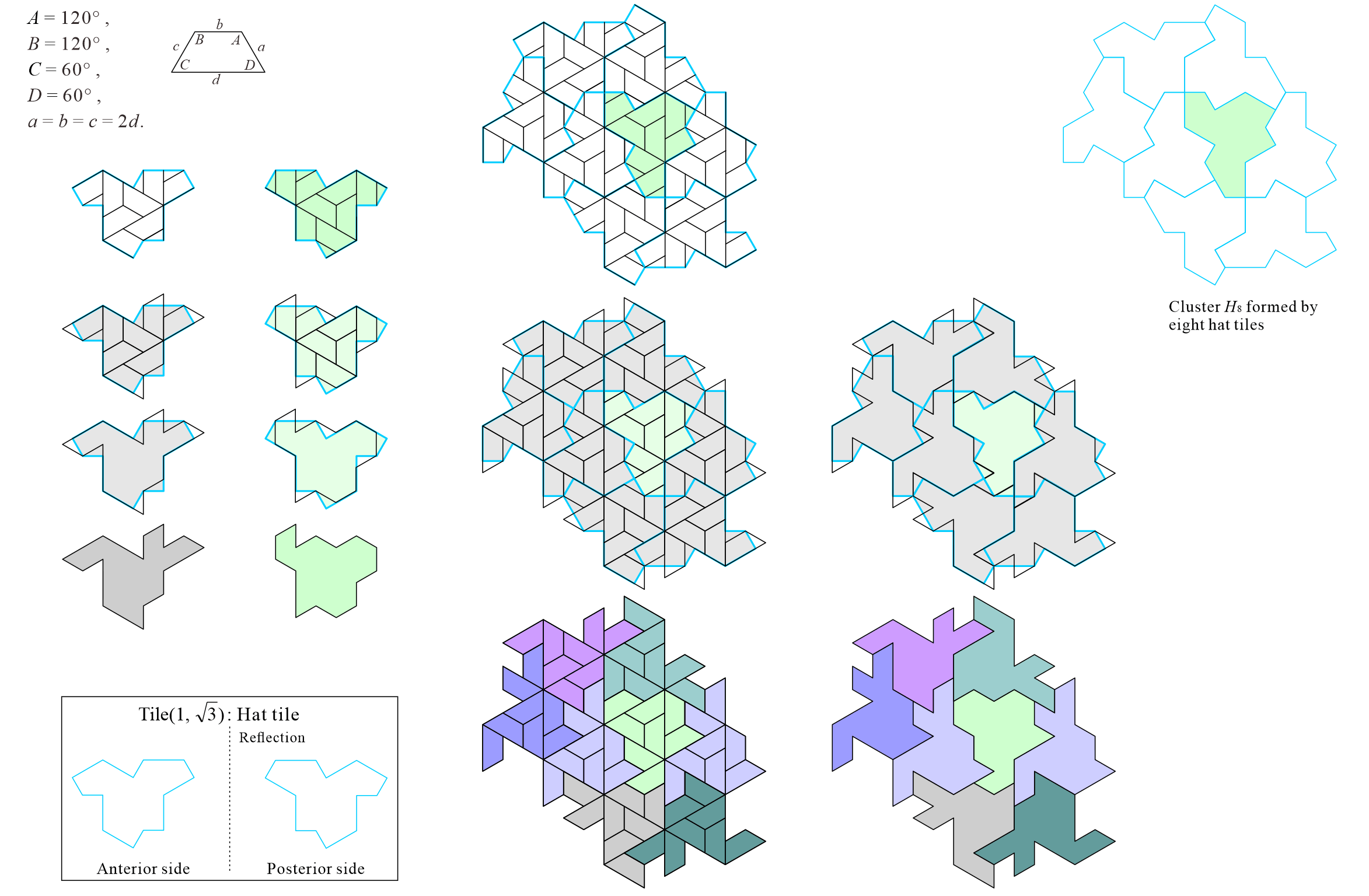} 
  \caption{{\small Periodic tiling obtained by applying Sugimoto's Perspective and 
  Amfirifma's Perspective to the periodic tiling in Figure~\ref{Fig.C-1} formed using the trapezoid, 
  thereby corresponding it to the hat tile.} }
\label{Fig.C-2}
\end{figure}

\begin{figure}[H]
 \centering\includegraphics[width=12cm,clip]{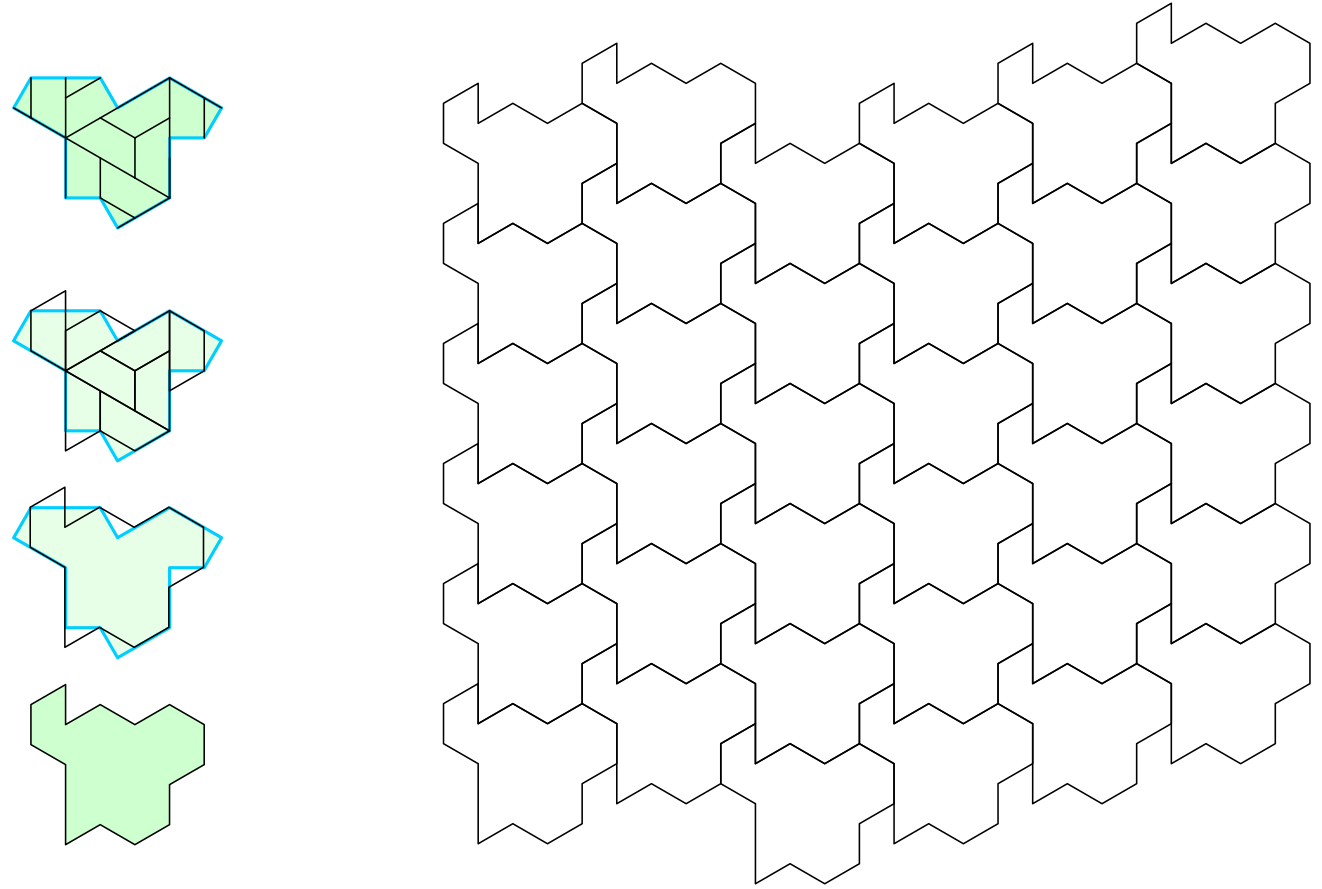} 
  \caption{{\small Periodic tiling generated by one of the concave polygons shown in Figure~\ref{Fig.C-2}.} }
\label{Fig.C-3}
\end{figure}

Figure~\ref{Fig.C-1} shows a periodic tiling with valence-six vertices, formed using a trapezoidal 
monotile that has line symmetry and can be divided into three congruent equilateral triangles. 
Figure~\ref{Fig.C-2} shows an example obtained by associating the tiling formed by this trapezoid 
with hat tiles according to Sugimoto's Perspective and Amfirifma's Perspective introduced in 
Section~\ref{section1}. As shown in the figure, the two concave polygons obtained in 
Figure~\ref{Fig.C-2} can form cluster $H_8$ and therefore can generate non-periodic tilings by 
repeated $H_7/H_8$ substitution. On the other hand, as shown in Figure~\ref{Fig.C-3}, one of 
these concave polygons admits periodic tilings.

From Figures~\ref{Fig.C-1}--\ref{Fig.C-3}, it can be confirmed that the concave polygons 
derived according to Sugimoto's Perspective and Amfirifma's Perspective do not necessarily 
have finite Heesch numbers.

\end{document}